\documentclass[a4paper,11pt]{article}
\usepackage{amssymb}
\usepackage{amsmath}
\usepackage{epsfig}
\usepackage{makeidx}
\usepackage[all]{xy}

\makeindex

\newtheorem{theorem}{Theorem}
\newtheorem{definition}[theorem]{Definition}
\newtheorem{conjecture}[theorem]{Conjecture}
\newtheorem{question}[theorem]{Question}

\newtheorem{proposition}[theorem]{Proposition}
\newtheorem{lemma}[theorem]{Lemma}
\newtheorem{corollary}[theorem]{Corollary}
\newtheorem{theorem-question}[theorem]{Theorem-Question}

\newenvironment{exafont}{\begin{bf}}{\end{bf}}
\newenvironment{example}{\vspace{0.3cm}\par\noindent\refstepcounter{theorem}\begin{exafont}Example
\thetheorem\end{exafont}\hspace{\labelsep}}{\vspace{0.3cm}\par}

\newenvironment{note}{\vspace{0.3cm}\par\noindent\refstepcounter{theorem}\begin{exafont}Note \thetheorem\end{exafont}\hspace{\labelsep}}{\vspace{0.3cm}\par}

\newenvironment{remark}{\vspace{0.3cm}\par\noindent\refstepcounter{theorem}\begin{exafont}Remark \thetheorem\end{exafont}\hspace{\labelsep}}{\vspace{0.3cm}\par}

\begin{document}

\baselineskip = 15pt

\title{\large Rock blocks.} 

\author{W. Turner}

\maketitle

\begin{center}
{\small \it Dedicated to the memory of Joe Silk.}
\end{center}

\bigskip

\begin{abstract}
Consider representation theory associated to 
symmetric groups, or to Hecke algebras in type A, 
or to $q$-Schur algebras, 
or to finite general linear groups
in non-describing characteristic.
Rock blocks are certain combinatorially defined
blocks appearing in such a representation theory, 
first observed by R. Rouquier.  
Rock blocks are much more symmetric than general blocks, 
and every block is derived equivalent to a Rock block.
Motivated by a theorem of J. Chuang and R. Kessar in the case of symmetric
group blocks of abelian defect, we pursue a structure 
theorem for these blocks.
\end{abstract}

\newpage

\begin{center}
{\bf \large 
Contents.}
\end{center}

{Introduction.}

\begin{description}

{\item I.
Highest weight categories, $q$-Schur algebras,
Hecke algebras, and finite general linear groups.}

{ 
\item II.
Blocks of $q$-Schur algebras,
Hecke algebras, and finite general linear groups.}

{ 
\item III.
Rock blocks of finite general linear groups and Hecke algebras,
when $w<l$.}

{ 
\item IV.
Rock blocks of symmetric groups, and the Brauer morphism.}

{ 
\item V.
Schur-Weyl duality inside Rock blocks of symmetric groups.}

{ 
\item VI.
Ringel duality inside Rock blocks of symmetric groups.}

{ 
\item VII.
James adjustment algebras for Rock blocks of symmetric groups.}

{ \item VIII.
Doubles, Schur super-bialgebras, and Rock blocks of Hecke algebras.}

{\item IX.
Power sums.}

{\item X.
Schiver doubles of type $A_\infty$.}

\end{description}

{References.}

\newpage

\begin{center}
{\bf \large 
Introduction}
\end{center}

\bigskip

Whilst the 20$^{th}$ century was still in its infancy, 
an article by F. G. Frobenius was published in
the Journal of the Berlin Science Academy,
which contained a description of 
the irreducible complex characters of all symmetric groups \cite{Frobenius}.
Since then, representation theory has evolved into a deep and sophisticated art, 
to the point where most papers in the subject are incomprehensible to 
the multitude of mathematicians.
After all this development however, some basic questions remain unanswered. 
Interrogate an expert on group representation theory over finite fields and 
you will quite soon witness a shrug of the shoulders, 
and a protestation of ignorance. 
The irreducible characters of symmetric groups, which Frobenius so casually exposed in characteristic zero, 
remain mysterious over fields of prime characteristic.
This monograph comprises a sequence of reflections 
surrounding the modular representation theory of symmetric groups.

Our approach to the subject is homological, inspired by M. Brou\'e's abelian defect group conjecture \cite{Broue}, and
encouraged by the proof of Brou\'e's conjecture for blocks of symmetric groups by 
J. Chuang, R. Kessar, and R. Rouquier.

The abelian defect group conjecture is the most homological of a menagerie of general conjectures in modular representation theory, 
each of which predicts a likeness between the representations 
of a finite group in characteristic $l$, and 
those of its $l$-local subgroups. 
It stakes that the derived category of any block $A$ of a finite group is equivalent to the derived category of its 
Brauer correspondent $B$, so long as the blocks have abelian defect groups. 
Such an equivalence should respect the triangulated structure of the derived category, 
and therefore descend from a two sided tilting complex of $A$-$B$-bimodules, by a theorem of J. Rickard \cite{Rickard}. 
Some have postured to prove the conjecture by induction, 
and encountered the difficulty of lifting an equivalence of stable categories to an equivalence of derived categories \cite{stable}. 
Others have tried to prove the conjecture for particular examples, such as symmetric groups. 

In 1991, R. Rouquier observed a certain class of blocks of symmetric groups, 
which he believed to possess a particularly simple structure. Indeed, Rouquier 
conjectured a beautiful structure theorem for such blocks 
of abelian defect, which was 
subsequently proved by J. Chuang, and R. Kessar \cite{KessCh}.
A corollary was a proof of Brou\'e's conjecture for this class of blocks, 
whose defect groups could be arbitrarily large.
In view of their history, these blocks should properly be called Rouquier, or 
Chuang-Kessar blocks. We use the curt abbreviation ``RoCK blocks''.
Such blocks can be defined in arbitrary defect, 
and in any species of type $A$ representation theory.

Chuang and Rouquier proved in a later work that all symmetric group blocks of identical defect 
possess equivalent derived categories which, in conjunction with the previous study of Rock blocks, 
established the truth of Brou\'e's conjecture for all blocks of symmetric groups 
\cite{ChRou}.

The will which motivated this text, was for a theorem like that of Chuang and
Kessar, describing Rock blocks whose defect groups 
are not necessarily abelian.
Such a result ought to be of broad interest since there is no known analogue of 
Brou\'e's conjecture in nonabelian defect. 

So far as symmetric groups are concerned, 
an algorithm of A. Lascoux, B. Leclerc and J-Y. Thibon gives a conjectural
description of all the decomposition numbers of blocks of abelian defect \cite{LLT}.
This algorithm was formulated upon the examination of numerous tables of decomposition numbers for small symmetric groups,
made by G. James \cite{JamKerb}. 
The algorithm has been proven to describe the decomposition numbers for Hecke algebras at a root of unity in characteristic zero, by S. Ariki \cite{Ariki}.   
In nonabelian defect however, little is currently known about the decomposition numbers, 
even conjecturally.
Symmetric group blocks of nonabelian defect are therefore beyond the influence of the general 
character theoretic predictions for algebraic groups made by G. Lusztig.

Over the length of this monograph, I hope to convince the Reader of the existence of 
a structure theorem for Rock blocks of arbitrary defect.
Indeed, a conjectural description of an arbitrary Rock block of a Hecke algebra  
appears in chapter 8 (conjecture~\ref{Rock}).
A generalized conjecture for blocks
of $q$-Schur algebras is given in chapter 9 (conjecture~\ref{SchurRock}).

In formulating these conjectures, I produced 
more than a cupful of imaginative sweat, and earlier
chapters of this booklet record theorems which 
point to the conjectures, and
give evidence for them.
For those readers with a fetish for decomposition numbers,
the saltiest of these theorems is probably a formula for the
decomposition matrix of a Rock block of a symmetric group, 
of arbitrary defect (theorem~\ref{nilpotent}).
 
\bigskip

It would be polite of me to be a little more precise.
Therefore, 
let us catalogue the more significant results of the article, 
and give some description of their 
character and logical intimacy,
before plunging into the depths of the text.

\bigskip

In the first chapter, we recall E. Cline, B. Parshall, 
and L. Scott's definition of a highest weight category.
We discuss quasi-hereditary algebras, 
and summarise C. Ringel's tilting theory for these algebras.
We recognise the $q$-Schur algebra ${\cal S}_q(n)$ as the graded dual 
of the quantized coordinate ring of a matrix algebra,
and recall S. Donkin's tilting theory for ${\cal S}_q(n)$. 
Schur-Weyl duality relates $q$-Schur algebras, and Hecke algebras ${\cal H}_q(\Sigma_r)$ 
associated to symmetric groups $\Sigma_r$.
We view this phenomenon in chapter one, and its relevance for the representation theory
of finite general linear groups $GL_n(\mathbb{F}_q)$ over a field $k$ of characteristic $l$, 
coprime to $q$.
Parametrizations of irreducible representations for all these algebras are given.
We assemble a few facts concerning wreath products of algebras.

\bigskip

Chapter two opens with a description of the abacus presentation of partitions, 
due to G. James.
We recall Nakayama's parametrization of blocks $k{\bf B}_{\tau,w} ^\Sigma$ of symmetric groups by their 
$l$-core $\tau$ and their weight $w$, 
as well as parametrizations
of blocks $k{\bf B}_{\tau,w} ^{{\cal S}_q}$ of $q$-Schur algebras,
blocks $k{\bf B}_{\tau,w} ^{{\cal H}_q}$ of Hecke algebras, 
and unipotent blocks $k{\bf B}_{\tau,w} ^{G_q}$ of finite general linear groups.
Blocks $k{\bf B}_{\tau,w} ^\Sigma$ have abelian defect groups if, and only if, $w<l$.
Chuang and Rouquier's general theory of $\mathfrak{sl}_2$-categorification implies that 
$k{\bf B}_{\tau,w} ^{\cal X}$ is derived equivalent to 
$k{\bf B}_{\tau',w} ^{\cal X}$, for a fixed weight $w$, fixed 
${\cal X} \in \{ \Sigma, {\cal S}_q, {\cal H}_q, G_q \}$, and various $\tau, \tau'$.

Blocks of symmetric groups are in one-one correspondence 
with weight spaces in the basic representation of 
$\widehat{\mathfrak{sl}_p}$.
Rock blocks are a particularly symmetric class of blocks,
distinguished combinatorially via their abacus presentation. 
When $w<l$, Chuang and Kessar's structure theorem states that a Rock block 
$k{\bf B}_{\rho,w} ^\Sigma$ of weight $w$ is
Morita equivalent to a wreath product $k{\bf B}_{\emptyset,1} ^\Sigma \wr \Sigma_w$ 
of a cyclic defect block
$k{\bf B}_{\emptyset,1} ^\Sigma$, and a symmetric group $\Sigma_w$ on $w$ letters.
A stated corollary of this theorem is a formula for the decomposition matrix of a 
symmetric group Rock block of abelian defect 
discovered by Chuang and K.M. Tan, and independently by B. Leclerc and H. Miyachi:
$$d_{\underline{\lambda} \underline{\mu}} =
\sum_{\alpha \in \Lambda^{p+1}  _w, \beta \in \Lambda^p _w}
\prod_{j=0} ^{p-1} c(\lambda^j; \alpha^j, \beta^j) 
c(\mu^j; \beta^j, (\alpha^{j+1})').$$
In this formula, 
$\underline{\lambda} = [\lambda^0,...,\lambda^{l-1}]$ is the $p$-quotient
of a partition with $p$-core $\rho$, relative to a certain abacus presentation;
the numbers $c(\lambda; \mu, \nu)$ are Littlewood-Richardson coefficients.

We state R. Paget's description of the Mullineux map on 
a Rock block of a Hecke algebra.
We proceed to describe Brauer correspondence for blocks of finite general linear groups, 
and then for blocks of symmetric groups. 
This chapter ends with the statement of a criterion of Brou\'e, for the lifting of a 
character correspondence between symmetric algebras to a Morita equivalence. 

\bigskip

The first honest mathematics appears in chapter three.
We sketch a proof, \`{a} la Chuang-Kessar, of a structure theorem for Rock blocks of
finite general linear groups of abelian defect.
More precisely, we show that a Rock block $k{\bf B}_{\rho,w}^{G_q}$ 
is Morita equivalent to the wreath product $k{\bf B}_{\emptyset,1} ^{G_q} \wr \Sigma_w$, when $w<l$
(theorem~\ref{main}).
A simple corollary is the Morita equivalence of $k{\bf B}_{\rho,w}^{{\cal H}_q}$ 
and $k{\bf B}_{\emptyset,1} ^{{\cal H}_q} \wr \Sigma_w$, 
so long as $w<l$, and $q \neq 1$ modulo $l$ (theorem~\ref{mainnew}). 

This structure theorem for Rock blocks of finite general linear groups of abelian defect
appeared in a previous paper of mine, 
and was written down independently by Miyachi \cite{Miyachi}. 
An application given in my paper is the revelation of Morita equivalences
between weight two blocks of finite general linear groups $GL_n(\mathbb{F}_q)$, as $q$ varies. 
Thanks to Chuang and Rouquier's theory it must now be possible to generalise this result,
and give comparisons between abelian defect blocks of $GL_n(\mathbb{F}_q)$ of arbitrary weight, 
as $q$ varies. We choose not to spend time chomping on this old pie, 
since we have become aware of dishes with a more exotic, and alluring aroma. 

\bigskip

In chapter four, we turn to Rock blocks of symmetric groups, of nonabelian defect.
We give a sweet proof that in characteristic two,
once a Rock block of a symmetric group has been localised at some idempotent,
Chuang and Kessar's theorem generalises to non-abelian defect.
To be precise, we prove that in characteristic two, $e k{\bf B}_{\rho,w} ^\Sigma e$ is Morita equivalent
to $k\Sigma_2 \wr \Sigma_w$, for some idempotent $e \in k{\bf B}_{\rho,w} ^\Sigma$ 
(theorem~\ref{viva}).
Our proof involves the Brauer morphism. 
This idea can be contorted and extended, to give a result in arbitrary
characteristic.
Indeed, to any Rock block $k{\bf B}_{\rho,w}^\Sigma$, 
we associate a natural $l$-permutation module $kM$, 
whose endomorphism ring $k{\cal E}$ is Morita equivalent to $k\Sigma_w$ 
(theorem~\ref{Wick}).

\bigskip

Schur-Weyl duality for $\Sigma_w$ is naturally visible inside a Rock block of weight $w$,
via the Brauer homomorphism, and our fifth chapter is dedicated to establishing this
truth. Formally, we prove that $k{\bf B}_{\rho,w} ^\Sigma / Ann(kM)$ is Morita
equivalent to the Schur algebra $S(w,w)$, and that the $S(w,w)$-$k\Sigma_w$-bimodule
corresponding to the $k{\bf B}_{\rho,w} ^\Sigma$-$k{\cal E}$-bimodule $kM$ is 
twisted tensor space (theorem~\ref{Comeon}).

\bigskip

Chapter six begins with a criterion for the lifting of a character correspondence
between quasi-hereditary algebras to a Ringel duality (theorem~\ref{Brother}), 
echoing Brou\'e's criterion
for a Morita equivalence between symmetric algebras. We use this result to prove the 
existence of Ringel dualities between certain subquotients
$kA_{(0,a_1,...,a_{p-2})}$ and $kB_{(a_1,....,a_{p-2},0)}$ of $k{\bf B}_{\rho,w} ^\Sigma$
(theorem~\ref{Ringel}). 
We call this collection of Ringel dualities a ``walk along the abacus'', 
because it is reminiscent of 
J. A. Green's observations on the homological algebra of the Brauer tree
\cite{Greentree}.

\bigskip

In the seventh chapter, we introduce the James adjustment algebra of a Hecke algebra
block. This is a quotient of the block by a nilpotent ideal, 
whose decomposition matrix is equal to the James 
adjustment matrix of the block. 
The principal result of this chapter is theorem~\ref{nilpotent}, 
which states that the James
adjustment algebra of $k{\bf B}_{\rho,w} ^\Sigma$ is Morita equivalent to a direct sum,
$$\bigoplus_{\substack{a_1,...,a_{p-1} \in \mathbb{Z}_{\geq 0} \\ \sum a_i =w}} 
\left( \bigotimes _{i=1} ^{p-1} {\cal S}(a_i,a_i) \right) ,$$
of tensor products of Schur algebras.
The decomposition matrix of $k{\bf B}_{\rho,w} ^\Sigma$ is a product 
of the matrix of sums of products of Littlewood Richardson coefficients
defined by Chuang $\&$ Tan, 
and Leclerc $\&$ Miyachi, and the decomposition matrix of
the James adjustment algebra.

\bigskip

At the entrance to chapter eight, we define a novel double construction. 
Indeed, given a bialgebra $B$ equipped with an algebra anti-endomorphism $\sigma$, 
which is also a coalgebra anti-endomorphism, and a dual bialgebra $B^*$, 
we project the structure of an associative algebra onto $B \otimes B^*$ 
(theorem~\ref{symmetric}).
Particular examples of these doubles 
show a remarkably close resemblance to Rock blocks.

If $Q$ is a quiver, let $P_Q$ be the path super-algebra of $Q$, 
modulo all quadratic relations.
Let $P_Q(n)$ be the super-algebra Morita equivalent to $P_Q$, 
all of whose irreducible modules have dimension $n$. The coordinate ring of $P_Q(n)$ is a
super-bialgebra, whose double we denote ${\cal D}_Q(n)$.
We call such an algebra a Schur quiver double, or Schiver double.
We prove that ${\cal D}_Q(n)$ is independent of the orientation of $Q$
(theorem~\ref{invariant}).
Conjecture~\ref{Rock} predicts that 
the Rock block $k{\bf B}_{\rho,w} ^\Sigma$ 
is Morita equivalent to a summand ${\cal D}_{A_{l-1}}(w,w)$
of the Schiver double ${\cal D}_{A_{l-1}}(w)$ associated to a quiver of type $A_{l-1}$. 

The principal obstacle to a proof of this conjecture
via the methods introduced here, 
has been my impotence
in producing a suitable grading on a Rock block. The Schiver doubles 
are naturally $\mathbb{Z}_+ \times \mathbb{Z}_+$-graded algebras,
the degree zero part being isomorphic to a tensor product  
$$\bigotimes_{v \in V(Q)} {\cal S}_v(n)$$ of classical Schur algebras
(remark~\ref{grading}).
Such gradings on the Rock blocks remain elusive.
Towards the end of chapter 8, 
we sketch a proof that the graded ring associated to 
a certain filtration on a Rock block,
resembles an algebra summand of a Schiver double. 

If the gradings proposed here on the Rock blocks do exist, then they will pass non-canonically
via derived equivalences to $\mathbb{Z} \times \mathbb{Z}$-gradings on arbitrary blocks \cite{grading}.   
Is it true that every symmetric group block can be positively graded,
so that the degree zero part is isomorphic to the 
James adjustment algebra of the block ?

\bigskip

In the ninth chapter we continue the study of Schiver doubles. 
Given any vertex $v$ of the quiver $Q$ with an arrow $a$ emanating from it, we
define a non-trivial complex ${\cal P}_r(a)$ for ${\cal D}_Q(n,r)$ whose homology groups are all
${\cal S}_v(n,r)$ modules, and
whose character is the power sum $p_r$ (theorem~\ref{pop}).

\bigskip

In the tenth and final chapter of this article, we consider Schiver doubles associated to quivers of type $A_\infty$, which 
enjoy a number of special homological properties. So long as $n \geq r$,
the module category of 
${\cal D}_{A_\infty}(n,r)$ is a highest weight category, and Ringel self-dual 
(theorem~\ref{exch}).
We speculate that any Rock block $k{\bf B}_{\rho,w} ^{{\cal S}_q}$ of a $q$-Schur algebra is 
Morita equivalent to a certain subquotient of ${\cal D}_{A_\infty}(w,w)$ 
(conjecture~\ref{SchurRock}).
We prove the existence of a long exact sequence of ${\cal D}_{A_\infty}(n,r)$-modules
$$......... \rightarrow {\cal D}_{A_\infty}(n,r)
\rightarrow {\cal D}_{A_\infty}(n,r)
\rightarrow {\cal D}_{A_\infty}(n,r) \rightarrow .........,$$ 
which generalises Green's walk along the Brauer tree for an infinite Brauer line. 
(theorem~\ref{walk})

\bigskip

Chapters one and two contain introductory material. Most of it should be familiar to students of
symmetric groups, $q$-Schur algebras, or the like. 
General aspects of finite group representation theory 
such as Brauer correspondence are often omitted 
from presentations of type A representation theory,
but we include a brief account of this correspondence here. 
I consider this to be important philosophically, 
as well as being necessary for some of our proofs. 
This article was conceived before the fire of modular representation theory
laid by R. Brauer, 
and his vision of a sympathy between global and local representations is 
bred in its bones.

The third, fourth, and fifth chapters all make use of local representation theory, 
and should properly be read consecutively. 
The appearance of Schur algebras in chapter five should not be a great 
surprise to students of semisimple algebraic groups familiar with Steinberg's tensor product theorem. 
However, I hope our approach via the Brauer morphism is at least provocative: 
it is so far unclear how to interpret J. Alperin's conjecture homologically in nonabelian defect.

Chapter six can be read independently of chapters three to five, 
and rests on the theory of quasi-hereditary 
algebras. With R. Paget's description of the Mullineux map on a Hecke algebra Rock block 
we cobble a pair of shiny black boots; 
wearing these we are able to comfortably
walk along the abacus. 

The description of the James adjustment algebra of a symmetric group Rock block in chapter seven 
relies on all the theory developed in earlier chapters. 
The results of chapter three allow one to understand some aspects of 
Rock blocks of Hecke algebras 
in characteristic zero, at a root of unity. 
The Schur algebra quotient of chapter five provides information which 
can be carried across the abacus
using the Ringel dualities of chapter six.

Beyond this complex crescendo come chapters eight, nine, and ten. 
The conjectures made here concerning Rock blocks
appear to be quite deep, 
and if proved, would envelop all the results of earlier chapters.  
However, their presentation is logically independent of
chapters three to seven, and carries a lighter burden of notation. 

The development of the article is in the direction of Time's arrow, 
so that more recent ideas appear towards the end of the monograph.

\bigskip

I am most grateful to Joe Chuang, to Karin Erdmann, 
and to Rowena Paget,
for encouraging me amongst these ideas, 
and to Steffen Koenig.
Hannah Turner supported me financially (partly), and libidinously (entirely).
The E.P.S.R.C. gave me some money, as well.
I thank the referee, for his careful reading of the manuscript, and useful comments.

This work, its morality, and the wilful emotions which  
dominated its creation, are dedicated to Joe Silk. 
He was a dear, demonic friend to me, and I wish ... to wish him farewell.

\newpage

\begin{center}
{\bf \large 
Chapter I

Highest weight categories, $q$-Schur algebras, Hecke algebras, and
finite general linear groups.}
\end{center}

\bigskip

We brutally summarise the representation theory of
the $q$-Schur algebra, of Hecke algebras of type A, and
of finite general linear groups in non-describing characteristic.

Although in later chapters, we will invoke such theory 
over more general commutative rings, for simplicity of presentation, 
in this chapter we only consider representation theory over a field $k$, \index{$k$}
of characteristic $l$. \index{$l$}

\begin{center}
{\large Highest weight categories.}
\end{center}

We state some of the principal definitions and results of E. Cline,
B. Parshall, and L. Scott's paper, \cite{CPS}.

\begin{definition} (\cite{CPS}, 3.1)
Let ${\cal C}$ be a locally Artinian, Abelian category over $k$, with enough injectives. 
Let $\Lambda$ be a partially ordered set, such that every interval $[\lambda, \mu]$ is finite, for $\lambda, \mu \in \Lambda$.
The category ${\cal C}$ is a \emph{highest weight category} with respect
to $\Lambda$ if,

\bigskip

(a) $\Lambda$ indexes a complete collection $\{ L(\lambda) \}_{\lambda \in \Lambda}$ 
of non-isomorphic simple objects of ${\cal C}$. \index{$L(\lambda)$}

\bigskip

(b) $\Lambda$ indexes a collection $\{ \nabla(\lambda) \}_{\lambda \in \Lambda}$ of
``costandard objects'' of ${\cal C}$, \index{$\nabla(\lambda)$}
for each of which there exists an embedding
$L(\lambda) \hookrightarrow \nabla(\lambda)$, such that all composition factors
$L(\mu)$ of $\nabla(\lambda)/L(\lambda)$ satisfy $\mu < \lambda$.

\bigskip

(c) For $\lambda, \mu \in \Lambda$, we have,
$dim_k Hom(\nabla(\lambda), \nabla(\mu)) < \infty$, and in addition,
$[ \nabla(\lambda) : L(\mu)] < \infty$.

\bigskip

(d) An injective envelope 
$I(\lambda) \in {\cal C}$ \index{$I(\lambda)$}
of $L(\lambda)$ possesses a filtration
$$0 = F_0(\lambda) \subset F_1(\lambda) \subset ...,$$ such that,

\begin{description}
\item (i) $F_1(\lambda) \cong \nabla(\lambda)$.
\item (ii) For $n>1$, we have $F_n(\lambda)/F_{n-1}(\lambda) \cong \nabla(\mu)$,
for some $\mu = \mu(n) > \lambda$.
\item (iii) For $\mu \in \Lambda$, we have $\mu(n) = \mu$ for finitely many $n$.
\item (iv) $I(\lambda) = \bigcup_i F_i(\lambda)$.
\end{description}
\end{definition}

\begin{definition} (\cite{CPS}, 3.6)
Let $S$ be a finite dimensional algebra over $k$.
Then $S$ is said to be quasi-hereditary if the category $S-mod$, of finitely
generated left $S$-modules is a highest weight category.
\end{definition}

For $M \in {\cal C}$, and $\Gamma \subset \Lambda$, 
let $M_\Gamma$ \index{$M_\Gamma$} be the largest subobject of $M$,
all of whose composition factors $L(\gamma)$ correspond to elements 
$\gamma \in \Gamma$.

\begin{theorem} \label{lick} (\cite{CPS}, 3.5)
Let ${\cal C}$ be a highest weight category with respect to $\Lambda$.
Let $\Gamma \subset \Lambda$ be a finitely generated ideal, and let
$\Omega \subset \Lambda$ be a finitely generated coideal. Suppose that
$\Gamma \cap \Omega$ is a finite set.

There exists a quasi-hereditary algebra $S(\Gamma \cap \Omega)$ 
\index{$S(\Gamma \cap \Omega)$}
with poset $\Gamma \cap \Omega$, unique up to Morita equivalence,
such that the derived category
$D^b(S(\Gamma \cap \Omega)-mod)$ may be identified as the full subcategory
of $D^b({\cal C})$ represented as complexes of finite sums of modules
$I(\omega)_\Gamma$, with $\omega \in \Gamma \cap \Omega$. $\Box$
\end{theorem}

\begin{remark} If $S$ is a quasi-hereditary algebra with module category
${\cal C}$, under the hypotheses of theorem~\ref{lick}
we may choose $S(\Gamma \cap \Omega)$ to be a subquotient
$i (S /SjS) i$ of $S$, where $i,j$ are certain idempotents in $S$. 
\end{remark}

\begin{theorem} \label{raf} (\cite{CPS}, 3.4, 3.11, \cite{CPS0}, 4.3b)
Let $S$ be a quasi-hereditary algebra, with respect to a 
poset $\Lambda$. Then,

\bigskip

(a) $\Lambda$ indexes a collection $\{ \Delta(\lambda) \}_{\lambda \in \Lambda}$ of
``standard objects'' of ${\cal C}$, for each of which there exists 
a surjection
$\phi_\lambda: \Delta(\lambda) \rightarrow L(\lambda)$, 
such that all composition factors $L(\mu)$ of
$ker(\phi_\lambda)$ satisfy $\mu < \lambda$.

\bigskip

(b) The projective cover $P(\lambda)$ \index{$P(\lambda)$}
of $L(\lambda)$ possesses a filtration,
$$P(\lambda) = G_0(\lambda) \supset G_1(\lambda) \supset ... \supset 
G_N(\lambda) = 0,$$ 
such that,
\begin{description}
\item (i) $G_0(\lambda)/G_1(\lambda) \cong \Delta(\lambda)$.
\item (ii) For $n>0$, we have $G_n(\lambda)/G_{n+1}(\lambda) \cong \Delta(\mu)$,
for some $\mu > \lambda$.
\end{description}

\bigskip

(c) For projective objects $P$ in ${\cal C}$,
the number $[P : \Delta(\lambda)]$ of objects $\Delta(\lambda)$
appearing in a filtration by standard objects, is independent of filtration, 
for $\lambda \in \Lambda$.

For injective objects $I$ in ${\cal C}$,
the number $[I : \nabla(\lambda)]$ of objects $\nabla(\lambda)$
appearing in a filtration by costandard objects, is independent of filtration, 
for $\lambda \in \Lambda$.

\bigskip

(d) $[I(\mu) : \nabla(\lambda)] = [\Delta(\lambda) : L(\mu)]$,
for $\lambda, \mu \in \Lambda$.

Dually, $[P(\mu) : \Delta(\lambda)] = [\nabla(\lambda) : L(\mu)]$,
for $\lambda, \mu \in \Lambda$. 

\bigskip

(e) The category $mod-S$ of right modules over $S$ is a highest weight category. $\Box$
\end{theorem}

\begin{remark}
Let $S$ be a finite dimensional algebra, whose simple modules 
$\{ L(\lambda) \}_{\lambda \in \Lambda}$ are parametrised by a poset 
$\Lambda$.
Suppose that $S$ satisfies conditions (a) and (b) of theorem~\ref{raf}.
Then $S$ is quasi-hereditary with respect to $\Lambda$, by a dual 
to theorem~\ref{raf}.
\end{remark}
 
\begin{definition}
The \emph{decomposition matrix} of a quasi-hereditary algebra $S$ is 
the matrix $(d_{\lambda \mu})$ \index{$(d_{\lambda \mu})$}
of composition multiplicities 
$([\Delta(\lambda) : L(\mu)])$, 
whose rows amd columns are indexed by $\Lambda$.
\end{definition}

Let $S$ be a quasi-hereditary algebra with respect to $\Lambda$.
We describe some elements of the theory of tilting modules for $S$, due
to C. Ringel \cite{Ringel} (see also \cite{Donkin}, A4).

\begin{definition}
A \emph{tilting module} for $S$ is a finite dimensional $S$-module, which 
may be filtered by standard modules, and may also be filtered by
costandard modules.
\end{definition}

\begin{theorem} (\cite{Donkin}, A4, theorem 1)
For $\lambda \in \Lambda$, there is an indecomposable 
tilting module $T(\lambda)$, unique up to isomorphism, such that
$[T(\lambda) : L(\lambda)] = 1$, and all composition factors $L(\mu)$ of
$T(\lambda)$ satisfy $\mu \leq \lambda$.

Every tilting module for $S$ is a direct sum of modules 
$T(\lambda), \lambda \in \Lambda$. $\Box$
\end{theorem}

\begin{definition}
A \emph{full tilting module} for $S$ is a tilting module, in which every $T(\lambda)$
occurs as a direct summand.

A \emph{Ringel dual} $S'$ of $S$ is defined to be $End_S(T)^{op}$, where $T$ is 
a full tilting module for $S$. 
\end{definition}

\begin{remark}
The Ringel dual of $S$ is unique, up to
Morita equivalence.

We say the bimodule $_S T_{S'^{op}}$ defines a Ringel duality between $S, S'$.
\end{remark}

\begin{theorem} (\cite{Donkin}, A4, theorem 2) 
The Ringel duals $S'$ of $S$ are quasi-hereditary algebras, with respect to
the poset $\Lambda^{op}$, opposite to $\Lambda$. 

The standard $S'$-module corresponding to $\lambda \in \Lambda^{op}$ is given by
$$\Delta'(\lambda) = Hom_S(T, \nabla(\lambda)).$$
Dually, the costandard $S'$-module corresponding to $\lambda \in \Lambda^{op}$ is given by
$$\nabla'(\lambda) = \Delta^r(\lambda) \otimes_S T,$$
where $\Delta(\lambda)^{r}$ denotes the standard right module for $S$.  $\Box$
\end{theorem}

\begin{center}
{\large $q$-Schur algebras.}
\end{center}

S. Donkin, and R. Dipper have associated,
to a natural number $n$, and a non-zero element $q \in k$, \index{$q$}
a bialgebra ${\cal A}_q(n)$ \cite{DipDon}.
This bialgebra is a $q$-deformation of the bialgebra 
${\cal A}(n) = k[M]$ \index{$M$} of regular functions on the associative algebra
$M = M_n(k)$, of $n \times n$ matrices over $k$.
Whilst the undeformed bialgebra ${\cal A}(n) = {\cal A}_1(n)$ is commutative,
${\cal A}_q(n)$ is noncommutative, for $q \neq 1$. 

In general, ${\cal A}_q(n)$ may be decomposed by degree as a direct sum,
$${\cal A}_q(n) = \bigoplus_{r \geq 0} {\cal A}_q(n,r)$$ 
\index{${\cal A}_q(n), {\cal A}_q(n,r)$}
of finite dimensional coalgebras.
Thus, upon writing ${\cal S}_q(n,r)$ 
for ${\cal A}_q(n,r)^*$, the bialgebra
${\cal S}_q(n)$, \index{${\cal S}_q(n), {\cal S}_q(n,r)$} defined to be the graded dual
of ${\cal A}_q(n)$, decomposes as a direct sum,
$${\cal S}_q(n) = \bigoplus_{r \geq 0} {\cal S}_q(n,r)$$ 
of finite dimensional algebras.
The algebras ${\cal S}_q(n,r)$ are the $q$-Schur algebras, first introduced
by R. Dipper and G. James \cite{DiJ}.

\bigskip

Let $\Lambda(n,r)$ \index{$\Lambda(n,r)$}
be the poset of partitions of $r$ with $n$ parts or fewer,
with the dominance ordering $\unlhd$.

\begin{theorem} (\cite{Donkin}, 0.22)
The $q$-Schur algebra ${\cal S}_q(n,r)$ is a highest weight category with
respect to $\Lambda(n,r)$. $\Box$
\end{theorem}

Let $p$ be the least natural number such that 
$1+q+...+q^{p-1}=0$, if such exists. Otherwise, let $p= \infty$.

\begin{theorem} (\cite{Donkin}, 4.3(7))
If $p= \infty$, then ${\cal S}_q(n,r)$ is semisimple for all $n,r$.
In this case, $\Delta(\lambda) \cong L(\lambda)$, for all 
$\lambda \in \Lambda(n,r)$.
\end{theorem}

\begin{remark}
The algebra ${\cal S}_q(n,r)$ possesses a natural 
anti-automorphism $\sigma$, inherited from the transpose operation on the matrix algebra $M$ 
(\cite{Donkin}, 4.1).

It is thus possible to twist a left/right 
${\cal S}_q(n,r)$-module by $\sigma$, to obtain a right/left ${\cal S}_q(n,r)$-module.
Composing this twist with the duality functor,
we may associate to a left/right ${\cal S}_q(n,r)$-module $N$,
a left/right ${\cal S}_q(n,r)$-module $N^*$, \index{$N^*$}
its \emph{contragredient dual}.

The contragredient dual of a standard module $\Delta(\lambda)$
is isomorphic to the costandard module, $\nabla(\lambda)$.
The contragredient dual of a costandard module $\nabla(\lambda)$ 
is isomorphic to the standard module, $\Delta(\lambda)$. 
\end{remark}

\begin{remark}
In general, ${\cal S}_q(n,1)$ is isomorphic to the matrix algebra $M$, and therefore 
has a unique irreducible left module $E$ \index{$E$}, 
and a unique irreducible right module $E^{op}$.

S. Donkin has defined the $q$-exterior powers 
$\bigwedge^r _q E$ 
(respectively $\bigwedge^r _q E^{op}$) of $E$ (respectively $E^{op}$), 
which are left (respectively right) ${\cal S}_q(n,r)$-modules of dimension  
$\left( \begin{array}{c} n \\
r \end{array} \right)$,
exchanged under the anti-automorphism $\sigma$ (\cite{Donkin}, 1.2). 
He has also defined 
$q$-exterior powers $\bigwedge^r _q M$ \index{$\bigwedge^r _q M$} 
of $M$, which are 
${\cal S}_q(n,r)$-${\cal S}_q(n,r)$-bimodules of dimension 
$\left( \begin{array}{c} n^2 \\ r \end{array} \right)$ 
(\cite{Donkin}, 4.1).

For a sequence $\alpha = (\alpha_1,...., \alpha_m)$ of natural numbers, whose
sum is $r$, let us
define 

\begin{center}
$\bigwedge_q^\alpha E =
\bigwedge_q^{\alpha_1} E \otimes ... \otimes \bigwedge_q ^{\alpha_m} E.$
\end{center}
\index{$\bigwedge_q^\alpha E$}
\begin{center}
$\bigwedge_q^\alpha E^{op} =
\bigwedge_q^{\alpha_1} E^{op} \otimes ... \otimes \bigwedge_q ^{\alpha_m} E^{op}.$
\end{center}
\end{remark}

\begin{theorem} \label{donut} (S. Donkin, \cite{Donkin}, 1.2, 4.1)

(a) For $\tau \in \Sigma_m = Sym \{ 1,...,m \}$, there is an ${\cal S}_q(n,r)$-module
isomorphism between the exterior powers,
$\bigwedge_q^{(\alpha_1,...,\alpha_m)} E$, and
$\bigwedge_q^{(\alpha_{\tau 1},..., \alpha_{\tau m})} E$.

\bigskip

(b) There is a non-degenerate bilinear form,

\begin{center}
$<,>: \bigwedge _q^\alpha E^{op} \times \bigwedge _q^\alpha E \rightarrow k,$
\end{center}

{\noindent such that $<x \circ s,y> = <x, s \circ y>$.
Therefore, $\bigwedge_q^\alpha E \cong (\bigwedge_q^\alpha E)^*$,
as left ${\cal S}_q(n,r)$-modules,
and $\bigwedge_q^\alpha E^{op} \cong (\bigwedge_q^\alpha E^{op})^*$,
as right ${\cal S}_q(n,r)$-modules.}

\bigskip

(c) Direct summands 
of $q$-exterior powers $\bigwedge_q^\alpha E$ are tilting modules for ${\cal S}_q(n,r)$.
The restriction of $\bigwedge_q ^r M$ to a left ${\cal S}_q(n,r)$-module 
is a full tilting module. 
The restriction of $\bigwedge_q ^r M$ to a right ${\cal S}_q(n,r)$-module 
is also a full tilting module. 

\bigskip

(d) Let $n \geq r$. The ${\cal S}_q(n,r)$-${\cal S}_q(n,r)$-bimodule
$\bigwedge^r _q M$ defines a Ringel duality between ${\cal S}_q(n,r)$ 
and ${\cal S}_q(n,r)^{op} = {\cal S}_q(n,r)$. 
\end{theorem}

\begin{center}
{\large Hecke algebras associated to symmetric groups.}
\end{center}

Let $\Sigma_r = Sym \{1,2,...,r \}$ \index{$\Sigma_r$}
be the symmetric group on $r$ letters.

\begin{definition} The Hecke algebra ${\cal H}_q(\Sigma_r)$ 
\index{${\cal H}_q(\Sigma_r)$}
associated to $\Sigma_r$, 
is the associative $k$-algebra with generators
$$\{ T_{i} | i = 1,...,r-1 \},$$ subject to the relations,
$$T_{i}T_{i+1}T_{i} =
T_{i+1}T_{i}T_{i+1}, \hspace{0.5cm} i=1,...,r-2,$$
$$T_{i}T_{j} = T_{j}T_{i}, \hspace{0.5cm} |j-i|>1,$$
$$(T_{i} - q)(T_{i} +1) = 0, \hspace{0.5cm} i=1,...,r-1.$$
\end{definition}

When $q=1$, the Hecke algebra ${\cal H}_q(\Sigma_r)$ 
is isomorphic to the group algebra
$k\Sigma_r$. Many of the properties of $k\Sigma_r$ generalise to
${\cal H}_q(\Sigma_r)$.
For example, we have the following theorem.

\begin{theorem}  (\cite{Humph}, chapter 7)

(a)  ${\cal H}_q(\Sigma_r)$ 
possesses an outer automorphism $\#$, the ``sign automorphism'', 
which takes $T_{i}$ to $-T_{i} + q-1$. \index{$\#$}

\bigskip

(b) Given $w \in \Sigma_r$, and a reduced expression
$w = s_1...s_{l(w)}$ as a product simple transpositions 
$s_j \in \{ (i,i+1), 1 \leq i \leq i-1 \}$, we may define an element 
$T_w = T_{s_1}...T_{s_{l(w)}} \in {\cal H}_q(\Sigma_r)$, where $T_{(i,i+1)} = T_i$.
The element $T_w$ is independent of the choice of reduced expression.

\bigskip

(c) The set $\{ T_{w} | w \in \Sigma_r \}$ is a basis for ${\cal H}_q(\Sigma_r)$.
$\Box$
\end{theorem}

\begin{remark} There is a one dimensional ``trivial module'' for ${\cal H}_q(\Sigma_r)$, denoted
$k$, on which $T_w$ acts as $q^{l(w)}$.
The twist of the trivial module by $\#$ is the ``sign module'', denoted $sgn$, on which $T_w$ acts 
as $(-1)^{l(w)}$. \index{$sgn$}
\end{remark}

There is an elementary relation between Hecke algebras, and $q$-Schur
algebras, which is commonly exploited to describe the representation theory of
the Hecke algebra. A proof, with references to various sources, is given in \cite{KSX}, 1.2.

\begin{theorem} \label{tensorspace} ``Schur-Weyl duality''

(a) The ${\cal S}_q(n,r)$-module $E^{\otimes r}$ is 
a tilting module.

\bigskip

(b) There is an algebra surjection 
${\cal H}_q(\Sigma_r) \twoheadrightarrow End_{{\cal S}_q(n,r)}(E^{\otimes r})$.

\bigskip

(c) ${\cal S}_q(n,r) = End_{{\cal H}_q(\Sigma_r)}(E^{\otimes r})$.
\end{theorem}

\begin{remark} 
If $\lambda$ is a partition, then $\lambda'$ denotes the conjugate partition.

A partition $\lambda$ is said to be \emph{$p$-regular} if, and only if, 
$\lambda$ does not have $\geq p$ identical parts.

A partition $\lambda$ is said to be \emph{$p$-restricted} if, and only if, 
$\lambda'$ is $p$-regular.
\end{remark}

The ${\cal S}_q(n,r)$-${\cal H}_q(\Sigma_r)$-bimodule $E^{\otimes r}$ 
is often referred to as \emph{tensor space}.
If $n \geq r$, then this bimodule is particularly regular.

\begin{theorem} \label{tensorspace2}
(J.A. Green \cite{Green}, S. Donkin, \cite{Donkin}, 2.1, 4.4) Let $n \geq r$.

(a) There is an idempotent 
$\xi_{\omega} \in {\cal S}_q(n,r)$, \index{$\xi_\omega$}
such that
$$E^{\otimes r} \cong {\cal S}_q(n,r) \xi_{\omega},$$ 
as ${\cal S}_q(n,r)$-modules.

\bigskip

(b) The ${\cal S}_q(n,r)$-module $E^{\otimes r}$
is projective, and injective. 

\bigskip

(c) ${\cal H}_q(\Sigma_r) = End_{{\cal S}_q(n,r)}(E^{\otimes r}) = 
\xi_\omega {\cal S}_q(n,r) \xi_\omega$.

\bigskip

(d) The \emph{Schur functor}, 
$$Hom_{{\cal S}_q(n,r)}(E^{\otimes r},-):
{\cal S}_q(n,r)-mod \rightarrow {\cal H}_q(\Sigma_r)-mod,$$
is exact, and takes simple modules either to simple modules, or to zero.
We obtain all simple ${\cal H}_q(\Sigma_r)$-modules
from the set of simple 
${\cal S}_q(n,r)$-modules, in this way.

\bigskip

(e) If $\lambda$ is a $p$-nonrestricted partition, 
then $Hom_{{\cal S}_q(n,r)}(E^{\otimes r},L(\lambda))=0$.
Otherwise, if $\lambda$ is a $p$-restricted partition, then
$D_\lambda = Hom_{{\cal S}_q(n,r)}(E^{\otimes r},L(\lambda))$ \index{$D_\lambda$}
is a
simple ${\cal H}_q(\Sigma_r)$-module.

\bigskip

(f) The set, 
$$\{ D_\lambda \quad | \quad \lambda \textrm{ a }p \textrm{-restricted partition of }
r \},$$ is a 
complete set of non-isomorphic irreducible ${\cal H}_q(\Sigma_r)$-modules. 
$\Box$
\end{theorem}

\begin{remark}
Let $D^{\lambda'} = D_\lambda ^{\#}$. \index{$D^\lambda$}
It follows from theorem~\ref{tensorspace2}(f) that
$$\{ D^\lambda \quad | \quad \lambda \textrm{ a }p \textrm{-regular partition of }
r \},$$ is a 
complete set of non-isomorphic irreducible ${\cal H}_q(\Sigma_r)$-modules. 
\end{remark}

\begin{definition}
The \emph{Specht module} associated to a partition $\lambda$ is the 
${\cal H}_q(\Sigma_r)$-module, 
$$S^\lambda = Hom_{{\cal S}_q(r,r)}(E^{\otimes r},\nabla(\lambda)).$$
\index{$S^\lambda$}
The ${\cal H}_q(\Sigma_r)$-module $S_\lambda$ is defined to be, 
$$S_\lambda = Hom_{{\cal S}_q(r,r)}(E^{\otimes r},\Delta(\lambda)).$$
\index{$S_\lambda$}
The \emph{decomposition matrix} of ${\cal H}_q(\Sigma_r)$ is 
the matrix $(d_{\lambda \mu})$ \index{$(d_{\lambda \mu})$}
of composition multiplicities $[S^\lambda: D^\mu]$, 
indexed by partitions
$\lambda$ of $r$, and $p$-regular partitions $\mu$ of $r$.
\end{definition}

\begin{remark} \label{dualSpecht}

(a) It follows from theorem \ref{tensorspace2} that $S^{\lambda \#} \cong S_{\lambda'} \cong S^{\lambda' *}$.
We call $S_\lambda$ a dual Specht module.
 
\bigskip

(b) When $q=1$, and $k$ is a field of characteristic zero, the Specht modules $S^\lambda$ form a complete
set of non-isomorphic simple $k\Sigma_r$-modules.
We write $\chi^\lambda$ \index{$\chi^\lambda$}
for the irreducible character of $\Sigma_r$ corresponding to 
$S^\lambda$.

\bigskip

(c) 
The standard ${\cal S}_q(r,r)$ module $\Delta((1^r)) = \bigwedge{^r}(E)$ is one dimensional (the "determinant representation").
The Specht modules $S^{(r)}$ and $S^{(1^r)}$ are also one dimensional (the "trivial representation", and the "sign representation").
\end{remark}

\begin{definition}
The \emph{Young module} associated to a partition $\lambda$ is the ${\cal H}_q(\Sigma_r)$-module
$$Y^\lambda = Hom_{{\cal S}_q(r,r)}(E^{\otimes r},I(\lambda)),$$
\index{$Y^\lambda$} where $I(\lambda)$ is the injective hull of $L(\lambda)$.
The \emph{twisted Young module} associated to $\lambda$ is the ${\cal H}_q(\Sigma_r)$-module $Y^{\lambda \#}$.
\end{definition}

\begin{remark}

(a) By Schur-Weyl duality, $Y^\lambda$ is an indecomposable ${\cal H}_q(\Sigma_r)$-module, with a filtration by Specht modules 
$S^\mu$, with $\mu \unrhd \lambda$, and a single section $S^\lambda$.
Thus, $Y^{\lambda \#}$ is an indecomposable ${\cal H}_q(\Sigma_r)$-module, 
filtered by dual Specht modules $S_\mu$, with $\mu \unlhd \lambda'$, and a single section $S_{\lambda'}$.

\bigskip

(b) It follows from theorem \ref{tensorspace2} that $Y^\lambda$ is projective if, and only if, $\lambda$ is $p$-restricted.

\bigskip

(c) A partition $\lambda = (\lambda_i)$ of $r$ defines a Young subgroup $\Sigma_\lambda = \times_i \Sigma_{\lambda_i}$ of $\Sigma_r$.
Inducing the trivial representation from ${\cal H}_q(\Sigma_\lambda) = \otimes_i {\cal H}_q(\Sigma_{\lambda_i})$ up to 
${\cal H}_q(\Sigma_r)$,
we obtain a module $M^\lambda$, which is a direct summand of tensor space, as a ${\cal H}_q(\Sigma_r)$-module. 
This is because induction from Young subalgebras corresponds to taking tensor products, under Schur-Weyl duality. 
It follows that the module $M^\lambda$ is a direct sum of Young modules. 
In a direct sum decomposition, 
$M^\lambda$ has a unique indecomposable summand isomorphic to $Y^\lambda$, 
with all other summands isomorphic to $Y^\mu$, with $\mu \unrhd \lambda$.

\bigskip

(d) As a ${\cal H}_q(\Sigma_r)$-module, tensor space is isomorphic to a direct sum of Young modules.

\end{remark}

\begin{center}
{\large Finite general linear groups.}
\end{center}

Throughout this section, $q$ is a prime power, coprime to $l$.
We consider representations of the finite general linear group $GL_n(q)$,
over the field $k$.

\bigskip
 
Let $V$ be an $n$ dimensional vector space over $\mathbb{F}_q$. 
Let $\underline{n} = (n_1,\ldots,n_N)$ be a sequence of natural
numbers whose sum is $n$. 
Let us fix an $\mathbb{F}_q$-basis of $V$. 
We define a collection of subgroups of $GL(V) = GL_n(\mathbb{F}_q) = GL_n(q)$
\index{$GL_n(q)$}
relative to this basis: 
\bigskip

A \emph{maximal torus} $T(q)$ - the subgroup of diagonal matrices. \index{$T(q)$}

A \emph{Borel subgroup} $B(q)$ - the subgroup of upper triangular matrices. \index{$B(q)$}

A \emph{Levi subgroup} $L_{\underline{n}}(q)$ \index{$L_{\underline{n}}(q)$}- viewed as, 

\begin{center}
$\left(\begin{array}{cccc}
GL_{n_1}(q) & 0 & \ldots & 0 \\
0 & \ddots & \ddots & \vdots \\
\vdots & \ddots & \ddots & 0 \\
0 & \ldots & 0 & GL_{n_N}(q)
\end{array}\right)$
\end{center}

A \emph{parabolic subgroup} $P_{\underline{n}}(q)$ \index{$P_{\underline{n}}(q)$}- 
viewed as, 

\begin{center}
$\left(\begin{array}{cccc}
GL_{n_1}(q) & * & \ldots & * \\
0 & \ddots & \ddots & \vdots \\
\vdots & \ddots & \ddots & * \\
0 & \ldots & 0 & GL_{n_N}(q)
\end{array}\right)$
\end{center}

A \emph{unipotent subgroup} $U_{\underline{n}}(q)$ \index{$U_{\underline{n}}(q)$}- viewed as, 

\begin{center}
$\left(\begin{array}{cccc}
1_{n_1} & * & \ldots & * \\
0 & \ddots & \ddots & \vdots \\
\vdots & \ddots & \ddots & * \\
0 & \ldots & 0 & 1_{n_N}
\end{array}\right)$
\end{center}

The \emph{Weyl group} $W \cong \Sigma_n$ - permutation matrices. \index{$W$}
The following lemma is well known.

\begin{lemma} \label{Bruhat}

(a) The normaliser of $U_{\underline{n}}(q)$ is 
$P_{\underline{n}}(q) = U_{\underline{n}}(q) \rtimes L_{\underline{n}}(q)$.

\bigskip

(b) We have Bruhat decomposition,
$$GL_n(q) = \coprod_{w \in W} 
B(q)wB(q). \Box$$
\end{lemma}

\begin{definition}
Harish-Chandra induction is the functor,
$$HCInd_{L_{\underline{n}}(q)}^{GL_n(q)}: 
kL_{\underline{n}}(q)-mod \rightarrow kGL_n(q)-mod,$$
\index{$HCInd_{L_{\underline{n}}(q)}^{GL_n(q)}$}
$$M \mapsto k[GL_n(q)/U_{\underline{n}}(q)] \bigotimes_{L_{\underline{n}}(q)} M.$$
Harish-Chandra restriction is the functor,
$$HCRes^{GL_n(q)} _{L_{\underline{n}}(q)}: GL_n(q)-mod \rightarrow L_{\underline{n}}(q)-mod,$$
\index{$HCRes^{GL_n(q)} _{L_{\underline{n}}(q)}$}
$$M \mapsto k[U_{\underline{n}}(q) \backslash GL_n(q)] \bigotimes_{GL_n(q)} M.$$
\end{definition}

\begin{remark}
Because $q$ and $l$ are coprime,
the $kGL_n(q)$-$kL_{\underline{n}}(q)$ bimodule
$$k[GL_n(q)/U_{\underline{n}}(q)],$$ 
is a summand of the ordinary induction
bimodule $_{GL_n(q)} k[GL_n(q)]_{L_{\underline{n}}(q)}$, 
and so Harish-Chandra induction is exact.
Its left and right adjoint, Harish-Chandra restriction,
is also exact.
\end{remark}

Assuming the notation of J. Brundan, R. Dipper, and A. Kleshchev's article \cite{BDK}, 
we introduce a set  
to parametrise conjugacy classes in $GL_n(q)$.

\bigskip

For $\sigma \in \bar{\mathbb{F}}_q ^*$ of degree $d_\sigma$ \index{$d_\sigma$} over 
$\mathbb{F}_q$,
let $(\sigma)$  be the companion matrix representing $\sigma$ in 
$GL_{d_\sigma}(q)$.

For a natural number $k$, let 
$$(\sigma)^k = diag((\sigma),...,(\sigma)) =
\left(\begin{array}{cccc}
(\sigma) & 0 & \ldots & 0 \\
0 & \ddots & \ddots & \vdots \\
\vdots & \ddots & \ddots & 0 \\
0 & \ldots & 0 & (\sigma)
\end{array}\right)$$
\index{$(\sigma)^k$}
be the block diagonal matrix embedding $k$ copies of $(\sigma)$ in $GL_{d_\sigma k}(q)$.

\bigskip

For $\sigma, \tau \in \bar{\mathbb{F}}_q ^*$ of degree $d$ 
over $\mathbb{F}_q$,
let us write $\sigma \sim \tau$ if $\sigma$ and $\tau$ have the same minimal
polynomial over $\mathbb{F}_q$. 

\begin{definition}

(a) Let ${\cal C}^{pre}_{ss}$ \index{${\cal C}^{pre}_{ss}$} be the set,
$$\{ ((\sigma_1)^{k_1},...,(\sigma_a)^{k_a}) \quad | \quad
\sigma_i \in \bar{\mathbb{F}}_q ^*, \sigma_i \nsim \sigma_j \textrm{ for }
i \neq j, \sum d_{\sigma_i}k_i = n \}.$$

(b) Let $\sim$ denote an equivalence relation on ${\cal C}^{pre}_{ss}$, given by 
$$((\sigma_1)^{k_1},...,(\sigma_a)^{k_a}) \sim
((\tau_1)^{m_1},...,(\tau_a)^{m_b})$$ 
so long as, 
\begin{description}
\item (i) $a=b$,
\item (ii) there exists $w \in \Sigma_a$, such that $k_{wi} = m_i$, 
and $\sigma_{wi} \sim \tau_i$, for $i=1,...,a$.
\end{description}

(c) Let ${\cal C}_{ss} = {\cal C}^{pre}_{ss} / \sim$. \index{${\cal C}_{ss}$}
\end{definition}

For $s = ((\sigma_1)^{k_1},...,(\sigma_a)^{k_a}) \in {\cal C}_{ss}$,
let $\kappa(s) = (k_1,...,k_a)$. \index{$\kappa(s)$}

If $\underline{\lambda} = (\lambda_1,...,\lambda_a)$ 
is a sequence of partitions whose degrees are described by the sequence 
$\kappa$, we write $\underline{\lambda} \vdash \kappa$.

The representation theory of the 
principal ideal domain $\mathbb{F}_q[x]$ implies the following result.

\begin{theorem} (Jordan decomposition)
The conjugacy classes of $GL_n(q)$ are in one-one correspondence with the
set, 
$$\{ (s, \underline{\lambda}) \quad | \quad s \in {\cal C}_{ss}, 
\underline{\lambda} \vdash \kappa(s) \}. \Box$$
\end{theorem}

\begin{theorem} (J.A. Green, \cite{Gr}, \cite{BDK}, 2.3) 
Fix an embedding 
$\bar{\mathbb{F}}_q ^* \hookrightarrow \bar{\mathbb{Q}}_l ^*$.
There is a one-one correspondence between the set of ordinary
irreducible characters of $GL_n(q)$ and the set,
$$\{ (s, \underline{\lambda}) \quad | \quad s \in {\cal C}_{ss}, 
\underline{\lambda} \vdash \kappa(s) \}. \Box$$
\end{theorem}

Given an embedding $\bar{\mathbb{F}}_q ^* \hookrightarrow \bar{\mathbb{Q}}_l ^*$, and
$ s \in {\cal C}_{ss}, 
\underline{\lambda} \vdash \kappa(s)$,
we write $\chi(s, \underline{\lambda})$ \index{$\chi(s, \underline{\lambda})$}
for the irreducible character of $GL_n(q)$  
corresponding to the pair $(s, \underline{\lambda})$ in the above theorem.

\begin{remark}
Let $s = ((\sigma_1)^{k_1},...,(\sigma_a)^{k_a})$, 
and $\underline{\lambda} \vdash \kappa(s)$.
The irreducible character $\chi(s, \underline{\lambda})$ 
corresponding to the pair $(s, \underline{\lambda})$, is equal to the 
Harish-Chandra induced character,
$$HCInd_{GL_{d_{\sigma_1}k_1}(q) \times...\times GL_{d_{\sigma_a}k_a}(q)}
^{GL_n(q)} \left( \chi((\sigma_1)^{k_1}, \lambda_1) \otimes... \otimes 
\chi((\sigma_a)^{k_a}, \lambda_a) \right).$$
Furthermore, if $\mu, \nu$ are partitions of $m,n$, and 
$\sigma \in \bar{\mathbb{F}}_q ^*$ has degree $d$, then, 
$$HCInd_{GL_{dm}(q) \times GL_{dn(q)}} ^{GL_{dm+dn}(q)}
\left( \chi((\sigma)^{m}, \mu) \otimes
\chi((\sigma)^{n}, \nu) \right)$$ 
$$= \sum c(\lambda; \mu, \nu) \chi( (\sigma)^{m+n}, \lambda),$$
where $c(\lambda; \mu, \nu)$ \index{$c(\lambda; \mu, \nu)$}
is the Littlewood-Richardson coefficient
associated to $\lambda, \mu, \nu$.
\end{remark}

Let ${\cal C}_{ss,l'} = 
\{ s \in {\cal C}_{ss} \quad | \quad \textrm{the order of }s$ is coprime to $l \}$.
\index{${\cal C}_{ss,l'}$}

\begin{theorem} (P. Fong and B. Srinivasan, \cite{BDK}, 2.4.6)
There is a decomposition of the group algebra $kGL_n(q)$ into a direct
sum of two sided ideals,
$$kGL_n(q) = \bigoplus_{s \in {\cal C}_{ss,l'}} kB_s,$$
where the characters in $B_s$ \index{$kB_s$} 
are given by the set,
$$\{ \chi(t, \underline{\lambda}) \quad | \quad \underline{\lambda} \vdash \kappa(t),
t \in {\cal C}_{ss} \textrm{ has }l \textrm{-regular part conjugate to }s \}.
\Box$$
\end{theorem}

\begin{definition}
The characters $\{ \chi(1, \lambda) \quad | \quad \lambda \vdash n \}$ of $GL_n(q)$ are the 
\emph{unipotent characters}.

The indecomposable components of the ${GL_n(q)}$ -$GL_n(q)$-bimodule $kB_1$ are the \emph{unipotent blocks}.
\end{definition}

Thanks to the following theorem, we may 
concentrate our curiosity on unipotent blocks.

\begin{theorem} (C. Bonnaf\'e and R. Rouquier, \cite{BoRou}, 11.8)
Every block of $GL_n(q)$ is Morita equivalent to a unipotent block. $\Box$
\end{theorem}

\begin{theorem} \label{Swd} 
(N. Iwahori, E. Cline, B. Parshall, L. Scott, see \cite{BDK}, 3.2d, 3.5a)
Let ${\cal M} = k[GL_n(q) / B(q)]$ \index{${\cal M}$}
denote the permutation module
for $kGL_n(q)$ on the coset space $GL_n(q) / B(q)$.

\bigskip

(a) If $q \neq 1$ modulo $l$, then ${\cal M}$ is a projective $kGL_n(q)$-module.

\bigskip

(b) ${\cal M}$ is a $kB_1$-module. 

\bigskip

(c) $End_{kB_1}({\cal M}) \cong {\cal H}_q(\Sigma_n)$.
Under this isomorphism, the endomorphism defined 
by the double coset $B(q)wB(q)$ maps to the element 
$T_w$ of ${\cal H}_q(\Sigma_n)$, for $w \in \Sigma_n$.

\bigskip

(d) The annihilator ideal $Ann_{kB_1}({\cal M})$ is nilpotent. 

\bigskip

(e) The quotient
$kB_1 / Ann_{kB_1}({\cal M})$ is Morita equivalent to ${\cal S}_q(n,n)$.

The standard ${\cal S}_q(n,n)$-module $\Delta(\lambda)$ corresponds,
under this Morita equivalence, 
to a $kB_1$-module, which is an $l$-modular reduction of the character
$\chi(1, \lambda)$ of $GL_n(q)$.

\bigskip

(f) The ${\cal S}_q(n,n)$-${\cal H}_q(\Sigma_n)$-bimodule corresponding to
the $kB_1$-${\cal H}_q(\Sigma_n)$-bimodule ${\cal M}$ under the Morita equivalence of (d), 
is isomorphic to twisted tensor space $E^{\otimes n \#}$.
$\Box$
\end{theorem}

\begin{definition}
By theorem~\ref{Swd}, simple $kB_1$-modules are in one-one correspondence with 
simple modules for the $q$-Schur algebra, which are in one-one correspondence with the set 
$\Lambda(n) = \Lambda(n,n)$. \index{$\Lambda(n)$}
We denote by $D(\lambda)$ \index{$D(\lambda)$} 
the simple $kB_1$-module corresponding to $\lambda \in \Lambda(n)$.
\end{definition}

\begin{remark}
By twisted tensor space $E^{\otimes n \#}$, we mean the tensor space of theorem~\ref{tensorspace}, 
on which the right action of ${\cal H}_q(\Sigma_n)$  
has been twisted by the signature
automorphism, $\#$.
\end{remark}
    
\begin{center}
{\large Wreath products.}
\end{center}

For a $k$-algebra $A$, and a natural number $w$,
the wreath product $A \wr \Sigma_w$, \index{$A \wr \Sigma_w$}
is defined to be the space
$A^{\otimes w} \otimes k\Sigma_w$, with associative multiplication,
$$(x_1 \otimes ... \otimes x_w \otimes \sigma).(y_1 \otimes... \otimes y_w 
\otimes \tau) =$$
$$x_1 .y_{\sigma^{-1} 1} \otimes ...\otimes x_w .y_{\sigma^{-1} w}
\otimes \sigma.\tau,$$
for $x_1,...,x_w,y_1,...,y_w \in A$, and $\sigma, \tau \in \Sigma_w$.
We collect a jumble of information concerning wreath products.
One reference for such stuff is an article of J. Chuang and K.M. Tan \cite{ChTan3}.  

\bigskip

For an $A$-module $L$, let us define the $A \wr \Sigma_w$-module $T^{(w)} L$ 
\index{$T^{(w)} L$}
to be
the $w$-fold tensor product $L^{\otimes w}$, on which the subalgebra
$A^{\otimes w}$ acts component-wise, and the symmetric group $\Sigma_w$ acts
by place permutations.

\bigskip

For an $A \wr \Sigma_w$-module $M$, and a $k\Sigma_w$-module $N$, we
define an $A \wr \Sigma_w$-module $M \otimes N$, by the action
$$(\alpha \otimes \sigma)(m \otimes n) = 
(\alpha \otimes \sigma)m \otimes \sigma n,$$ 
for $\alpha \in A^{\otimes w}, \sigma \in \Sigma_w, m \in M, n \in N$.

\bigskip

Let $I$ be a finite set.
Given a set of natural numbers $\{ w_i, i \in I \}$, 
whose sum is $w$, there is a Young subgroup of $\Sigma_w$, isomorphic to
$\times_{i \in I} \Sigma_{w_i}$. Accordingly, there is
a subalgebra of $A \wr \Sigma_w$, isomorphic to
$\bigotimes_{i \in I} A \wr \Sigma_{w_i}$.

\bigskip

Let $\Lambda_l$ be the set of all $l$-regular partitions.
Let $\Lambda_{l,w} ^I$ \index{$\Lambda_{l,w} ^I$} 
be the set of $I$-tuples 
$(\lambda^i)_{i \in I}$ of $l$-regular partitions,
whose orders $(w_i)_{i \in I}$ sum to $w$.

\begin{theorem} (I. MacDonald, \cite{MacD})
Let $\{ S(i), i \in I \}$ be a complete set of non-isomorphic
simple $A$-modules.
For $\underline{\lambda} \in \Lambda_{l,w} ^I$, let
$$S(\underline{\lambda}) =
Ind_{\otimes_{i \in I} A \wr \Sigma_{w_i}} 
^{A \wr \Sigma_w} \left( \bigotimes_{i \in I} 
T^{w_i}(S(i)) \otimes D^{\lambda^i} \right) .$$
\index{$S(\underline{\lambda})$}
The set
$\{ S(\underline{\lambda}) \quad | \quad \lambda \in \Lambda_{l,w} ^I \}$,
is a complete
set of non-isomorphic simple $A \wr \Sigma_w$-modules.
$\Box$
\end{theorem}

Let $q$ be a prime power, coprime to $l>0$.
Let $B(q)$ be a Borel subgroup of $GL_n(q)$. 
Thus, the direct product $\times^w B(q)$ of $w$ copies of $B(q)$ is a subgroup of the base group
$\times^w GL_n(q)$ inside $GL_n(q) \wr \Sigma_w$.

\bigskip

The following theorem is a straightforward generalisation of theorem~\ref{Swd}.

\begin{theorem} \label{Kill}
Let ${\cal M}_w$ \index{${\cal M}_w$} 
be the $kB_1 \wr \Sigma_w$-module 
$k[GL_n(q) \wr \Sigma_w / \times^w B(q)]$. 

\bigskip

(a) If $q \neq 1$ modulo $l$, then ${\cal M}_w$ is a projective 
$kGL_n(q) \wr \Sigma_w$-module.

\bigskip

(b) $End_{GL_n(q) \wr \Sigma_w} ({\cal M}_w) \cong
{\cal H}_q(\Sigma_n) \wr \Sigma_w$.

\bigskip 

(c) The simple modules sent to zero by the functor, 
$$Hom_{GL_n(q) \wr \Sigma_w}({\cal M}_w, -),$$ are precisely those simples
$S(\underline{\lambda}), \underline{\lambda} \in \Lambda_{l,w} ^{\Lambda(n)}$,
for which some entry, 
$\lambda^\mu$ 
in $\underline{\lambda}$, 
is non-zero,
for some $p$-nonrestricted partition $\mu$. $\Box$
\end{theorem}

\begin{proposition} \label{pentagon}
The wreath product ${\cal H}_q(\Sigma_n) \wr \Sigma_w$ 
is a symmetric algebra.
\end{proposition}
Proof:

Define a form $<,>$ on ${\cal H}_q(\Sigma_n)$, 
extending the following form bilinearly:
$$< T_u,T_v > =
\left\{ 
\begin{array}{cc}
q^{l(u)} & \textrm{ if }u=v^{-1} \\
0 & \textrm{otherwise.}
\end{array}
\right.$$
R. Dipper and G. James have proved (\cite{JamDip3}, (2.3)) 
that this is a symmetric, associative bilinear form.
It is also non-degenerate, since $q$ is a unit in
$k$.

In case $n=w,q=1$, one obtains a symmetric bilinear form $<,>_\Sigma$
on $k\Sigma_w$.
Now define a bilinear form on ${\cal H}_q(\Sigma_n) \wr \Sigma_w$, by the 
formula,
$$< x_1 \otimes ... \otimes x_w \otimes \sigma, y_1 \otimes... \otimes y_w 
\otimes \tau >_\wr =$$
$$<x_1, y_{\sigma^{-1} 1}> ...<x_w, y_{\sigma^{-1} w}><\sigma, \tau>_\Sigma,$$
for $x_1,...,x_w,y_1,...,y_w \in {\cal H}_q(\Sigma_n)$, and 
$\sigma,\tau \in \Sigma_w$.
Then $<,>_\wr$ is an associative, symmetric, 
non-degenerate bilinear form on ${\cal H}_q(\Sigma_n) \wr \Sigma_w$.
$\Box$

\bigskip

\newpage

\begin{center}
{\bf \large 
Chapter II

Blocks of $q$-Schur algebras, Hecke algebras, and finite general linear
groups.}
\end{center}

Let $k$ be a field of characteristic $l$. We summarise some elements of 
block theory for the $q$-Schur algebras, the Hecke algebras, 
and the finite general linear groups in non-describing characteristic.
We introduce the Rock blocks, and state some known results concerning them.

\begin{center}
{\large Abacus combinatorics.}
\end{center}

Let $p$ be a natural number. It is often
convenient for us to picture partitions on an abacus with
$p$ runners, following G. James (\cite{JamKerb}, pg. 78).

We thus label the runners on an abacus $0,...,p-1$, from left to right,
and label its rows $0,1,...$, from the top downwards.
If $\lambda=(\lambda_1,\lambda_2,...)$ is a partition 
with $m$ parts or fewer then we may represent $\lambda$ on 
the abacus with $m$ beads: 
for $i=1,...,m$, write $\lambda_i+m-i=s+pt$, with 
$0 \leq s \leq p-1$, and place a bead on the $s^{th}$
runner in the $t^{th}$ row. For example, 

\bigskip

\begin{picture}(170,90)(-150,10)
\put(0,100){\makebox(0,0){0}}
\put(30,100){\makebox(0,0){1}}
\put(60,100){\makebox(0,0){2}}
\put(0,95){\line(1,0){60}}
\put(0,95){\line(0,-1){75}}
\put(30,95){\line(0,-1){75}}
\put(60,95){\line(0,-1){75}}
\put(0,90){\circle*{5}}
\put(0,80){\circle*{5}}
\put(0,70){\circle*{5}}
\put(30,90){\circle*{5}}
\put(30,80){\circle*{5}}
\put(30,70){\circle*{5}}
\put(30,60){\circle*{5}}
\put(30,50){\circle*{5}}
\put(60,90){\circle*{5}}
\put(60,80){\circle*{5}}
\put(60,70){\circle*{5}}
\put(60,60){\circle*{5}}
\put(60,50){\circle*{5}}
\put(60,40){\circle*{5}}
\put(60,30){\circle*{5}}
\end{picture}

{\noindent is an abacus representation of the partition $(6,4,2^2,1^2)$, when $p=3$.}

Sliding a bead one row upwards on its runner, into a vacant position,
corresponds to removing a $p$-hook from the rim of the partition $\lambda$.
Given an abacus representation of a partition,
sliding all the beads up as far as possible produces an 
abacus representation of the $p$-core of that partition,
a partition from which no further rim $p$-hooks can be removed.
The pictured example $(6,4,2^2,1^2)$ is therefore a $3$-core.
The $p$-core of a partition is independent of the way in which $p$-hooks have 
been removed.
The $p$-weight of a partition is the number of $p$-hooks 
removed to obtain the $p$-core.

Fix an abacus representation of a partition $\lambda$,
and for $i=0,...,p-1$, let $\lambda^i_1$ be the number
of unoccupied positions on the $i^{th}$ runner which
occur above the lowest bead on that runner.
Let $\lambda^i_2$ be the number of unoccupied positions
on the $i^{th}$ runner which occur above the second lowest
bead on that runner, etc.
Then $\lambda^i = (\lambda^i_1, \lambda^i_2,...)$ is a partition,
and the $p$-tuple $[\lambda^0,...,\lambda^{p-1}]$ \index{$[\lambda^0,...,\lambda^{p-1}]$}
is named
the $p$-quotient of $\lambda$.
The $p$-quotient depends on the number of beads in the abacus 
representation of $\lambda$, but is well defined up to a cyclic permutation.
The $p$-weight of $\lambda$ is equal to the sum $|\lambda^0|+...+
|\lambda^{p-1}|$.
 
The partitions with a given $p$-core $\tau$ and $p$-weight $w$ \index{$w$}
may be parametrized by $p$-quotients: \index{$\tau$}

Fix $m$ so that any such partition has $m$ parts or fewer.
Representing the partitions on an abacus with $m$ beads,
there is a $p$-quotient for each one.
We thus 
introduce a bijection between the set of partitions with
$p$-core $\tau$ and weight $w$, and the set of $p$-tuples
$[\sigma^0,...,\sigma^{p-1}]$ such that 
$|\sigma^0|+...+|\sigma^{p-1}| = w$.
It is a convention of this article, frequently to label 
a partition by its $p$-quotient.

Whenever we represent partitions with a given $p$-core 
on an abacus,
we assume that $m$ is fixed as above.

\begin{center}
{\large Block parametrisation.}
\end{center}

Let $A_{{\cal S}_q}(n) = {\cal S}_q(n,n)$, let 
$A_{{\cal H}_q}(n) =  {\cal H}_q(\Sigma_n)$, and let $A_{G_q}(n) = kB_1$, 
a direct sum of the unipotent blocks of $kGL_n(q)$.

We give 
parametrizations for the blocks of $A_{\cal X}(n)$ \index{$A_{\cal X}(n)$}, for 
${\cal X} \in \{ {\cal S}_q, {\cal H}_q, G_q \}$.
The parametrizations depend on an invariant $p$.

\begin{definition}
If ${\cal X} \in \{ {\cal S}_q, {\cal H}_q \}$, then $p = p({\cal X})$ \index{$p$}
is defined to be the least natural number such that,
$$1+q+...+q^{p-1} = 0,$$ 
in $k$ if such exists, and $p = \infty$, otherwise.

If ${\cal X} = G_q$, for some prime power $q$, then we insist that $l>0$. 
In this case, $p = p({\cal X})$ is defined to be the 
multiplicative order of $q$, modulo $l$.
\end{definition}

\begin{remark}
Let $l>0$.
If $q$ is a prime power such that $q \neq 1$ modulo $l$, 
then $p(G_q) = p({\cal S}_q)$. 

However, if $q=1$ modulo $l$,
then
$p(G_q) = 1$, while $p({\cal S}_q) = l$.
\end{remark}

\begin{theorem} \label{block}
(R. Brauer, P. Fong, B. Srinivasan \cite{FoSri}, 7A, G. James, R. Dipper, \cite{JamDip3}, 4.13, 
S. Donkin, A. Cox, \cite{Cox})
The blocks of $A_{\cal X}(n)$ are in one-one correspondence
with pairs $(w, \tau)$, where $w \in \mathbb{N}_0$, and $\tau$ is a $p$-core
of size $n-wp$.
\end{theorem}

We write $b_{\tau,w} ^{\cal X}$ \index{$b_{\tau,w} ^{\cal X}$} for the block idempotent in 
$A_{\cal X}(n)$ corresponding to the pair $(w, \tau)$, defined over ${\cal O}$.
We write $k{\bf B}_{\tau,w} ^{\cal X}$ \index{$k{\bf B}_{\tau,w} ^{\cal X}$} 
for the corresponding block algebra.
We say that $b_{\tau,w} ^{\cal X}$ (respectively $k{\bf B}_{\tau,w} ^{\cal X}$) 
is a block of weight $w$, with core $\tau$.

\bigskip

In this book, it will be convenient for us to work over
an $l$-modular system $(K, {\cal O}, k)$  \index{$(K, {\cal O}, k)$}.
Thus, ${\cal O}$ is a discrete valuation ring with field of fractions $K$ of characteristic zero, 
maximal ideal $\wp$, and residue field $k = {\cal O}/ \wp$ of characteristic $l$.
For an ${\cal O}$-module $M$, we write $KM = K \otimes_{\cal O} M$, and $kM = k \otimes_{\cal O} M$.

Suppose $A_{\cal X}(n)$ is defined over ${\cal O}$.
Abusing notation, we then write $b_{\tau,w}^{\cal X}$ for the lift over ${\cal O}$ of the corresponding block idempotent defined over $k$.
We write $R{\bf B}_{\tau,w} ^{\cal X}$ for the corresponding block algebra over $R$, for $R \in \{ K, {\cal O}, k \}$.

\begin{remark} 
Here, and in the sequel, we write ``block'' in 
abbreviation of either ``block idempotent'', or ``block algebra''.
\end{remark}

\begin{theorem} 
(R. Brauer, P. Fong, B. Srinivasan \cite{FoSri}, 7A, G. James, R. Dipper, \cite{JamDip3}, 4.13, 
S. Donkin, A. Cox, \cite{Cox})
Let $\tau$ be a partition of $t$, and a $p$-core.
Let $n = pw+ t$. 

\bigskip

(a) A standard, costandard, or simple ${\cal S}_q(n,n)$-module 
indexed by the partition $\lambda$ lies in $k{\bf B}_{\tau,w} ^{{\cal S}_q}$ 
if, and only if, $\lambda$ has $p$-core $\tau$, and $p$-weight $w$.

\bigskip

(b) A Specht, or simple ${\cal H}_q(\Sigma_n)$-module indexed by 
the partition $\lambda$ lies in $k{\bf B}_{\tau,w} ^{{\cal H}_q}$ 
if, and only if, $\lambda$ has $p$-core $\tau$, and $p$-weight $w$.
 
\bigskip

(c) The irreducible character,
$$\chi \left( ((1)^{k_1},(\sigma_2)^{k_2},...,(\sigma_a)^{k_a}), 
(\lambda_1,...,\lambda_a) \right),$$ 
of $GL_n(q)$ lies in $k{\bf B}_{\tau,w} ^{G_q}$ 
if, and only if, $\lambda_1$ has $p$-core $\tau$, and the order of
$\sigma_i$ is a power of l, for $i \geq 2$. $\Box$
\end{theorem}

We write $k{\bf B}_{\tau,w} ^{\cal S}$ for the block
$k{\bf B}_{\tau,w} ^{{\cal S}_1}$ of the Schur algebra ${\cal S}(n,n)$.
We write $k{\bf B}_{\tau,w} ^{\Sigma}$ for the block
$k{\bf B}_{\tau,w} ^{{\cal H}_1}$ of the symmetric group $\Sigma_n$.

The following theorem is standard.
Part (a) follows from the studies of P. Fong and B. Srinivasan \cite{FoSri}.
Part (b) follows from work of R. Dipper and G. James \cite{JamDip1}.

\begin{theorem}
(a) 
The blocks $k{\bf B}_{\tau,w} ^{\Sigma}$ and $k{\bf B}_{\tau,w} ^{G_q}$ have abelian
defect groups if, and only if, $w<l$.

\bigskip

(b) If $w<l$, then a character 
$\chi \left( ((1)^{k_1},(\sigma_2)^{k_2},...,(\sigma_a)^{k_a}), 
(\lambda_1,...,\lambda_a) \right)$ in 
$k{\bf B}_{\tau,w} ^{G_q}$, is constrained by the 
condition that $\sigma_i$ is of degree $p$, for $i \geq 2$. $\Box$
\end{theorem}

\begin{center}
{\large Derived equivalences, and Rock blocks.}
\end{center}

\begin{theorem} (J. Chuang, R. Rouquier, \cite{ChRou})
Let ${\cal X}$ be an element of the set
$\{ {\cal S}_q, {\cal H}_q, G_q, {\cal S}, \Sigma \}$. Let $\tau, \tau'$ 
be $p$-cores. Then the bounded derived categories 
$D^b(k{\bf B}_{\tau,w} ^{\cal X})$, and $D^b(k{\bf B}_{\tau',w} ^{\cal X})$, are 
equivalent. $\Box$
\end{theorem}

\begin{definition}
Suppose $p,w$ are fixed.
A $p$-core $\rho$ \index{$\rho$}
is said to be a \emph{Rouquier core} if it has an abacus presentation,
on which there are at least $w-1$ more beads on runner $i$, than on runner
$i-1$, for $i=1,...,p-1$.
\end{definition}

\begin{example}
Let $p=w=3$. 
Then the partition $(6,4,2^2,1^2)$ pictured on the abacus at the beginning of the chapter
is a Rouquier core.
\end{example}

R. Rouquier conjectured the following structure theorem \cite{Rouq}.

\begin{theorem} \label{structsymm}
(J. Chuang, R. Kessar, \cite{KessCh})
Let $w < l$, and let $\rho$ be a Rouquier core.
The block $k{\bf B}_{\rho,w} ^\Sigma$ of a symmetric group is Morita
equivalent to the wreath product $k{\bf B}_{\emptyset,1} ^\Sigma \wr \Sigma_w$.
$\Box$
\end{theorem}

Let ${\cal X} \in \{ {\cal S}_q, {\cal H}_q, G_q, {\cal S}, \Sigma \}$.

\begin{definition} 
A \emph{Rock block} is any block $b_{\rho,w} ^{\cal X}$ (respectively $k{\bf B}_{\rho,w} ^{\cal X}$), 
where $\rho$ is a Rouquier core for $p,w$.
\end{definition}

\begin{remark}

(a) ``RoCK block'' is slang for ``Rouquier, or Chuang-Kessar block''.

\bigskip

(b) For fixed ${\cal X}, p, w$, all Rock blocks are Morita equivalent (cf. lemma \ref{cock}).

\bigskip

(c) Throughout this article, the letter $\rho$ represents a Rouquier core.
Other $p$-cores will be represented by different letters, such as $\tau$.
\end{remark}

\bigskip 

The next couple of results were deduced from theorem~\ref{structsymm}.

\begin{theorem} 
(J. Chuang, K.M. Tan, \cite{ChTan5})
Let $w<l$.
The Rock block
$k{\bf B}_{\rho,w} ^{\cal S}$ is Morita
equivalent to $k{\bf B}_{\emptyset,1} ^{\cal S} \wr \Sigma_w$.
$\Box$
\end{theorem}

For a natural number $a$, let $\Lambda^a _w$ \index{$\Lambda^a _w$}
be the set of $a$-tuples
$(\alpha^1,...,\alpha^a)$ of partitions $\alpha^i$, 
such that $\sum_i | \alpha^i | = w$.

Let 
$c(\lambda; \mu, \nu)$ be the 
Littlewood-Richardson coefficient corresponding to
partitions $\lambda, \mu, \nu$, taken to
be zero whenever $| \lambda| \neq |\mu| + | \nu|$.
 
\begin{theorem} (J. Chuang, K.M. Tan, \cite{ChTan5}, H. Miyachi, \cite{Miyachi})
Let $w < l$.
The decomposition matrix of the Rock block $k{\bf B}_{\rho,w} ^{\cal S}$ is equal to the
decomposition matrix of $k{\bf B}_{\emptyset,1} ^{\cal S} \wr \Sigma_w$. It is 
given by,
$$d_{\underline{\lambda} \underline{\mu}} =
\sum_{\alpha \in \Lambda^{p+1}  _w, \beta \in \Lambda^p _w}
\prod_{j=0} ^{p-1} c(\lambda^j; \alpha^j, \beta^j) 
c(\mu^j; \beta^j, (\alpha^{j+1})'). \Box$$
\end{theorem}

The following theorem was also observed, 
upon viewing the canonical basis
for the basic representation of $\widehat{\mathfrak{sl}_p}$.

\begin{theorem} \label{stolen} (J. Chuang and K.M. Tan, \cite{ChTan},
B. Leclerc, and H. Miyachi, \cite{Miy})
Let $l=0$, and let $q$ be a $p^{th}$ root of unity in $k$. 
The decomposition matrix of the Rock block $k{\bf B}_{\rho,w} ^{{\cal S}_q}$ is equal to the
decomposition matrix of $k{\bf B}_{\emptyset,1} ^{{\cal S}_q} \wr \Sigma_w$. It is 
given by,
$$d_{\underline{\lambda} \underline{\mu}} =
\sum_{\alpha \in \Lambda^{p+1}  _w, \beta \in \Lambda^p _w}
\prod_{j=0} ^{p-1} c(\lambda^j; \alpha^j, \beta^j) 
c(\mu^j; \beta^j, (\alpha^{j+1})'). \Box$$
\end{theorem}

\begin{lemma}
The $p$-regular partitions 
$[\lambda^0,\lambda^1...,\lambda^{p-1}]$ 
with $p$-core $\rho$, are those for which $\lambda^0$ is empty. $\Box$
\end{lemma}

\begin{theorem} \label{Rowena}
(R. Paget, \cite{Paget1}) Let $\rho$ be a Rouquier core. Then
$$D^{[\emptyset, \lambda^1,...,\lambda^{p-1}]} \cong 
D_{[\emptyset, \lambda^{p-1},...,\lambda^1]},$$
as $k{\bf B}_{\rho,w}^{{\cal H}_q}$-modules. $\Box$
\end{theorem}

\begin{center}
{\large Brauer correspondence for 
blocks of finite general linear groups.}
\end{center}

If $a \in {\cal O}$, we denote the $l$-modular reduction of $a$ by $\bar{a}$, an element of $k$. 
We insist that $K$ is a splitting field for $G$.

\begin{definition}
(see \cite{Bro1})
For a finite group $G$, with an $l$-subgroup $P$,
and an ${\cal O}G$-module $M$,
the Brauer homomorphism is defined to be the quotient map  
$$ Br^G _P: M^P \rightarrow
M(P) = M^P/( \sum_{Q<P} Tr_Q^P(M^Q) + 
\wp M^P).$$ \index{$M^P$} \index{$M(P)$} \index{$Br^G_P$}
{\noindent
The Brauer quotient $M(P)$ is the quotient 
of $P$-fixedpoints of $M$ by relative
traces from proper subgroups, reduced modulo $l$;
$M(P)$ is a $kN_G(P)$-module.}
\end{definition}

The Brauer homomorphism factors through the epimorphism ${\cal O}M^P \rightarrow kM^P$. 
Abusing notation, we write $Br_P^G$ for the induced map from $kM^P$ to $M(P)$.

\bigskip

In case that $G$ acts on $M = {\cal O}G$ by conjugation,
without doubt ${\cal O}G(P) = kC_G(P)$, 
and the quotient map
$$
Br^G_P: ({\cal O}G)^P \rightarrow kC_G(P)
$$
{\noindent
is an algebra homomorphism, the classical
Brauer homomorphism with 
respect to $P$, given by the rule}
$$
Br^G_P(\sum_{g \in G} a_g g) = \sum_{g \in C_G(P)} 
\bar{a}_g g.
$$

\begin{theorem} (L. Scott)
Let $M$ be a $kG$-permutation module.
Then $M(P) \neq 0$ 
if, and only if, $M$ has a direct summand with a vertex 
containing $P$. $\Box$
\end{theorem}

We present Brauer correspondence for unipotent blocks of 
finite general linear groups.
Thus, let $l>0$, and let $q$ be a prime power, coprime to $l$. 

\bigskip

Let $t = n -wp \geq 0$. Let
$$L_1 = \times^w GL_p(q) < GL_{wp}(q),$$
$$L_2 = GL_{wp}(q) \times GL_{t}(q) < GL_n(q),$$
be Levi subgroups.
Let $D$ be a Sylow $l$-subgroup of $GL_{wp}(q)$. 

\bigskip

Let $\tau$ be a partition of $t$, and a $p$-core.
Let $b_{\tau,w} ^{G_q}$ (respectively $b_{\tau,0} ^{G_q}$) be the 
unipotent block of $GL_n(q)$
(respectively $GL_{t}(q))$ with $p$-core $\tau$.

The centralizer and normalizer of $D$ in $GL_n(q)$ are contained in $L_2$:

\begin{lemma} (\cite{FoSri}, 1A, 3D, 3E)
$$C_{GL_n(q)}(D) = C_{GL_{wp}(q)}(D) \times GL_{t}(q) \leq L_1 \times GL_{t}(q),$$
and, 
$$N_{GL_n(q)}(D) = N_{GL_{wp}(q)}(D) \times GL_{t}(q) \leq
(GL_p(q)\wr \Sigma_w) \times GL_{t}(q). \Box$$
\end{lemma}

\begin{theorem} \label{BrPuig}
(M. Brou\'e, \cite{Bro2}, 3.5)  
The block $b_{\tau,w} ^{G_q}$ has $D$ as a defect group, and 
$$Br_D^{GL_n(q)}(b_{\tau,w} ^{G_q}) = 1_{kC_{GL_{wp}(q)}D} \otimes b_{\tau,0} ^{G_q},$$ 
where $Br_D^{GL_n(q)}$ is the classical Brauer morphism. 
$\Box$
\end{theorem}

We record some information on
these defect groups $D$ of unipotent blocks of $GL_n(q)$, and on
their
centralizers and normalizers.

\bigskip

Let $D_1$ be a Sylow $l$-subgroup of $GL_p(q)$.

\begin{lemma} (\cite{FoSri}, 3D, 3E)

(a) $\times^w D_1$ is a Sylow $l$-subgroup of $L_1$. 

\bigskip

(b) $D_1$ is isomorphic to a cyclic group of order $l^a$, 
where 
$$a = max \{ i \in \mathbb{Z}_{\geq 0} | l^a \textrm{ divides } q^p-1 \}.$$
 
(c) $D \cong \times^w D_1$ 
if, and only if, $w<l$.
If $w \geq l$, then $D$
is non-abelian.

\bigskip

(d) The normalizer of $D$ in $GL_{wp}(q)$ is 
contained in the subgroup,
$$L_1 \rtimes \Sigma_w \cong GL_p(q) \wr
\Sigma_w,$$ 
where $\Sigma_w$ is the group of block permutation matrices, whose
conjugation action
on $\times^w GL_p(q)$ permutes the $GL_p(q)$'s.
$\Box$
\end{lemma}

\begin{center}
{\large Brauer correspondence for 
blocks of symmetric groups.}
\end{center}

We now describe Brauer correspondence for blocks of symmetric groups.
Therefore, for the rest of this section, $l=p>0$.

\bigskip

Let $t = n-wp \geq 0$.
Let
$$Y_1 = \times^w \Sigma_p < \Sigma_{wp},$$ 
$$Y_2 = \Sigma_{wp} \times \Sigma_{t} < \Sigma_n$$ be Young subgroups. 
Let $D$ be a Sylow $p$-subgroup of $\Sigma_{wp}$.

\bigskip

Let $\tau$ be a partition of $t$, and a $p$-core.
Let $b_{\tau,w} ^\Sigma$ (respectively $b_{\tau,0} ^\Sigma$) 
be the block of $\Sigma_n$
(respectively $\Sigma_t)$ with $p$-core $\tau$.

\begin{lemma} (\cite{JamKerb}, 4.1.19,4.1.25)
$$C_{\Sigma_n}(D) = C_{\Sigma_{wp}}(D) \times \Sigma_{t} \leq
Y_1 \times \Sigma_{t},$$
and, 
$$N_{\Sigma_n}(D) = N_{\Sigma_{wp}}(D) \times \Sigma_{t} \leq
(\Sigma_p \wr \Sigma_w) \times \Sigma_{t}. \Box$$
\end{lemma}

\begin{theorem}\label{Brauermorph} (M. Brou\'e, L. Puig, \cite{Bro2}, 2.3)  
The block $b_{\tau,w} ^\Sigma$ has $D$ as a defect group, and 
$$Br_{D'}^{\Sigma_n}(b_{\tau,w} ^\Sigma) = 1_{kC_{\Sigma_{wp}}(D')} \otimes b_{\tau,0} ^\Sigma,$$
for any $p$-subgroup $D'$ of $\Sigma_{wp}$ without fixed points in the
set $\{1,2,...,wp \}$.
\end{theorem}

We offer some information on
these defect groups $D$, and their
centralizers and normalizers.

\bigskip

Let $D_1$ be a Sylow $p$-subgroup of $\Sigma_p$.
The following lemma is easily checked.

\begin{lemma}

(a) $\times^w D_1$ is a Sylow $p$-subgroup of $Y_1$

\bigskip

(b) $D_1 \cong C_p$, a cyclic group of order $p$.

\bigskip

(c) $D \cong \times^w D_1$ 
if, and only if, $w<p$.
If $w \geq p$, then $D$ is non-abelian.

\bigskip

(d) The normalizer of $D$ in $\Sigma_{wp}$ is 
contained in the subgroup
$$Y_1 \rtimes \Sigma_w \cong \Sigma_p \wr \Sigma_w,$$ 
where $\Sigma_w$ is the group of block permutation matrices, whose
conjugation action
on $\times^w \Sigma_p$ permutes the $\Sigma_p$'s.
$\Box$
\end{lemma}

\begin{center}
{\large Morita equivalence for symmetric algebras.}
\end{center}

\begin{definition}
An ${\cal O}$-order $A$, of finite rank, is said to be \emph{symmetric} 
if there exists an 
${\cal O}$-linear form 
$\phi: A \rightarrow {\cal O}$ such that,
\begin{description}
\item (i) $\phi(aa') = \phi(a'a)$ for all $a,a' \in A$.
\item (ii) The map 
$\hat{\phi}: A \rightarrow Hom_{\cal O}(A,{\cal O})$ defined by 
$\hat{\phi}(a)(a') = \phi(aa')$, for $a,a' \in A$, is an isomorphism of ${\cal
O}$-modules.
\end{description}
\end{definition}

M. Brou\'e has given a sufficient condition for two ${\cal O}$-algebras 
to be Morita equivalent:

\begin{theorem}\label{balls} (\cite{Bro3}, 2.4)
Let $A$ and $B$ be symmetric ${\cal O}$-orders of finite rank, and let $M$ be an
$A$-$B$ bimodule which is projective of finite rank
at the same time as a left $A$-module and
as a right $B$-module.
Suppose that the functor 
$$KM \otimes_{KB} - : KB-mod \rightarrow KA-mod,$$  
is an equivalence of categories.
Then the functor 
$$M \otimes_B - : B-mod \rightarrow A-mod,$$ 
is an equivalence of categories.
$\Box$
\end{theorem}

\newpage

\begin{center}
{\bf \large 
Chapter III

Rock blocks of finite general linear groups and Hecke algebras, when $w<l$.}
\end{center}

\bigskip

We sketch a proof of the 
Morita equivalence of Rock blocks of finite general linear
groups of abelian defect and their local blocks. 
Our proof of this theorem imitates Chuang and Kessar's proof
of theorem~\ref{structsymm}.
We subsequently deduce an analogous result for Rock blocks of 
Hecke algebras.

\begin{center}
{\large Rock blocks of finite general linear groups.}
\end{center}

Let $(K, {\cal O}, k)$ be an $l$-modular system.
Let $q$ be a prime power, coprime to $l$.
Let $p= p(G_q)$ be the multiplicative order of $q$, modulo $l$.

\begin{theorem} \label{main} (W. Turner \cite{Turner}, H. Miyachi \cite{Miyachi})
Let $w<l$.
The Rock block ${\cal O}b_{\rho,w} ^{G_q}$ is Morita equivalent 
to the principal block, ${\cal O}b_{\emptyset,1} ^{G_q} \wr \Sigma_w$, of 
$GL_p(q) \wr \Sigma_w$.
\end{theorem}

For pedagogical purposes, we sketch the proof of this result. 
This will allow us to understand the similarities and differences as we approach Rock blocks of nonabelian defect ($w \geq l$), in later chapters.
We include a proof of proposition \ref{rank}, since this is relevant for the proof of theorem  \ref{mainnew}.
Further details can be found elsewhere \cite{Turner}.

\bigskip

Let the Rouquier core $\rho = \rho(p,w)$ have size $r$. 
Let $v = wp +r$.  \index{$v$} \index{$r$} 

The principal block of $GL_p(q)\wr \Sigma_w$ is Morita equivalent to that
block tensored with the defect 0 block ${\cal O}{\bf B}_{\rho,0} ^{G_q}$, 
a block of $(GL_p(q) \wr \Sigma_w) \times
GL_r(q)$.
We prove theorem~\ref{main} by showing that Green correspondence
induces a Morita equivalence between this block of 
$(GL_p(q) \wr \Sigma_w) \times GL_r(q)$,
and the block ${\cal O}{\bf B}_{\rho,w}^{G_q}$ of $GL_v(q)$, when $w<l$.

\bigskip

Let $GL_v(q) = G = G_0 > G_1 > ... > G_w = L$ \index{$G$} \index{$G_i$} \index{$L$} 
be a sequence of Levi subgroups of $GL_v(q)$, where 
$$G_i = GL_p(q)^i \times GL_{v-ip}(q).$$ 
Let $$P_1 >...> P_w$$ be a sequence of parabolic subgroups
of $G_0 >...> G_{w-1}$ with Levi subgroups $$G_1 >...> G_w,$$ and unipotent
complements $$U_1 >... >U_w,$$
such that $P_i = U_i \rtimes G_i < G_{i-1}$.

\bigskip

Let $U_i^+ = \frac{1}{|U_i|} \sum_{u \epsilon U_i} u$, 
a central idempotent in ${\cal O} U_i$. 

Note that $|U_i|$ is a power of $q^p$, so equal to 1 
(modulo $l$). Note also that $G_i$ commutes with $U_i^+$. 

\bigskip

Let $a = b_{\emptyset,1}^{G_q}$ be the principal block idempotent
of $GL_p(q)$. Let
$$b_i = a_1 ^{\otimes i} \otimes b_{\rho,w-i} ^{G_q},$$  
\index{$b_i$} 
a block idempotent of $G_i$, for $1 \leq i \leq w$.
We set $G = G_0$, $b = b_0$, \index{$b$} and $f = b_w$. \index{$f$}

\bigskip

Let $\Sigma_w$ be the subgroup of permutation matrices of $GL_v(q)$ 
whose conjugation action 
on $L$ permutes the factors of 
$GL_p(q)^i$.

We define $N$ \index{$N$} to be the semi-direct product of $L$ and $\Sigma_w$, a 
subgroup of $GL_v(q)$ isomorphic to $(GL_p(q) \wr \Sigma_w) \times GL_r(q)$.

\bigskip

To prove theorem~\ref{main}, 
we show that ${\cal O}Nf$ and ${\cal O}Gb$ are
Morita equivalent, so long as $w<l$.
It is not clear how to define the corresponding 
${\cal O}Gb$-${\cal O}Nf$-bimodule directly.
However, we can describe the
${\cal O}Gb$-${\cal O}Lf$-bimodule obtained by restriction.

\bigskip

Let $Y = _GY_L = {\cal O} Gb_0U^+_1b_1...U^+_wb_w$, \index{$Y$}
an $( {\cal O} Gb, {\cal O} Lf)$-bimodule. 
The functor $Y \otimes_L -$
from $L$-mod to $G$-mod is 
$$HCInd_{G_1b_1}^{G_0b_o} ... HCInd_{G_wb_w}^{G_{w-1}b_{w-1}},$$
where $HCInd$ is Harish-Chandra induction.

\begin{proposition} \label{rank}
The algebras ${\cal O}Nf$, and $End_G(Y)$ have the same ${\cal O}$-rank,
equal to $w! dim_K(KLF)$.
\end{proposition}
Proof:

The algebra $End_G(Y)$ is ${\cal O}$-free.
It is therefore enough to compute the dimensions over $K$ of $End_G(KY)$ and $KNf$.
One of these is straightforward -
$_NKNf$ is the induced module $Ind_L^N(KLf)$, so certainly has
dimension $w! dim (KLf)$.
The proposition will be proved when we have shown that
$End_G(KY)$ has the same dimension.   

We calculate,
$$_G KY = _G KY \otimes_{L} KLf =
HCInd_{G_1b_1}^{G_0b_0} ... HCInd_{G_wb_w}^{G_{w-1}b_{w-1}}(_LKLf).$$ 
This is to be done by first computing
$HCInd_{G_1b_1}^{G_0b_0} ... HCInd_{G_wb_w}^{G_{w-1}b_{w-1}}(\psi)$, 
when $\psi$
is an irreducible character of $KLf$,
using the Littlewood-Richardson rule.
The relevant combinatorics have already been 
described by Chuang and Kessar, and
we record these below as lemma~\ref{Ch-Kess}.
Here, if $\lambda = (\lambda_1 \geq \lambda_2 \geq ...)$ and
$\mu = (\mu_1 \geq \mu_2 \geq...)$ are partitions,
we write $\mu \subset \lambda$ exactly when $\mu_i \leq \lambda_i$
for $i=1,2,...$.
An abacus is fixed so that the 
relevant $p$-quotients are well-defined.

\begin{lemma} \label{Ch-Kess}(J. Chuang, R. Kessar \cite{KessCh}, Lemma 4) 
Let $\lambda$ be a partition with
$p$-core $\rho$ and weight $u \leq w$. Let $\mu \subset \lambda$ be
a partition with $p$-core $\rho$ and weight $u-1$. 
Then
there exists $m$ with $0 \leq m \leq p-1$ such that $\mu^{l} =
\lambda^{l}$ for $l \neq m$ and $\mu^m \subset
\lambda^{m}$ with $|\mu^m| = |\lambda^m| -
1$. Moreover the complement of the Young diagram of $\mu$ in that of
$\lambda$ is the Young diagram of the hook partition $(m
+1,1^{p-m-1})$.    
$\Box$
\end{lemma}

In terms of character theory, by the Littlewood-Richardson rule, this
means that Harish-Chandra induction from 
$K{\bf B}_{\emptyset,1} ^{G_q} \otimes K{\bf B}_{\rho,u-1} ^{G_q}$ to
$K{\bf B}_{\rho,u} ^{G_q}$
takes the character $\chi(1,(m+1,1^{p-m-1})) \otimes \chi(1,\mu)$
to the sum of $\chi(1,\lambda)$'s, such that $\lambda$ is obtained from
$\mu$ by moving a bead down the $m^{th}$ runner.

Let us count the number of ways of sliding 
single beads down the $e^{th}$
runner of a core $j$ times, so that on the resulting runner
the lowest bead has been lowered
$\mu^{e}_1$ times, 
the second top bead has been
lowered $\mu^e_2$ times, etc., so that $\mu^e_1 \geq
\mu^e_2 \geq ...$ and $\sum_i \mu^e_i = j$.
It is equal to the number of 
ways of writing the numbers $1,...,j$ in the
Young diagram of $(\mu^e_1,\mu^e_2,...)$ so that numbers increase
across rows and down columns - 
that is, the degree of the character $\chi^{\mu^e}$ of the
symmetric group $\Sigma_j$ (see \cite{JamKerb}, 7.2.7).

\bigskip

The characters in the block $KLf$ are of the form,  
$$\psi = \chi(s_1,\lambda_{1}) 
\otimes ... \otimes \chi(s_w,\lambda_{w}) \otimes \chi(1,\rho),$$
where either $s_i=1$, and $\lambda_i$ is a $p$-core,
or else $\lambda_i = (1)$, and $s_i = (\sigma_i)$ 
is given by a degree $p$ element $\sigma_i \in \bar{\mathbb{F}}_q^*$,
whose order is a power of $l$, for $1 \leq i \leq w$.

\bigskip

A combinatorial description of the multiplicity of a
given irreducible
character of $KGb$ in $KY \otimes_{L} \psi$
is now visible: 

Suppose that $s_1,...,s_{r_0}$ are all equal to $1$,
that $\lambda_i = (1)$ for $i \geq r_0+1$,
that $\sigma_{r_0+1} \sim ... \sim \sigma_{r_0+r_1} =:\theta_1$ 
are conjugate 
elements of $\bar{\mathbb{F}}_q^*$, that
$\sigma_{r_0+r_1+1} \sim ... \sim \sigma_{r_0+r_1+r_2}=:\theta_2$ are
are conjugate 
elements of $\bar{\mathbb{F}}_q^*$,
not conjugate to $\theta_1$, etc. etc.
Also suppose that $\lambda_i$, for $i=1,...,r_0$ is a $p$-hook, and that,
$$\lambda_1=...=\lambda_{e_0}=(1^p),$$
$$\lambda_{e_0+1}=...=\lambda_{e_0+e_1}=(2,1^{p-1}),$$
$$...$$ 
$$\lambda_{e_0+..+e_{p-2}+1}=...=\lambda_{e_0+..+e_{p-1}}=(p),$$ 
where $e_0+..+e_{p-1} = r_0$.
Then $KY \otimes_{L} \psi$ is
equal to the character sum, 
$$\sum
( dim \chi^{\mu^0}... dim \chi^{\mu^{p-1}}
. dim \chi^{\nu^1}. dim \chi^{\nu^2}...)$$ 
$$\chi \left( 
(1, (\theta_1)^{r_1}, (\theta_2)^{r_2},...),(\mu,\nu^1,\nu^2,...)
\right)$$
Here, the summation is over partitions
$\mu=[\mu^0,..,\mu^{p-1}]$ of $|\rho|+r_0p$ with
core $\rho$, such that $(|\mu^0|,...,|\mu^{p-1}|)=
(e_0,e_1,...,e_{p-1})$;
over partitions $\nu^1$ of $r_1$;
over partitions $\nu^2$ of $r_2$, etc.

If, when we selected a character of $L$, we had permuted some
of the $(s_i,\lambda_{i})$'s 
(there are $(w!/e_0!..e_{p-1}!r_1!r_2!...)$ different arrangements),
we would have seen the same character when we applied $_G KY \otimes_L -$.
So the character of $_G KY$ is the sum of characters,

$$\sum_{|\mu^i|=e_i,|\nu^i|=r_i} \bigg[ (w!/e_0!..e_{p-1}!r_1!r_2!...)$$
$$\times
dim \chi^{\mu^0}... dim\chi^{\mu^{p-1}}.
dim \chi^{\nu^1}. dim \chi^{\nu^2}...$$
$$\times
dim \chi(1,(1^p))^{e_0}.
dim \chi(1,(2,1^{p-2}))^{e_1}
... dim \chi(1,(p))^{e_{p-1}}$$
$$\times 
dim  \chi(\rho_{X-1}) \times 
dim \chi({\theta_1}, (1))^{r_1}.
dim \chi({\theta_2}, (1))^{r_2}
...$$
$$\times
\chi \left( (1, (\theta_1)^r_1, (\theta_2)^{r_2},...),
([\mu^0,...,\mu^{p-1}],\nu^1,\nu^2,...) \right)
\bigg].$$ 
What is the dimension of the semisimple algebra 
$End_G(KY)$ ? 
Remembering that
$\sum_{|\sigma|=m}$ $|\chi^{\sigma}|^2 = m!$, it is,
$$\sum_{e_0+..+e_{p-1}+r_1+r_2+...=w} \bigg[
 (w!/e_0!..e_{p-1}!r_1!r_2!...)^2$$ 
$$\times e_0!..e_{p-1}!r_1!r_2!...$$
$$\times
\left( dim \chi(1,(1^p)) \right) ^{2e_1}.
\left( dim \chi(1,(2,1^{p-2})) \right) ^{2e_2}
... \left( dim \chi(1,(p)) \right) ^{2e_p}$$
$$\times
\left( dim \chi(\rho_{X-1}) \right) ^2 \times
\left( dim \chi({\theta_1},(1)) \right) ^{2r_1}.
\left( dim \chi({\theta_2},(1)) \right)^{2r_2}
... \bigg]$$

$$= w! \sum_{e_0+..+e_{p-1}+r_1+r_2+...=w}
\bigg[ (w!/e_0!..e_{p-1}!r_1!r_2!...)$$              
$$\times 
\left( dim \chi(1,(1^p)) \right) ^{2e_0}.
\left( dim\chi(1,(2,1^{p-2})) \right) ^{2e_1}
...\left( dim \chi(1,(p)) \right) ^{2e_{p-1}}$$
$$\times
\left( dim \chi(\rho_{X-1}) \right) ^2 \times
\left( dim \chi({\theta_1},(1)) \right) ^{2r_1}.
\left( dim \chi({\theta_2},(1)) \right) ^{2r_2}...\bigg]$$

$$= w!.\textrm{dim}(KLf). \Box$$

\bigskip

Let $D = D^1 \times ...\times D^w$ be a Sylow $l$-subgroup of
$GL_p(q)^1 \times ... \times GL_p(q)^w$.

\begin{lemma} \label{one} Let $w<l$.

(a) $D$ is a defect group of ${\cal O} G_ib_i$, for $i=0,...,p-1$.

\bigskip

(b) $Br^G_D(b_i) = 1_{kD} \otimes b_{\rho,0}^{G_q}$, and $Br^G_D(U_i^+) = 1$.

\bigskip

(c) $N$ stabilizes $f$, and as an ${\cal O} (N \times L)$-module, ${\cal
    O} Nf$ is
    indecomposable with diagonal vertex $\Delta D$ \index{$\Delta D$}. 
In particular, ${\cal O} Nf$ is a
    block of $N$.

\bigskip

(d) ${\cal O} Gb$ and ${\cal O} Nf$ both have defect group $D$, and are Brauer
    correspondents.  $\Box$
\end{lemma}

\bigskip

By the Brauer correspondence, 
the ${\cal O}(G \times
G)$-module ${\cal O} Gb$ and the ${\cal O} (N \times N)$-module 
${\cal O} Nf$ both have vertex
$\Delta D$ and are Green correspondents. 
Let $X$ \index{$X$} be the Green correspondent of ${\cal O} Gb$ in $G \times N$, an
indecomposable summand of $Res^{G \times G}_{N \times N}
( {\cal O} Gb)$ with
vertex $\Delta D$.  
Because ${\cal O}Nf$ is a direct summand of $Res^{G \times N}_{N \times
N}(X)$, we have $Xf \neq 0$, so $Xf = X$ and $X$ is an $( {\cal O} Gb,
{\cal O} Nf)$-bimodule.

Theorem~\ref{main} is a consequence of the following:

\begin{proposition}
Let $w<l$. Then $_GY_L \cong {_G X_L}$, and ${\cal O}Nf \cong End_G(X) \cong End_G(Y)$ as algebras.
The left ${\cal O}G$-module $_GX$ is a progenerator for ${\cal O}Gb$.
Hence $X \otimes_N -$ induces a
Morita equivalence between ${\cal O} Nf$ and ${\cal O} Gb$. $\Box$
\end{proposition}

\begin{remark}
The correspondence between indecomposable modules of $kNf$ and
indecomposable modules of $kGb$ given by Theorem~\ref{main} 
above is exactly
Green correspondence between $G$ 
and $N$.
For if $M$ is an indecomposable of $kNf$ with vertex $Q$, the
$kG$-module $X \otimes_{kN} M$ cannot have a smaller vertex than $Q$,
as then $M = X^* \otimes_{kG} X \otimes_{kN} M$ 
would have a smaller vertex than
$Q$.
\end{remark}

\begin{remark} The proof of theorem~\ref{main} fails when $w \geq l$, because
the indecomposable module $_G X_N$ does not restrict to an 
indecomposable module, $_G X_L$.
\end{remark} 

\begin{center}
{\large Rock blocks of Hecke algebras, when $w<l$.}
\end{center}

Let us persist with the notation of the last section. Thus,  
$(K, {\cal O}, k)$ is an $l$-modular system, and
$q$ is a prime power, coprime to $l$.

\begin{lemma} \label{referee}
Suppose $A$ and $B$ are ${\cal O}$-algebras, which are free ${\cal O}$-modules of finite rank, and
suppose $A$ and $B$ are Morita eqivalent via $F: A-mod \rightarrow B-mod$.
Let $e,f$ be idempotents in $A,B$. Then the following are equivalent:

(i) $eAe$ and $fBf$ are Morita equivalent.

(ii) $eS = 0$ if, and only if, $fF(S) = 0$, for all irreducible $kA$-modules $S$.  
\end{lemma}

Theorem~\ref{main} squashes to Hecke algebras as follows:

\begin{theorem} \label{mainnew} Let $w<l$.
Let $p = p(G_q) = p({\cal H}_q)>1$. 
The Rock block ${\cal O}{\bf B}_{\rho,w} ^{{\cal H}_q}$ is
Morita equivalent to ${\cal O}{\bf B}_{\emptyset,1} ^{{\cal H}_q} \wr \Sigma_w$.
\end{theorem}
Proof:

Let $B_n(q)$ be a Borel subgroup of $GL_n(q)$ for $n=p,v$.
Since $p>1$, for $n=p,v$, the ${\cal O}GL_n(q)$-module
${\cal O}[GL_n(q)/B_n(q)]$ is projective,
having been
induced from the projective trivial ${\cal O}B_n(q)$-module.

\bigskip

We now construct a Morita 
bimodule from bimodules already available to us.  
Let $G=GL_v(q)$ be as in theorem~\ref{main}.
Let $H_w = \times^w GL_p(q)$.
Let $G_w = H_w \rtimes \Sigma_w = GL_p(q) \wr \Sigma_w$.
Let $_GZ_{G_w}$ \index{$Z$} be a bimodule
inducing a Morita
equivalence between ${\cal O}{\bf B}_{\rho,w} ^{G_q}$ and 
${\cal O}{\bf B}_{\emptyset,1} ^{G_q} \wr \Sigma_w$; 
such a bimodule exists by theorem~\ref{main}.
Let 
$$\eta = \frac{1}{|B_p(q)|^w} \sum_{x \in B_p(q) \times .. \times B_p(q)} x,$$ 
$$\xi = \frac{1}{|B_v(q)|} \sum_{b \in B_v(q)} b.$$ 
Then  
$$\xi Z \eta$$
is an ${\cal H}_q(\Sigma_v)$-${\cal H}_q(\Sigma_p) \wr \Sigma_w$ bimodule.  

\bigskip

We need to show
that under the Morita equivalence from 
${\cal O}{\bf B}_{\rho,w} ^{G_q}$ to 
${\cal O}{\bf B}_{\emptyset,1} ^{G_q} \wr \Sigma_w$;
given by $Z$, simple $kG$-modules killed by $\eta$
become simple $kG_w$-modules killed by $\xi$, and vice-versa.
The truth of the theorem then follows from lemma \ref{referee}. 

\bigskip

The characters of $G$ killed by $\xi$ are the non-unipotent
characters, by theorem~\ref{Swd}.
Under $_{H_w}Z^* \otimes_G -$, these become sums of characters 
$\chi^1 \otimes ...\otimes \chi^w$, such that one of the $\chi^i$'s is a
non-unipotent character of $GL_p(q)$ 
(according to the proof of proposition~\ref{rank}).
These are all killed by $\eta$.
Conversely, by theorem~\ref{Kill},
the characters of $G_w$ killed by $\eta$ are all
composition factors of characters 
$Ind_{H_w}^{G_w}(\chi^1 \otimes...\otimes \chi^w)$, where one of the
$\chi^i$'s is a non-unipotent character of $GL_p(q)$.
And $_G Z \otimes_{G_w} -$ sends this induced character to 
a sum of non-unipotent characters of $G$, all of which are
killed by $\xi$.
  
\bigskip

The simples for $G$ which vanish under multiplication 
by $\xi$ are those
$D(\lambda)$ indexed by
$p$-singular partitions.
These are simple composition factors 
of non-unipotent characters $\chi((1, (\sigma)^{|\nu|}),(\lambda_{0}, \nu))$,
where $\nu$ is a non-empty partition, and $\sigma$ is an $l$-element 
of $\bar{\mathbb{F}}_q ^*$ of degree $p$ (\cite{BDK}, Theorem 4.4d).
But these characters correspond (under the Morita equivalence) to
characters for $G_w$ which are killed by $\eta$. 
Conversely, simple modules for $G_w$ which are killed by $\eta$
are composition factors of induced characters
$Ind_{H_w}^{G_w}(\chi^1 \otimes...\otimes \chi^w)$, where one of the
$\chi^i$'s is a non-unipotent character of $GL_p(q)$.
This becomes a character of $G$ sent to zero by $\xi$
on application of $Z
\otimes_{G_w} -$.
$\Box$

\newpage
\begin{center}
{\bf \large 
Chapter IV

Rock blocks of symmetric groups, and the Brauer morphism.}
\end{center}

\bigskip

The structure theorem of Chuang and Kessar states that Rock blocks of symmetric groups of abelian defect $w$ are Morita equivalent to 
$k{\bf B}_{\emptyset,1}^\Sigma \wr \Sigma_w$ (theorem \ref{structsymm}). 
The corresponding statement is false in nonabelian defect: the global and local blocks have different numbers of simple modules.
Furthermore, the techniques developed to study Brou\'e's abelian defect group conjecture \cite{Broue}
give little clue as to how to formulate an analogous result in nonabelian defect, let alone how to prove it.

Alperin's weight conjecture \cite{Alperin}
suggests a deep uniformity in representation theory, which exists for all blocks of finite groups, and not only those of abelian defect.
It therefore makes sense to search for an analogue of Chuang and Kessar's result, which is true for blocks arbitrary defect,
and to develop techniques which may lead towards a proof.
That is the dominant concern of this monograph. 
In chapter $8$, we describe a conjecture, which predicts the structure of Rock blocks of symmetric groups of arbitrary defect.
In other chapters, we give various numerical and structural results which point towards this conjecture, although none of them confirm it.

From a finite group representation theoretical perspective, our methods are somewhat eccentric, 
involving a peculiar application of the Brauer homomorphism, the theory of quasi-hereditary algebras, 
Ringel duality, cross-characteristic comparisons, quivers, Schur bialgebras, doubles, etc..
Standard tools, such as Green correspondence, appear to be all but useless in our situation.
The resulting conjecture (conjecture \ref{Rock}) is simple in essence, and appears to be quite deep.
It can be seen to be true in abelian defect, by comparison with Chuang and Kessar's structure theorem.

In this chapter, we introduce notation, for the study of Rock blocks of symmetric groups,
of arbitrary defect. We show that in characteristic two,
a Rock block is isomorphic to the group algebra of $\Sigma_2 \wr \Sigma_w$, 
once it has been cut by a certain idempotent (theorem~\ref{viva}).
Our proof involves an unusual application of the Brauer homomorphism.
We use this method to obtain a weaker result for Rock blocks of symmetric groups
in arbitrary characteristic: we prove that the endomorphism ring of a certain 
$l$-permutation module $M$ for $k{\bf B}_{\rho,w} ^\Sigma$ is Morita equivalent to $k\Sigma_w$
(theorem~\ref{Wick}).

\begin{center}
{\large Notation for Rock blocks of symmetric groups.}
\end{center}

Throughout this chapter, we consider blocks of symmetric groups.
Therefore, $l=p$.

Let $(K, {\cal O}, k)$ be a $p$-modular system.
Let $w$ be a natural number.
Let $\rho = \rho(p,w)$ be a Rouquier core, and a partition of $r$. 
Let $v=wp+r$.

Let $\Sigma_{v} = Sym\{1,2,...,v \}$.
Let $b_{\rho,w}^\Sigma$ be the Rock
block of $\Sigma_v$.
Let
$$L = \Sigma_p^1 \times ... \times \Sigma_p^w \times \Sigma_r^0 \leq \Sigma_v,$$
\index{$L$}
where $\Sigma_p^i = Sym \{(i-1)p+1,...,ip \}$, 
and $\Sigma_r^0 = Sym \{wp+1,...,wp+r \}$.
Let 
$$f = b_{\emptyset,1} ^\Sigma \otimes ... \otimes b_{\emptyset,1} ^\Sigma 
\otimes b_{\rho,0} ^\Sigma,$$ \index{$f$}
a block of $L$.
Let $$D = C_p^1 \times ... \times C_p^w \leq L,$$ \index{$D$}
where $C_p^i$ is the group of order $p$ generated by the single 
element $((i-1)p+1 ... ip)$.
Let
$$N = L \rtimes \Sigma_w \cong (\Sigma_p \wr \Sigma_w) \times \Sigma_r^0,$$ \index{$N$}
the normaliser of $L$ in $\Sigma_v$.
Since the centralizer $C_N(D) \leq L$,  
the idempotent $f$ is also a block of $N$ (\cite{Alp} 15.1(5)). 

\bigskip

Let $e$ \index{$e$} be an idempotent of $k\Sigma_v$,  
defined to the product $e_w.e_{w-1}...e_0$ of block idempotents, 
$$e_i = b_{\emptyset,1} ^\Sigma \otimes ... \otimes b_{\emptyset,1} ^\Sigma
\otimes b_{\rho,i}^\Sigma,$$
of $\Sigma_p^1 \times .. \times \Sigma_p^{w-i} \times \Sigma_{ip+r}$.

\bigskip

Aping proposition~\ref{rank}, and its proof, we have,

\begin{proposition} \label{should}
The dimensions $dim_k(ek\Sigma_ve)$ and $dim_k(kNf)$ are both
equal to $w!.dim_k(kLf)$.
$\Box$
\end{proposition}

\begin{conjecture} \label{cutiso}
There is an algebra isomorphism 
$ek\Sigma_ve \cong kNf$, for arbitrary $w,p$.
\end{conjecture}

\begin{remark}
The conjecture is true for $w<p$, by Chuang and Kessar \cite{KessCh}. 
The method of Chuang and Kessar fails for $w>p$, because the $k\Sigma_v$-$kNf$ bimodule $k\Sigma_v e$ fails to have diagonal vertex in this case, 
and consequently the Green correspondence cannot be used to give an abstract description of the bimodule.
We prove the conjecture for $p=2$ below, using the Brauer homomorphism.

R. Paget \cite{Paget2} has computed the 
projective summands of $_{\Sigma_v} k\Sigma_ve$, and
shown that $ek\Sigma_ve$, and $kNf$ have the same decomposition matrix. 
Her proof uses theorem~\ref{nilpotent}, of this monograph.
\end{remark}

\begin{center}
{\large Endomorphism rings.}
\end{center}

Let $P$ be the subgroup of $D$ of order $p$ 
generated by the element,
$$z = (1 ... p)(p+1 ... 2p)...((w-1)p+1 ... wp).$$
Let, $$C=C_{\Sigma_v}(P)=D \rtimes \Sigma_w \times \Sigma_r^0 \leq N.$$ \index{$C$}

We consider the classical Brauer homomorphism - the surjective  
algebra homomorphism from
$(k\Sigma_v)^P$ to $kC$ which truncates elements of the group algebra
at $C$.
The images of block idempotents under the Brauer morphism 
are given by theorem~\ref{Brauermorph}.
In particular, 
$$Br^{\Sigma_v}_D(e)=
Br^{\Sigma_v}_P(e) = Br^{\Sigma_v}_P(f) = 
1_{kD} \otimes b_{\rho,0}^\Sigma.$$
Therefore,
$$Br^{\Sigma_v}_P((ek\Sigma_ve)^P) = 
Br^{\Sigma_v}_P(e).kC.Br^{\Sigma_v}_P(e) = kC.b_{\rho,0}^\Sigma.$$
If $p=2$, then $N=C$, and $f = Br^{\Sigma_v}_P(f) = Br^N_P(f)$. 
In this case, we can prove conjecture \ref{cutiso}:

\begin{theorem} \label{viva}
If $p=2$, then the Brauer homomorphism 
$$Br^{\Sigma_v}_P: (ek\Sigma_ve)^P \rightarrow kC$$  
restricts to an isomorphism of algebras $(ek\Sigma_ve)^P \cong kNf$.
Furthermore, $(ek\Sigma_ve)^P = ek\Sigma_ve$. 
\end{theorem}
Proof:

From lemma~\ref{should}, 
we know that $dim_k(ek\Sigma_ve) = dim_k(kNf)$.
Also $N=C$, so that 
$$dim_k(ek\Sigma_ve)^P \leq dim_k(ek\Sigma_ve) = dim_k(kNf) = 
dim_k(kCf).$$
But the Brauer homomorphism
$Br^{\Sigma_v}_P: (ek\Sigma_ve)^P \rightarrow kCf$ is a surjection.
So the dimensions above are all equal, and the theorem is proven.
$\Box$

\begin{remark}
Theorem \ref{viva} is false for $p>2$, because $kNf$ and  $(ek\Sigma_ve)^P$ are both strictly larger than $kCf$.
\end{remark}

If we cut down the module $k\Sigma_ve$, and the algebras that we are
considering,
it is possible to generalise the proof of theorem~\ref{viva} 
to $p \geq 2$.
More precisely, let $\zeta_w = \sum_{x \in \Sigma_p^1 \times ... \times
\Sigma_p^w}x$. \index{$\zeta_w$}

\bigskip

Let $M = {\cal O} \Sigma_v e \zeta_w$. \index{$M$}
Note that $M$ is not a projective ${\cal O} \Sigma_v$-module.
Let ${\cal E} = End_{\Sigma_v}(M)$. \index{${\cal E}$}

\bigskip

Since $M$ is a $p$-permutation module, it is ${\cal O}$-free, and 
$End_{\Sigma_v}(RM) \cong R{\cal E}$, for
$R \in \{ K, {\cal O}, k \}$
(\cite{Land}, 1.14.5).

\bigskip

To complete this chapter, we present the following theorem, 
as the consequence of a triad of lemmas. 
We expect the theorem to be true for $R = {\cal O}$, but we are unable to prove it, since our method involves the Brauer morphism, 
which takes values in a vector space over $k$.

\begin{theorem} \label{Wick}
Let $R \in \{ K,k \}$. Then,
$$R{\cal E} \cong R\Sigma_w \otimes R{\bf B}_{\rho,0}^\Sigma.$$   
\end{theorem}

First, some notation.
Let $\Theta_u$ \index{$\Theta_u$} 
be the set of partitions with core $\rho$ of weight $u$
obtained by moving beads only up the rightmost runner of the abacus 
representation 
of $\rho$.
These correspond to partitions of $u$, via
$[\emptyset,..,\emptyset,\nu] \leftrightarrow \nu$.
In the usual notation for partitions, this correspondence is 
$$(\rho_1 + p\nu_1,\rho_2 + p\nu_2,...,\rho_u +p\nu_u,
\rho_{u+1},...) \leftrightarrow (\nu_1,\nu_2,...,\nu_u).$$

\begin{lemma} \label{comb1}
The character of $_{\Sigma_v} KM$ is equal to 
$$dim(\chi^\rho).\sum_{\nu} dim
(\chi^\nu). \chi^{[\emptyset,..,\emptyset,\nu]}.$$
Its endomorphism ring is isomorphic to 
$K\Sigma_w \otimes K{\bf B}_{\rho,0}^\Sigma$.
\end{lemma}
Proof:

We have $_{\Sigma_v} KM \cong K\Sigma_v e 
\otimes_{L} KLf \zeta_w$. 
This module is the $KLf$-module with character 
$dim(\chi^\rho). (\chi^{(p)} \otimes...\otimes \chi^{(p)} \otimes \chi^\rho)$,  
induced to 
$(\bigotimes^{w-1} K\Sigma_p) \otimes K{\bf B}_{\rho,1}^\Sigma$, then induced to 
$(\bigotimes^{w-2} K\Sigma_p) \otimes K{\bf B}_{\rho,2}^\Sigma$, etc., etc..
That its character is 
as stated follows from a symmetric group analogue  
of proposition~\ref{rank}, and its proof.
 
The endomorphism ring is isomorphic to
$K\Sigma_w \otimes K{\bf B}_{\rho,0}^\Sigma$, since all algebras concerned are
semisimple. 
$\Box$

\bigskip

From now on, let us fix
an isomorphism between the 
endomorphism ring of $_{\Sigma_v} KM$ and
$K\Sigma_w \otimes K{\bf B}_{\rho,0}^\Sigma$,  
so that under the Morita 
equivalence, given by $KM$,
between $K{\bf B}_{\rho,w}^\Sigma /Ann(KM)$ and $K\Sigma_w \otimes
K{\bf B}_{\rho,0}^\Sigma$,  
the characters $\chi^{[\emptyset,...,\emptyset,\lambda]}$ and
$\chi^\lambda \otimes \chi^\rho$ correspond.

\bigskip

Let $\zeta_D = \sum_{x \in D} x$. \index{$\zeta_D$}
The sum $\zeta_D$ is the image of $\zeta_w$ under $Br_P$.
Note that $C$ commutes with $\zeta_D$, 
and so $kC$ acts on the right of 
$kC.b_{\rho,0} ^\Sigma.\zeta_D$. Under this action, we have:

\begin{lemma} \label{Karin2}
We have an isomorphism of algebras, 
$$kC.b_{\rho,0} ^\Sigma
/Ann(kC.b_{\rho,0} ^\Sigma.\zeta_{D}) \cong k\Sigma_w \otimes k{\bf B}_{\rho,0} ^\Sigma.$$
\end{lemma}
Proof:

It is true that $C \cong (D \rtimes \Sigma_w) \times \Sigma_r^0$, so that the
top subgroup $\Sigma_w$ commutes with $\zeta_D$.
Thus,  $$k\Sigma_w \otimes k{\bf B}_{\rho,0} ^\Sigma \cong 
kC.b_{\rho,0} ^\Sigma.\zeta_D$$ 
as vector spaces, via $x \mapsto x \zeta_D$.
The annihilator
$Ann_{kC}(kC.b_{\rho,0} ^\Sigma.\zeta_D)$ is thus equal to
$I(D).k\Sigma_w.b_{\rho,0} ^\Sigma$ as a right $kC$-module, 
where $\Sigma_w$ is the top group, 
and where $I(D)$ is the augmentation ideal of $kD$.
Furthermore, the quotient of $kC.b_{\rho,0} ^\Sigma$ by the annihilator
acts on $k\Sigma_w \otimes k{\bf B}_{\rho,0} ^\Sigma$ by the
right regular action, and thus the quotient is actually
isomorphic to $k\Sigma_w \otimes k{\bf B}_{\rho,0} ^\Sigma$, as an algebra. 
$\Box$

\bigskip

We may prove a version of lemma~\ref{comb1} over $k$, using $Br_P$:

\begin{lemma} We have an isomorphism of algebras,
$$k{\cal E} \cong k\Sigma_w \otimes
k{\bf B}_{\rho,0}^\Sigma.$$ 
The implicit action of $k\Sigma_w \otimes k{\bf B}_{\rho,0}^\Sigma$ on $kM$ is,
$$m \circ (x \otimes y) = m.exye,$$
for $x \in k\Sigma_w$, $y \in k{\bf B}_{\rho,0} ^\Sigma$.   
\end{lemma}

Proof:

First observe that $k{\bf B}_{\rho,w} ^\Sigma =$ $ek\Sigma_ve\oplus *$ as $L$-$L$-bimodules, 
where $*$ has no summand with
vertex containing $\Delta P$,
since, 
$$ek\Sigma_ve(\Delta P) = Br_P(e).kC.Br_P(e) =$$
$$kC.b_{\rho,0}^\Sigma 
= k{\bf B}_{\rho,w} ^\Sigma(\Delta P).$$
{\noindent
Likewise, $kNf$ is a summand of $k{\bf B}_{\rho,w} ^\Sigma$ as $L$-$L$-bimodules,
where the complement}
has summands whose vertices do not contain $\Delta P$.
However, a vertex of the $L$-$L$-bimodule $\sigma kLf$ is 
$\Delta D^{(\sigma,1)} \geq \Delta P$ (for
$\sigma \in \Sigma_w$). 
So, as an $L$-$L$-bimodule, $kNf$ is the sum of summands of $k{\bf B}_{\rho,w } ^\Sigma$ 
with vertices containing $\Delta P$.
Hence, $ek\Sigma_ve$ is a summand of $kNf$ as an $L$-$L$-bimodule, and since 
their dimensions are equal by lemma~\ref{should}, there are
isomorphisms of $L$-$L$-bimodules, 
$$ek\Sigma_ve \cong kNf 
\cong \bigoplus_{\sigma \in \Sigma_w} (\sigma kLf).$$
It follows that $dim_k(ek\Sigma_ve \zeta_w) = w!.dim_k(k{\bf B}_{\rho,0} ^\Sigma)$.

\bigskip

Let $A$ be the subalgebra of $(ek\Sigma_ve)^P$ generated by $ekCe$.
Note that $A$ acts on the right of $kM = k\Sigma_ve \zeta_w$ 
by multiplication, thus commuting 
with the left action of $k\Sigma_v$.
In fact, $A$ acts on the right of $(ek\Sigma_ve)^P \zeta_w$ by
multiplication. 
Applying the Brauer homomorphism $Br_P$, we realise
an action of $Br_P(A) = kC.b_{\rho,0} ^\Sigma$ 
on the right of 
$$Br_P((ek\Sigma_ve)^P \zeta_w) =
Br_P((ek\Sigma_ve)^P).Br_P(\zeta_w) =
kC.b_{\rho,0} ^\Sigma \zeta_D.$$
by multiplication.
Thus there is an algebra homomorphism, given by the composition,
$$A \rightarrow Br_P(A) \rightarrow 
kC.b_{\rho,0} ^\Sigma /Ann(kC.b_{\rho,0} ^\Sigma.\zeta_{D}),$$
{\noindent in the kernel of which lies $Ann_A((ek\Sigma_ve)^P \zeta_w)$.  
Lemma~\ref{Karin2} implies the existence of a surjection,}
$$A/Ann_A((ek\Sigma_ve)^P \zeta_w) \rightarrow k\Sigma_w\otimes
k{\bf B}_{\rho,0} ^\Sigma,$$ 
{\noindent as well as the natural surjection, }
$$A/Ann_A(k\Sigma_ve \zeta_w) \rightarrow
A/Ann_A((ek\Sigma_ve)^P \zeta_w).$$
{\noindent In addition, there is a sequence of natural injections,} 
$$A/Ann_A(k\Sigma_ve \zeta_w) \rightarrow
End_{\Sigma_v}(k\Sigma_ve \zeta_w) \rightarrow$$
$$Hom_{\Sigma_v}(k\Sigma_ve,k\Sigma_ve \zeta_w) \rightarrow
ek\Sigma_ve \zeta_w.$$
{\noindent But we have already agreed that $ek\Sigma_ve \zeta_w$ 
has dimension
equal to,}
$$w!.dim_k k{\bf B}_{\rho,0} ^\Sigma =
dim_k (k\Sigma_w \otimes k{\bf B}_{\rho,0} ^\Sigma),$$
{\noindent
so in fact all of the above homomorphisms are isomorphisms.
In particular,  
$$End_{\Sigma_v}(kM) =
End_{\Sigma_v}(k\Sigma_ve \zeta_w) \cong k\Sigma_w \otimes k{\bf B}_{\rho,0} ^\Sigma.
\Box$$

\newpage
\begin{center}
{\bf \large 
Chapter V

Schur-Weyl duality inside Rock blocks of symmetric groups.}
\end{center}

In this chapter, we strengthen theorem~\ref{viva}, and show that the Brauer
homomorphism reveals Schur-Weyl duality between ${\cal S}(w,w)$, and 
$k\Sigma_w$,
structurally inside the
Rock block $k{\bf B}_{\rho,w} ^\Sigma$ (theorem~\ref{Comeon}).

An alternative proof of the existence of a quotient of
a symmetric group block, Morita equivalent to ${\cal S}(w,w)$ 
has been given by Cline,
Parshall, and Scott, using Steinberg's tensor product theorem 
for the algebraic group $GL_n$ (\cite{CPS4}, 5.3).
If the Reader inclines towards an understanding of blocks of finite groups, 
the proof given here should be interesting, 
because it is independent 
of algebraic group theory.

\begin{center}
{\large A Schur algebra quotient.}
\end{center}

We assume the notation of chapter 4.
Let $I = Ann_{Rb_{\rho,w} ^\Sigma}(M)$, \index{$I$}
and let $B_{\rho,w}^\Sigma = {\bf B}_{\rho,w} ^\Sigma/I$.
Then $I$ is ${\cal O}$-pure in ${\bf B}_{\rho,w} ^\Sigma$, and so $RB_{\rho,w} ^\Sigma = R{\bf B}_{\rho,w} ^\Sigma/RI$,
for $R \in \{ K ,{\cal O}, k \}$.

\begin{theorem} \label{Comeon}
Let $R \in \{ K,k \}$. 
Then $R{\cal E} = End_{\Sigma_v}(RM)$ 
is Morita equivalent to $R\Sigma_w$, and
$RB_{\rho,w}^\Sigma$ is Morita equivalent to the Schur algebra
${\cal S}(w,w)$, defined over $R$.
  
The ${\cal S}(w,w)$-$R\Sigma_w$-bimodule corresponding to the 
$RB_{\rho,w} ^\Sigma$-$R{\cal E}$-bimodule $RM$
is twisted tensor space, $E^{\otimes r \#}$.
\end{theorem}

Theorem~\ref{Comeon} may be seen as a module theoretic 
interpretation of a theorem of Erdmann \cite{Erd2}, 
which realises the
decomposition matrix of the Schur algebra ${\cal S}(w,w)$ 
as a submatrix of the
decomposition matrix of $k{\bf B}_{\rho,w} ^\Sigma$.
Indeed, we have the following interesting corollary:

\begin{corollary} (``Converse to Schur-Weyl duality'')
Every block of polynomial $GL_n(k)$-modules is Morita equivalent
to a quotient of some symmetric group algebra, localised at an
idempotent.
$\Box$
\end{corollary}

Theorem \ref{Comeon} is clearly
analogous to theorem \ref{Swd}
concerning general linear groups in non-describing
characteristic.
However, our methods hardly
resemble those of Brundan, Dipper, and Kleshchev. 

\bigskip

The proof of theorem~\ref{Comeon} is the length of this chapter. 
We first find summands of $kM$ as a right 
$k\Sigma_w \times \Sigma_r^0$-module
which are twisted Young modules for $k\Sigma_w$, tensored with 
the block $k{\bf B}_{\rho,0} ^\Sigma$. 
We then observe that all the indecomposable summands of $kM$ are isomorphic to such twisted Young modules.
We finally show that the map 
$k{\bf B}_{\rho,w}^\Sigma \rightarrow End_{k{\cal E}}(kM)$ is 
surjective by an explicit calculation.

\bigskip

Let us state and prove some combinatorial preliminaries.

\begin{lemma} \label{comb2}
Let $\mu \notin \Theta_w$ be a partition of weight $w$ with core
$\rho$. Let $\theta \in \Theta_w$.
Then $\mu \ntrianglerighteq \theta$.
\end{lemma}
Proof:
In an abacus representation of $\mu$ obtained from 
$\rho$, at least one bead must be moved up a runner other than the
rightmost runner.
This means that $\mu_i$ is greater than $\rho_i$, for some $i>w$.  
Thus, $\mu_1 + \mu_2 + ... + \mu_w$ is at most
$\rho_1 + \rho_2 + ... +\rho_w + (wp-1)$.
But $\theta_1 + \theta_2 + ... + \theta_w = 
\rho_1 + \rho_2 + ... +\rho_w + wp$.
$\Box$

\begin{lemma} \label{comb3}
Let $\theta \in \Theta_u$ be equal 
to $[\emptyset,..,\emptyset,\lambda]$, where $\lambda$ is a 
partition of $u$.
Let $t$ be such that $u+t \leq w$.
Then the character summand of
$$Ind_{(\times^p \Sigma_t) \times \Sigma_{up+r}}^{\Sigma_{(t+u)p+r}}
(\chi^{(1^t)} \otimes...\otimes \chi^{(1^t)} \otimes \chi^\theta)$$
{\noindent obtained by removing all character summands indexed by 
partitions outside $\Theta_{u+t}$ is given by}
$$\sum_{\nu} l_\nu \chi^{[\emptyset,...,\emptyset,\nu]},$$
{\noindent where
$Ind_{\Sigma_t \times \Sigma_u}^{\Sigma_{t+u}}
(\chi^{(1^t)} \otimes \chi^\lambda)
= \sum_{\nu} l_\nu \chi^\nu$.}
\end{lemma}
Proof:

We have $\theta = (\rho_1 + p\lambda_1,\rho_2 + p\lambda_2
,...,\rho_u +p\lambda_u, \rho_{u+1},...)$.
By the Littlewood-Richardson rule, the only characters obtained by 
inducing
$\chi^{(1^t)} \otimes \chi^{\theta}$ to $\Sigma_{t+up+r}$ are obtained 
by adding nodes onto $t$ different rows of $\theta$.
Repeating this process $p$ times, the only way of obtaining a
character indexed by a partition in $\Theta_{t+u}$ is by adding nodes
onto the same $t$ rows $p$ times, in such a way that the resulting
partitions lie in $\Theta_{t+u}$. The ways of doing this correspond
exactly to the ways of adding a node to $t$ different
rows of $\lambda$, so
that the resulting composition is a partition.
These correspond exactly to the character summands of 
$Ind_{\Sigma_t \times \Sigma_u}^{\Sigma_{t+u}}
(\chi^{(1^t)} \otimes \chi^\lambda)$, 
by the Littlewood-Richardson rule.
$\Box$

\begin{center}
{\large Permutation modules for $\Sigma_v$.}
\end{center}
 
Suppose that $\lambda = (\lambda_1, \lambda_2,..., \lambda_N)$
is a partition of $w$.
Let $j^1_\lambda = 0$, and let $j^i _\lambda$ \index{$j^i _\lambda$}
be the sum $\sum_{r=1} 
^{i-1} \lambda_r$  
($i=2,3,...,N$).
Let $J^i _\lambda$ \index{$J^i _\lambda$}
($i=1,2,...,N$) be the subgroup,
$$Sym \{ 
pj_\lambda^i+1,pj_\lambda^i+p+1,...,pj_\lambda^i+p(\lambda_i-1)+1
\}$$
$$\times Sym \{ 
pj_\lambda^i+2,pj_\lambda^i+p+2,...,pj_\lambda^i+p(\lambda_i-1)+2 \}
$$ 
$$\times...\times Sym \{ 
pj_\lambda^i+p,pj_\lambda^i+p+p,...,pj_\lambda^i+p\lambda_i \}$$
{\noindent of $\Sigma_v$, a product of $p$ symmetric groups, each of which
is isomorphic to $S_{\lambda_i}$.
Note that $J_\lambda^i$ has support $\{ m \in \mathbb{Z} |
pj_\lambda^i+1 \leq m \leq pj_\lambda^{i+1} \}$.}
Let $$J_\lambda = J^1 _\lambda 
\times J^2 _\lambda \times...\times J^N_\lambda,$$ \index{$J_\lambda$} a
subgroup of $\Sigma_v$ isomorphic to $\times^p \Sigma_\lambda$, 
a direct copy of
$p$ copies of the Young subgroup $\Sigma_\lambda$ of $\Sigma_w$.

\bigskip

Recall that $\Sigma^0_r$ is defined to be 
$Sym \{ wp+1,...,wp+r \}$.
We define $$\Sigma^0_{r+pj_\lambda^i} = 
Sym \{ 1,2,3,...,pj_\lambda^i,wp+1,wp+2,...,wp+r \}.$$
This is the symmetric group whose support is equal to the
support of the direct product
$J_\lambda^1 \times ... \times J_\lambda^{i-1} \times
\Sigma_r^0$.

\bigskip

Let $\xi^i _{\lambda} = \sum_{y \in J^{i} _{\lambda}} y$, 
\index{$\xi_\lambda, \xi^i_\lambda$}
for $i=1,2,...,N$.
Let $$\xi_\lambda = (\xi^1 _\lambda \xi^2 _\lambda...\xi^N_\lambda)  
b_{\rho,0}^\Sigma.b_{\rho,\lambda_1}^\Sigma.
b_{\rho,\lambda_1+\lambda_2}^\Sigma...
b_{\rho,w}^\Sigma,$$  
an element of ${\cal O}\Sigma_v$. 
In this formula, we take  
$b_{\rho,j_\lambda^i}^\Sigma$ to be an element of 
${\cal O}\Sigma_{r+pj_\lambda^i}^0$.
Let $\eta^i _{\lambda} = \sum_{y \in J^{i} _{\lambda}} sgn(y).y$,
\index{$\eta_\lambda$, $\eta^i_\lambda$} 
for $i=1,2,...,N$. 
Let $$\eta_\lambda = (\eta^1 _\lambda \eta^2 _\lambda...\eta^N_\lambda)  
b_{\rho,0}^\Sigma.b_{\rho,\lambda_1}^\Sigma.
b_{\rho,\lambda_1+\lambda_2}^\Sigma...
b_{\rho,w}^\Sigma,$$ 
an element of ${\cal O}\Sigma_v.$

\bigskip

Consider ${\cal O}\Sigma_v \xi_\lambda$. 
This ${\cal O}\Sigma_v$-module may be constructed as follows:

\bigskip

Take the projective module ${\cal O}\Sigma_r$-module 
${\cal O}{\bf B}_{\rho,0}^\Sigma$.
Tensor this module with $p$ copies of the Young module
$Y^{(\lambda_1)}$, 
each of which is
isomorphic to the trivial module for ${\cal O}\Sigma_{\lambda_1}$.
Cut the resulting $p$-permutation
module off at the block of $\Sigma_{r+p\lambda_1}$ with
core $\rho$.
Now tensor this module with $p$ copies of $Y^{(\lambda_2)}$ and cut off
at the block of $\Sigma_{r+p\lambda_1+p\lambda_2}$ with core $\rho$.
Repeat this process until a $p$-permutation
module for ${\cal O}\Sigma_v$ is
obtained, in the block with core $\rho$.

\bigskip

Why can ${\cal O}\Sigma_v \xi_\lambda$ be constucted in this way ? Because
$\xi_\lambda^i$ generates the trivial module for $J_\lambda^i$, which
is isomorphic to $p$ copies of the symmetric group $\Sigma_\lambda^i$.
And because the idempotent $b_{\rho,j_\lambda^i} ^\Sigma$ commutes
with the subgroup $\Sigma^0_{r+pj_\lambda^i}$. 

We have a similar construction of ${\cal O}\Sigma_v \eta_\lambda$.
This time, 
rather than tensoring with $p$ copies of the trivial module
each time before inducing, 
you should tensor with $p$ copies of the alternating module.
\bigskip

When searching for summands of the right $k\Sigma_w \times
\Sigma_r^0$-module $kM$, 
we will be interested in the projective part of 
${\cal O}\Sigma_v \eta_\lambda = ({\cal O}\Sigma_v \xi_\lambda)^\#$.

Pursuing this, let us first note that ${\cal O}\Sigma_v \xi_\lambda$ is a
$p$-permutation module.
This is because a direct summand of a module 
induced from a trivial
module tensored with a sum of Young modules, from a Young subgroup; 
the module ${\cal O}{\bf B}_{\rho,0} ^\Sigma$ for ${\cal O}\Sigma_r^0$ 
is isomorphic to a sum
of Young modules, since that block has defect zero.

Hence it is a direct sum of Young modules. It is known that a Young
module $Y^\mu$ is projective if and only if $\mu$ is
$p$-restricted. Also, a Young module $Y^\mu$ has a Specht filtration,
and the Specht subquotients $S^\gamma$ occuring satisfy 
$\gamma \trianglerighteq \mu$.

So any Young module summand of ${\cal O} \Sigma_v \xi_\lambda$ which is not
projective only has Specht quotients $S^\gamma$ where 
$\gamma \trianglerighteq \mu$ 
for some $p$-restricted partition $\mu$.
 
Now, ${\cal O}\Sigma_v \eta_\lambda = 
({\cal O}\Sigma_v \xi_\lambda)^\# \cong
{\cal O}\Sigma_v \xi_\lambda \otimes sgn$, where $sgn$ is the signature
representation.
And we know that 
$S^\gamma \otimes sgn \cong (S^{\gamma'})^*$, by Remark \ref{dualSpecht}.
Furthermore, tensoring with $sgn$ preserves indecomposability and
takes projectives to projectives and non-projectives to
non-projectives.

Thanks to the above discussion, we know the following fact:

\begin{lemma} \label{constituent}
Let ${\cal O}\Sigma_v \eta_\lambda = U \oplus V$, where $U$ is
projective and $V$ has no projective summand. 
Then the character of $KV$ has irreducible constituents $\chi^\delta$ marked 
by $\delta \trianglelefteq \beta$ 
where $\beta = \beta(\delta)$ is $p$-singular. $\Box$
\end{lemma}

Let $K^i _\lambda$ ($i=1,2,...,N$) \index{$K^i _\lambda$}
be the diagonal subgroup of the direct product $J_\lambda^i$
(here we fix isomorphisms between the factors
of $J^i_\lambda$
which preserve the ordering on $\{ 1,2,..,wp \}$), 
a group isomorphic to $\Sigma_{\lambda_i}$. 
Let $$K_\lambda = K^1 _\lambda \times 
K^2 _\lambda \times...\times K^N_\lambda,$$ \index{$K_\lambda$}
a subgroup of $\Sigma_v$ isomorphic to a Young subgroup 
$\Sigma_\lambda$ of $\Sigma_w$.

Note that $N \cong L \rtimes K_{(w)}$, and $K_\lambda < K_{(w)}$
normalises $D$.
The following lemma is straightforward:

\begin{lemma}

(a) The direct product $J_\lambda ^i \cong \times_p \Sigma_{\lambda_i}$  
is normalised by $P$. The generator $z$ of $P$ acts on $J_\lambda ^i$
by circulating the $p$ direct factors of $J_\lambda ^i$.

\bigskip

(b) The direct product $J_\lambda \cong \times_p \Sigma_{\lambda}$  
is normalised by $P$. The generator $z$ of $P$ acts on $J_\lambda$
by circulating the $p$ direct factors of $J_\lambda$.

\bigskip

(c) The diagonal subgroup $K_\lambda^i$ of $J_\lambda^i$ is equal to
$C_{\Sigma_v}(P) \cap J_\lambda^i$. The diagonal subgroup
$K_\lambda$ of $J_\lambda$ is equal to 
$C_{\Sigma_v}(P) \cap J_\lambda$. 

\bigskip

(d) We have 
$Br_P(\xi_\lambda^1 \xi_\lambda^2 ... \xi_\lambda^N) =
\sum_{K_\lambda} y$, and 

{\noindent $Br_P(\eta_\lambda^1 \eta_\lambda^2 ... \eta_\lambda^N) =
\sum_{K_\lambda} sgn(y).y$
$\Box$}
\end{lemma}

\begin{center}
{\large Summands of $M_{\cal E}$ as Young modules.}
\end{center}

Recall from the proof of theorem~\ref{Wick} that the algebra $A$ acts
on the right of the module $kM$ and its image is the endomorphism
ring $k{\cal E}$. The algebra $A$ commutes with the conjugation action of $P$, and thus
preserves $(kM)^P$ in its action. Hence, 
$k{\cal E} \cong k\Sigma_w \otimes k{\bf B}_{\rho,0}^\Sigma$ preserves $(kM)^P$
in its action.  

\begin{lemma} \label{Karin3}

(a) We have $\eta_{\lambda} \in ({\cal O}\Sigma_v)^P$, and,
$$Br_P(\eta_{\lambda}) = (\sum_{K_\lambda}
sgn(y).y).b_{\rho,0}^\Sigma.$$ 
(b) $Br_P(\eta_\lambda (kM)^P) \cong
Ind_{\Sigma_{\lambda} \times \Sigma^0 _r}^{\Sigma_w \times \Sigma^0 _r}
(sgn \otimes k{\bf B}_{\rho,0}^\Sigma )$, as a right
$\Sigma_w\times \Sigma_r^0$-module. 
\end{lemma}
Proof:

Note that $P$ commutes with $b_{\rho,j_\lambda^i} ^\Sigma$ for every $i$,
and 
$Br_P(b_{\rho,j_\lambda^i}^\Sigma) = b_{\rho,0}^\Sigma$. 
Since the Brauer map is
an algebra homomorphism, we compute the image under $Br_P$ of the
product $\eta_\lambda$ to be
$$Br_P(\eta) = (\sum_{K_\lambda}
sgn(y).y).b_{\rho,0}^\Sigma.$$
This completes the proof of (a).
Thus, 
$$Br_P(\eta_{\lambda}.(kM)^P)$$
$$= Br_P(\eta_{\lambda}(k\Sigma_v)^Pe \zeta_w)$$
$$= (\sum_{y \in K_{\lambda}} sgn(y)y).kC.b_{\rho,0}^\Sigma.\zeta_D$$
$$= \zeta_D.(\sum_{y \in K_{\lambda}} 
sgn(y)y).kC.b_{\rho,0}^\Sigma,$$
{\noindent which is isomorphic to 
$Ind_{\Sigma_{\lambda} \times \Sigma^0 _r}^
{\Sigma_w \times \Sigma^0 _r}(sgn 
\otimes k{\bf B}_{\rho,0}^\Sigma)$ 
as a right $k\Sigma_w \times \Sigma^0 _r$-module.}
$\Box$

\begin{proposition} \label{Q1}
Let $\lambda$ be a partition of $w$. There is an idempotent
$f_\lambda$ \index{$f_\lambda$} in ${\cal O}\Sigma_v$ such that, 

(a) $f_\lambda KM$ is isomorphic to
$Ind_{\Sigma_{\lambda} \times \Sigma_r^0}^{\Sigma_w \times
\Sigma_r^0}(sgn \otimes K{\bf B}_{\rho,0}^\Sigma)$, as a right $K\Sigma_w \times
\Sigma_r^0$-module. 

\bigskip

(b) $f_\lambda kM$ is isomorphic to 
$Ind_{\Sigma_{\lambda} \times \Sigma_r^0}^{\Sigma_w \times \Sigma_r^0}
(sgn \otimes k{\bf B}_{\rho,0}^\Sigma)$, as a right 
$k\Sigma_w \times \Sigma_r^0$-module.

\bigskip

(c) $f_\lambda kM = \eta_\lambda kM = \eta_\lambda (kM)^P \cong
Br_P( \eta_\lambda (kM)^P)$. 
\end{proposition}
Proof:

By lemma~\ref{constituent}, we can write the right ${|cal O}\Sigma_v$-module, 
$\eta_\lambda {\cal O}\Sigma_v = U \oplus V$, where $U$ is
projective and $V$ has no projective summand.
And that the character of $KV$ only has constituents labelled
by $\delta \trianglelefteq \beta$ 
where $\beta = \beta(\delta)$ is $p$-singular.
Since $\eta_\lambda {\cal O}\Sigma_v$ is a quotient of ${\cal O}\Sigma_v$, the
submodule $U$ may be generated by an idempotent. Let $f_\lambda$ 
be such an idempotent, so that $U = f_\lambda {\cal O}\Sigma_v$. We prove (a) and (b) for this idempotent.

Note that the character of the left module $K \Sigma_v f_\lambda$ is the same as the character of the right module $f_\lambda K \Sigma_v$,
since all symmetric group characters are real, and therefore self-dual.
Similarly, the character of the left module $K \Sigma_v \eta_\lambda$ is the same as the character of the right module $\eta_\lambda K \Sigma_v$,

First we show that for $\gamma \in \Theta_w$, we have
$$ (\chi^\gamma,char(K\Sigma_v \eta_\lambda))
= (\chi^\gamma,char(K\Sigma_vf_\lambda)),$$
where $char(X)$ denotes the character of a $K\Sigma_v$-module $X$.

To see this, only observe that for $\gamma \in \Theta_w$, the
formula $(\chi^\gamma,char(KV))=0$ holds, 
for the complement $V$ of $U$ defined above. 
This is a consequence of lemma~\ref{comb2}, for $\Theta_w$ contains
no $p$-singular partitions.

\bigskip

As our second task, we show that $\eta_\lambda KM = f_\lambda KM$,
Indeed, 
$$ \eta_\lambda KM \cong Hom_{\Sigma_v}(K\Sigma_v \eta_\lambda, KM),$$
(note that the algebra $K\Sigma_v$ is semisimple). And by
lemma~\ref{comb1}, $char(KM)$ only has constituents in $\Theta_w$. As
has already been established, the part of $K\Sigma_v \eta_\lambda$ in
$\Theta_w$ is the same as the part of $K\Sigma_v f_\lambda$ in $\Theta_w$,
so 
$$Hom_{\Sigma_v}(K\Sigma_v \eta_\lambda,KM) \cong$$ 
$$Hom_{\Sigma_v}(K\Sigma_v f_\lambda,KM) \cong f_\lambda KM.$$
So $f_\lambda KM \subseteq \eta_\lambda KM$, and the two spaces have the same dimension.
Therefore $f_\lambda KM = \eta_\lambda KM$, as required.

\bigskip

As our third task, we show that $\eta_\lambda kM = f_\lambda kM$. Note that the embedding
$f_\lambda M \hookrightarrow \eta_\lambda M$ splits via left multiplication by $f_\lambda$, and therefore $f_\lambda M$ is ${\cal O}$-pure in $\eta_\lambda M$.
However, $f_\lambda M$ and $\eta_\lambda M$ have the same dimension over $K$, and therefore over ${\cal O}$.
Therefore $f_\lambda M = \eta_\lambda M$, and $\eta_\lambda kM = f_\lambda kM$.

\bigskip

As a fourth task, we should convince ourselves
that the character component of $K\Sigma_v
\eta_\lambda$ which corresponds to $\Theta_w$ is equal to
$$dim(\chi^\rho).\sum_{\nu \vdash w} 
l_\nu \chi^{[\emptyset,..,\emptyset,\nu]} , \leqno{(*)}$$
where $l_\nu$ is defined by 
$Ind_{\Sigma_\lambda}^{\Sigma_w}(sgn \otimes ...\otimes sgn) = 
\sum_\nu l_\nu \chi^\nu$.

From this formula, part (a) of this proposition
follows directly - recall that $- \otimes_{\Sigma_v} KM$ 
matches the character 
$\chi^{[\emptyset,...,\emptyset,\lambda]}$ of $K{\bf B}_{\rho,w}^\Sigma$
with the character $\chi^\lambda \otimes \chi^\rho$ of 
$K\Sigma_w \otimes K{\bf B}_{\rho,0}^\Sigma$.

To see $(*)$, consider our construction of 
${\cal O}\Sigma_v \eta_\lambda$ as a module, tensored, induced and cut,
tensored, induced and cut,....
At the same time meditate upon lemma~\ref{comb3}. The formula is
then visible, inductively.

\bigskip

Our final task is to prove part (b) of the proposition.
The Brauer map gives (by lemma~\ref{Karin2}) a surjection
$$\eta_\lambda (kM)^P \rightarrow 
Ind_{\Sigma_{\lambda} \times \Sigma^0 _r}^{\Sigma_w \times \Sigma^0 _r}
(sgn \otimes k{\bf B}_{\rho,0}^\Sigma).$$
{\noindent
In addition, there is an injection}
$$\eta_\lambda (kM)^P 
\rightarrow \eta_\lambda kM = f_\lambda kM.$$
But part (a) of the proposition shows that
the dimension of $f_\lambda kM$ 
is equal to the dimension
of $Ind_{\Sigma_{\lambda} \times \Sigma^0 _r}^{\Sigma_w \times \Sigma^0 _r}
(sgn \otimes k{\bf B}_{\rho,w}^\Sigma)$, so these two maps must be isomorphisms.
This completes the proof of (b) and (c).
$\Box$

\begin{lemma} \label{Q2}
The right $k\Sigma_w \times \Sigma_r^0$-module $kM$ is
isomorphic to a sum of summands of $\bigoplus_\lambda f_\lambda kM$.
\end{lemma}
Proof:

Let us write ${\cal O}\Sigma_v$ as a direct sum of projective indecomposable modules.
Let $j$ be an idempotent, such that $_{\Sigma_v}{\cal O}\Sigma_v j$ is the sum of projective indecomposables in this decomposition
with simple tops  $\{ D^\lambda | \lambda \in \Theta_w \}$.

Note that ${\cal O}\Sigma_v f_\lambda$ has the projective cover of
$D^{[\emptyset,...,\emptyset,\lambda']}$ as a summand - 
because its character is 
$\chi^{[ \emptyset,...,\emptyset,\lambda']} + 
($ a sum of $\chi^\mu$'s,
$\mu \ntrianglerighteq [\emptyset,...,\emptyset,\lambda]$), by
the formula $(*)$ and lemma~\ref{comb2}.
Thus, $\bigoplus_{\lambda \vdash w} f_\lambda M$ has every summand of
$j M$ as a summand.

But $j M = M$, since
$$i {\cal O}\Sigma_v e \Sigma_w = Hom_{\Sigma_v}({\cal O}\Sigma_v i, M) = 0$$ for
any projective indecomposable module ${\cal O}\Sigma_v i$ with simple top
outside $\{ D^\lambda | \lambda \in \Theta_w \}$ (recall $KM$
has character summands corresponding to elements of 
$\Theta_w$, and $K\Sigma_v i$ has character summands corresponding to
partitions outside $\Theta_w$ by lemma~\ref{comb2}).

The lemma holds ! 
$\Box$

\bigskip

Let $G$ be a finite group, with subgroups $H$ and $K$.
Let $\zeta_H = \sum_{h \in H} h$ \index{$\zeta_H, \zeta_K$}
and let $\zeta_K = \sum_{k \in K}
k$ be corresponding sums in the group algebra $kG$.

\begin{lemma} \label{Q3}

(a)As $kG$-modules, $kG. \zeta_H \cong k[G/H]$ and $kG. \zeta_K
\cong k[G/K]$.

\bigskip

(b)
As vector spaces, $Hom_G(kG \zeta_H , kG \zeta_K) \cong 
k[H \backslash G / K]$ 

\bigskip

(c)
Any element of $Hom_G(kG. \zeta_H , kG. \zeta_K)$ can be written as
right multiplication by some element of $kG$.
\end{lemma}
Proof:

The map $g \zeta_H \mapsto gH$, for $g \in G$ defines a $kG$-module
isomorphism $kG \zeta_H \cong k[G/H]$.
Furthermore,
$$k[H \backslash G /K] \cong Hom_G(k[G/H], k[G/K]),$$ 
$$HgK \mapsto (H \mapsto \sum_{x \in HgK/K} xK).$$
In other words, $Hom_G(kG \zeta_H , kG \zeta_K)$ has a basis $T_s$
indexed by elements $s \in H \backslash G / K$, where 
$$T_{HgK}: \zeta_H \mapsto 
\sum_{x \in HgK} x = \zeta_H .g. 
\sum_{y \in {\cal T}} y,$$
{\noindent and ${\cal T}$ is a set of representatives for 
$(g^{-1}Hg \cap K) \backslash K$.}
So $T_{HgK}$ can be defined as right multiplication by 
$\sum_{y \in {\cal }} y$. $\Box$

\begin{lemma} \label{Q4}
The natural map $f_\lambda k\Sigma_v f_\mu \rightarrow 
Hom_{k\cal E}( f_\mu kM , f_\lambda kM)$ is surjective.
\end{lemma}
Proof:

Since $f_\lambda kM = \eta_\lambda kM$, left multiplication by $x \in
k\Sigma_v$ on $f_\lambda kM \subseteq k\Sigma_v$ is equivalent to left
multiplication by $xf_\lambda \in k\Sigma_v$.
Similarly, if $x.y \in \eta_\mu kM = f_\mu kM$, then $f_\mu x.y = x.y$.
To prove the lemma then, it is sufficient to show that every element
of $Hom_{\cal E}(\eta_\lambda kM, \eta_\mu kM)$ is given (as left
multiplication) by an element of $k\Sigma_v$.

\bigskip

Recall from lemma~\ref{Q1} and the proof of lemma~\ref{Karin3} that, 
$$\eta_\lambda kM = \eta_\lambda (kM)^P \cong Br_P(\eta_\lambda (kM)^P)$$ 
$$= \zeta_D.(\sum_{y \in K_{\lambda}} 
sgn(y)y).kC.b_{\rho,0}^\Sigma$$
$$\cong 
(\sum_{y \in K_\lambda} sgn(y).y) kK_{(w)} \otimes k{\bf B}_{\rho,0} ^\Sigma.$$  
Thanks to lemma~\ref{Q3}, any homomorphism between,
$$(\sum_{y \in K_\mu} sgn(y).y) kK_{(w)}, \textrm{ and}$$ 
$$(\sum_{y \in K_\lambda} sgn(y).y) kK_{(w)},$$
is given as left multiplication by an element of $kK_{(w)}$.
See this by applying the
automorphism $\#$.
Thus, any homomorphism between,
$$\eta_\mu kM \cong (\sum_{y \in K_\mu} sgn(y).y) kK_{(w)} \otimes k{\bf B}_{\rho,0}^\Sigma, 
\textrm{ and}$$
$$\eta_\lambda kM \cong
(\sum_{y \in K_\lambda} sgn(y).y) kK_{(w)} \otimes k{\bf B}_{\rho,0}^\Sigma,$$ 
is given as left multiplication by  an element $z$ of $kK_{(w)} \otimes k{\bf B}_{\rho,0}^\Sigma$.

To complete the proof of the lemma, we need to show we can choose a preimage of $z \in kK_{(w)}$ under $Br_P$ which sends $\eta_\mu (kM)^P$ to $\eta_\lambda (kM)^P$.
There is a diagonal embedding $\Delta: \Sigma_w \hookrightarrow \times^p \Sigma_w \cong J_{(w)}$, whose image in $J_{(w)}$ is $K_{(w)}$. 
The group $P$ acts on $J_{(w)}$ by rotating the $p$ copies of $\Sigma_w$.
Thinking of $z$ as an element of $k\Sigma_w$, we can consider the element $t = z^{\otimes p} \in kJ_{(w)}$.
Note that $t$ commutes with $P$. Writing $z = \sum_{g \in \Sigma_w} a_g g$ as a linear combination of group elements $g$, we compute
$$Br_P(t) = \sum (a_g)^p \Delta(g).$$
Since $\mathbb{F}_p$ is a splitting field, we may take $k = \mathbb{F}_p$, so that $(a_g)^p = a_g$, for all $a_g$, and thus
$$Br_P(t) = \sum a_g \Delta(g) = z \in kK_{(w)}.$$
Multiplying on the left by $i.t$, where $i$ is the idempotent factor of $\eta_\mu$ defines a map $\eta_\mu (kM)^P \rightarrow \eta_\lambda (kM)^P$.
Computing the effect of this map via the Brauer morphism, we find it corresponds to multiplication by $z$ inside $K_{(w)}$, as required.
$\Box$

\begin{corollary}
The natural algebra homomorphism from $k{\bf B}_{\rho,w}^\Sigma$ 
to the endomorphism ring
$End_{k{\cal E}}(kM)$ is a surjection.
\end{corollary}
Proof:

The algebra $k{\cal E}$ is defined to be the endomorphism ring
$End_{\Sigma_v}(kM)$.
So $k\Sigma_v$ maps to $End_{k{\cal E}}(kM)$. 

By proposition \ref{Q1} (b), and lemma \ref{Q2}, and Schur-Weyl duality, $End_{k{\cal E}}(kM)$ is Morita equivalent to the Schur algebra,
since the indecomposable summands of tensor space as a symmetric group module are precisely the twisted Young modules.  
The collection of isomorphism classes of irreducible $End_{k{\cal E}}(kM)$ -modules is therefore in bijection with the collection of partitions of $w$.

The algebra $k{\bf B}_{\rho,w}^\Sigma/kI$ is Morita equivalent to a quotient of the Schur algebra ${\cal S}(v,v)$, with isomorphism classes of irreducible modules
in bijection with $\Theta_w$, which is in bijection with the collection of partitions of $w$.
The algebras $End_{k{\cal E}}(kM)$ and $k{\bf B}_{\rho,w}^\Sigma/kI$ therefore have the same number of isomorphism classes of irreducible modules.

Thanks to lemma~\ref{Q2}
and lemma~\ref{Q4}, for any primitive idempotents 
$i,j \in k{\bf B}_{\rho,w}^\Sigma/kI$, the 
natural map 
$$i (k{\bf B}_{\rho,w}^\Sigma/kI) j \rightarrow 
Hom_{\cal E}( j kM , i kM)$$ 
{\noindent is surjective.}
Since $k{\cal E}$ and $k{\bf B}_{\rho,w}^\Sigma/kI$ have the same number of isomorphism classes of irreducible modules, the map
from $k\Sigma_v$ to $End_{k{\cal E}}(kM)$ is surjective.
$\Box$
  
\begin{remark} In the light of theorem~\ref{Q1} and lemma~\ref{Q2}, 
theorem~\ref{Comeon} is proven. $\Box$
\end{remark}

\begin{center}
{\large Remarks and questions.}
\end{center}

Recall that the set of irreducible
characters for the Schur algebra $K{\cal S}(w,w)$ may be parametrized
$\{ \chi(\lambda) | \lambda$ \index{$\chi(\lambda)$}
a partition of $w \}$, 
and the set of simple $k{\cal S}(w,w)$-
modules $\{ L(\lambda) | \lambda$ a partition of $w \}$.
In this way, $\chi(\lambda)$ corresponds under (non-twisted)
Schur-Weyl duality to
the character $\chi^\lambda$ of $\Sigma_w$.
In addition,
$\chi(\lambda)$ has a single composition factor
isomorphic to $L(\lambda)$, and for any other composition factor 
$L(\mu)$, we have $\lambda \vartriangleright \mu$. 

\bigskip

Let ${\cal A}= End_{{\cal E}} (M)$. \index{${\cal A}$} 
Then ${\cal A}$ is ${\cal O}$-free, and $K{\cal A} \cong K{\bf B}_{\rho,w}^\Sigma/KI$ 
is Morita equivalent to the Schur algebra
$K{\cal S}(w,w)$, by theorem~\ref{Comeon}.
By convention, we match by this Morita equivalence
the character $\chi(\lambda)$ of
$K{\cal S}(w,w)$ with the character
$\chi^{[\emptyset,...,\emptyset,\lambda']}$ of $K\Sigma_v$.

Theorem~\ref{Comeon} also informs us that 
$k{\cal A} = k{\bf B}_{\rho,w}^\Sigma/kI$
is Morita equivalent to the Schur algebra $k{\cal S}(w,w)$.
How do simple modules correspond ?

\begin{proposition} \label{Hahaha}
Under the Morita equivalence between $k{\cal S}(w,w)$ 
and $k{\bf B}_{\rho,w}^\Sigma/kI$ of
theorem~\ref{Comeon}, the simple module $L(\lambda)$ corresponds to 
$D^{[\emptyset,...,\emptyset,\lambda']}$.
\end{proposition}
Proof:

For $\lambda, \mu$ partitions of $w$, let the multiplicities 
$m^\lambda_\mu$ \index{$m^\lambda_\mu$}
be defined by,
$$Ind_{\Sigma_\lambda} ^{\Sigma_w} K = 
\sum_\mu m^\lambda _\mu .\chi^\mu.$$
We know for that $m^\lambda _\mu$ is zero for all partitions
$\lambda \vartriangleright
\mu$ of $w$, and that $m^\lambda _\lambda = 1$ for all partitions
$\lambda$ of $w$.

\bigskip

The right $K\Sigma_w \times \Sigma_r$-module 
$f_\lambda KM$ has character,
$$dim(\chi^\rho)
\sum_\mu m^\lambda_\mu (\chi^{\mu'} \otimes \chi^\rho).$$
Consider $Hom_{\cal E}(M, f_\lambda M)$, 
a left ${\cal A}$-module, which has character
$$dim(\chi^\rho)
\sum_\mu m^\lambda _\mu \chi^{[\emptyset,...,\emptyset,\mu']},$$ 
Reduced modulo $p$, this is the
$k{\cal A}$-module, 
$$Hom_{\Sigma_w \times \Sigma_r}(kM,f_\lambda kM),$$ 
{\noindent
which via Morita equivalence, corresponds to
to dim$(\chi^\rho)$ copies of the 
$k{\cal S}(w,w)$-module,   }
$$E^{\otimes r \#}.(\sum_{\sigma \in \Sigma_\lambda} sgn(\sigma) \sigma).$$
{\noindent
In other words, we have $dim(\chi^\rho)$ copies of the 
$k{\cal S}(w,w)$-module,     }
$$E^{\otimes r}.(\sum_{\sigma \in \Sigma_\lambda} \sigma).$$
This module is
isomorphic to a $p$-modular reduction of the
$K{\cal S}(w,w)$-module with character
$dim(\chi^\rho)
\sum_\mu m^\lambda _\mu.\chi(\mu)$.

Induction according to the dominance ordering now implies that the
$p$-modular reduction of an $K{\cal S}(w,w)$-module with character
$\chi(\lambda)$, mapped under Morita equivalence to a $k{\cal
A}$-module, has the same composition factors as a $p$-modular
reduction of a $K\Sigma_v$-module with character
$\chi^{[\emptyset,...,\emptyset,\lambda']}$.
The lower unitriangularity of the decomposition matrices of 
$k{\cal S}(w,w)$ and
$k\Sigma_v$ implies the result.
$\Box$

\begin{corollary} 
The decomposition matrix of $k{\cal S}(w,w)$ 
is a submatrix of the decomposition
matrix of $k{\bf B}_{\rho,w} ^\Sigma$, where row and column $\lambda$ corresponds
to row and column $[\emptyset,...,\emptyset,\lambda']$.
$\Box$
\end{corollary}

\begin{remark} It is not possible that
$k\Sigma_v/Ann(k\Sigma_ve) \cong End_{kNf}(k\Sigma_ve)$ when $p=2$.

Because on one hand (by \cite{Land}, 1.14.5)
the $k$-dimension of $End_{ek\Sigma_ve}(k\Sigma_ve)$ 
is at least as great as the
$K$-dimension of $End_{KNf}(K\Sigma_ve) \cong b_\rho(K\Sigma_v)$ (this
isomorphism still holds when $w \geq p$, by the character
calculation in \cite{KessCh}).
On the other hand, 
the dimension of $k\Sigma_v/Ann(k\Sigma_ve)$ is strictly smaller than the
dimension of $K{\bf B}_{\rho,w}^\Sigma$, because
there are elements of $k{\bf B}_{\rho,w}^\Sigma$ which act on $k\Sigma_ve$ as zero.  
For example, let $xk\Sigma_v$ be a simple right 
$k{\bf B}_{\rho,w} ^\Sigma$-module, 
not in the top of $ek\Sigma_v$ (such exist, since $ek\Sigma_ve$ and
$k{\bf B}_{\rho,w} ^\Sigma$ are not Morita equivalent).
In this situation, $xk\Sigma_ve \cong Hom(ek\Sigma_v,xk\Sigma_v) = 0$.
\end{remark}  

\begin{remark}

(a) The modules in the top (and socle) of $_{\Sigma_v} kM$ are those $D^
    {[\emptyset,...,\emptyset,\lambda]}$'s such that $\lambda$ is
    $p$-regular, by theorem~\ref{Comeon} and proposition~\ref{Hahaha}.

\bigskip

(b) The Young module summands of $_{\Sigma_v} kM$ are those $Y^
    {[\emptyset,...,\emptyset,\lambda]}$'s such that $\lambda$ is
    $p$-regular, by theorem~\ref{Comeon} and proposition~\ref{Hahaha}.
\end{remark}

\begin{question}
The following are open:

When $p=2$, what are the summands of $k\Sigma_ve_N$ ? 
What are their vertices and what are their sources ?

When $p=2$, what is the vertex of 
$_{\Sigma_v} k\Sigma_ve _N$ ? What is the source ?
\end{question}

\newpage
\begin{center}
{\bf \large 
Chapter VI

Ringel duality inside Rock blocks of symmetric groups.}
\end{center}

We prove that a Rock block of a symmetric group, of arbitrary defect,
possesses a family of internal symmetries, given as Ringel dualities between
various subquotients (theorem~\ref{Ringel}). 
Since this sequence of symmetries resembles 
J.A. Green's walk around the Brauer tree \cite{Greentree}, we name it
``a walk along the abacus''.

\begin{center}
{\large A criterion for Ringel duality.}
\end{center}

Let $(K, {\cal O}, k)$ be an $l$-modular system.
We prove a sufficient condition for Ringel
duality between two split
quasi-hereditary ${\cal O}$-algebras.
 
\bigskip

Let $R$ be a commutative Noetherian ring.
Cline, Parshall, and Scott have defined split quasi-hereditary $R$-algebras (\cite{CPS}, 3.2).
For example, the Schur algebra ${\cal S}(n,r)$, defined over $R$,
is a split quasi-hereditary algebra, with respect to the
poset $\Lambda(n,r)$.
(see \cite{CPS}, 3.7).

More generally, for $\Gamma$ an ideal of $\Lambda(n,r)$, and
$\Omega$ a coideal of $\Lambda(n,r)$,
the generalised Schur algebra ${\cal S}(\Gamma \cap \Omega)$ is
a split quasi-hereditary subquotient of ${\cal S}(n,r)$, with respect
to the poset $\Gamma \cap \Omega$.

\begin{definition}
For a split quasi-hereditary ${\cal O}$-algebra $A$, 
let us define a $K$-$k$- \emph{tilting module} $T$ to be a
finitely generated $A$-order which is a tilting module
over $K$, as well as a tilting module over $k$.
\end{definition}

The following result resembles M. Brou\'e's theorem~\ref{balls}, 
which gives sufficient conditions for a
Morita equivalence between symmetric ${\cal O}$-algebras:
 
\begin{theorem} \label{Brother}
Let $A, B$ be split quasi-hereditary algebras over
${\cal O}$ with respect to posets $\Lambda, \Upsilon$. 
Suppose that $KA, KB$ are semisimple.
Let $T$ be an $A$-$B$-bimodule which is a $K$-$k$-tilting
module at the same time as a left $A$-module and as a right
$B$-module.
Suppose that the functors 
$$KT \otimes_{KB} - : KB-mod \rightarrow KA-mod$$
$$- \otimes_{KA} KT : mod-KA \rightarrow mod-KB$$ 
are equivalences of categories, such that the resulting bijections
of irreducible modules define order-reversing maps
between $\Lambda$ and $\Upsilon$.
Then $kT$ defines a Ringel duality between $kA$ and $kB^{op}$.
Furthermore, $B \cong End_A(T)$, and $A \cong End_B(T)$.
\end{theorem}
Proof:

Let $\{ P(\mu), \mu \in \Upsilon \}$ be the set of
non-isomorphic principal indecomposable 
$B$-modules.
Let $\Psi$ be the order-reversing
map from $\Upsilon$ to $\Lambda$ defined by $KT$.

The $A$-module $T$ has $|\Upsilon|$ distinct
summands $T \otimes_{B} P(\mu)$.
Viewed over $K$, such a summand
has a single composition factor $K\Delta(\Psi \mu)$, and all other
composition factors $K\Delta(\lambda), \lambda < \Psi \mu$.
This is because $\Psi$ is order-reversing.
Viewed over $k$, such a summand has a 
single composition factor $L(\Psi \mu)$,
and all other composition factors $L(\lambda), \lambda < \Psi \mu$.   
In other words, $kT \otimes_{kB} P(\mu)$ is the indecomposable
tilting module $T(\Psi \mu)$ for $kA$. 
Thus, $kT$ is a full tilting module as a left $kA$-module.

Likewise, $kT$ is a full tilting module as a right $kB$-module.
It follows that $kA$ and $kB$ act faithfully on $kT$.

Note that there are natural isomorphisms
$KB \cong End_{KA}(KT)$ and $KA \cong End_{KB}(KT)$ of algebras.
 
The endomorphism ring
$End_A(T)$ is an ${\cal O}$-order, so that
$kEnd_A(T)$ injects into $End_{kA}(kT)$, and such
that
$KEnd_A(T)$ surjects onto 
$End_{KA}(KT)$, which is isomorphic to $KB$.
The dimension of $End_{kA}(kT)$ is given (\cite{Donkin}, A.2.2(ii))
by the formula,
$$\sum_\nu [kT : \Delta_k(\nu)][kT : \nabla_k(\nu)].$$
The dimension of $End_{KA}(KT)$ is given by the formula,
$$\sum_\nu [KT : \Delta_K(\nu)][KT : \nabla_K(\nu)].$$

It follows from the equality of these formulae, 
that the rank of $End_A(T)$ is equal to the rank of $B$,
and $kEnd_A(T) \cong End_{kA}(kT)$.
Indeed, $kEnd_A(T) \cong End_{kA}(kT) \cong kB$,
since $kB$ acts faithfully on $kT$.

We have now
verified that $kA$ and $kB$ are in Ringel duality.
Since the natural map from $kB$ to $kEnd_A(T)$ is an isomorphism, 
$B$ is ${\cal O}$-pure in $End_A(T)$, and so
$B \cong End_A(T)$.
Likewise, $A \cong End_B(T)$.
This completes the proof of theorem~\ref{Brother}.
$\Box$

\begin{center}
{\large Combinatorial preliminaries.}
\end{center}

Let $p$ be a prime number and $w$ any natural number.
Let $\rho = \rho(p,w)$ be a minimal Rouquier core. Thus,
in an abacus presentation, $\rho$ has precisely $w-1$ more beads on the
$i^{th}$ runner than on the $i-1^{th}$ runner.

\begin{lemma} \label{cock}
Let $\tilde{\rho}$ be an arbitrary Rouquier core.
Then ${\cal O}{\bf B}_{\tilde{\rho},w}^\Sigma$ is
Morita equivalent to ${\cal O}{\bf B}_{\rho,w}^\Sigma$.
The resulting correspondence of partitions preserves $p$-quotients. 
\end{lemma}
Proof:

Scopes' isometries \cite{Scopes}
correspond to the motion of $k$ beads one runner
leftwards on the abacus, where $k \geq w$.
Suppose that $\tilde{\rho}$ has $N_i \geq (w-1)$ more beads on runner
$i$ than on runner $i-1$. Suppose that $N_j > (w-1)$.
Then there is a Scopes isometry 
which moves those $N_j$ beads from
runner $j$ to runner $j-1$, which may be followed by a Scopes isometry
which moves $N_j+N_{j+1}$ beads from runner $j+1$ to runner $j$,....., 
which may be followed by a Scopes isometry which moves 
$N_j+N_{j+1}+...+N_{p-1}$ beads from runner $p-1$ ro runner $p-2$.
If $j>1$, then the result is a $p$-core with $N'_{j-1} > (w-1)$ 
more beads on runner $j-1$ than on runner $j-2$.

We may thus proceed with a further $p-j$ Scopes isometries, moving
beads from runners $j-1,...,p-2$ leftwards.
Following this procedure to its natural conclusion (always pushing
beads leftwards), and at last circulating the ordering of runners on
the abacus so that
runner $0$ becomes runner $j$, we obtain a $p$-core
$\tilde{\rho}_0$, smaller than $\tilde{\rho}$ which still satisfies the
condition that there are $\geq (w-1)$ more beads on runner $i$ than
on runner $i-1$.
By induction, lemma~\ref{cock} is proven. 
$\Box$

\bigskip

Let 
$${\cal I} =
\{ \lambda \quad | \quad \lambda \textrm{ has core } 
\rho  \textrm{ and weight }w \}$$ \index{${\cal I}$}
be the indexing poset of $k{\bf B}_{\rho,w} ^{\cal S}$.
We wish to describe the dominance order on ${\cal I}$. 
To $\lambda \in {\cal I}$ with $p$-quotient 
$[\lambda^0,\lambda^1,...,\lambda^{p-1}]$, let us associate an element
$<\lambda> \in \mathbb{N}^{wp}$, \index{$\lambda$} given by
$$(\lambda^{p-1}_1,\lambda^{p-1}_2,...,\lambda^{p-1}_w,
\lambda^{p-2}_1,\lambda^{p-2}_2,...,\lambda^{p-2}_w,...,
\lambda^0_1,\lambda^0_2,...,\lambda^0_w).$$
Let us place the dominance order on ${\mathbb N}^{wp}$.
We then have:

\begin{proposition} \label{ohdear}
Let $\lambda, \mu \in {\cal I}$.
In the dominance order,
$\mu \unlhd \lambda$ if and only if $<\mu> \unlhd <\lambda>$.
\end{proposition}
Proof:

Following lemma~\ref{cock}, and 
applying a series of Scopes isometries (which preserve the
dominance order), we may replace $\rho$ by $\tilde{\rho}$, where
$\tilde{\rho}$ has at least $N>2w$ beads on runner $i$ than runner
$i-1$ (for $i=1,...,p-1$).
Let 
$$\tilde{\cal I} =
\{ \lambda \quad | \quad 
\lambda \textrm{ has core }\tilde{\rho} \textrm{, and weight }
w \}.$$
{\noindent Suppose that $\lambda \lhd \mu$ are neighbours in
$\tilde{\cal I}.$
There is a sequence
$\lambda=\lambda_0 \lhd \lambda_1 \lhd... \lhd \lambda_m = \mu$, where
the Young diagram of $\lambda_{j-1}$ is obtained from the Young diagram of
$\lambda_j$ by removing a box and placing it lower in the diagram
(for $j=1,...,m$).
We may assume that no box is removed from the core in this sequence of 
motions.}

Thus, at each step
we must remove a box from one of the $p$-hooks which has been
added to $\tilde{\rho}$ to create $\mu$.
Since $\tilde{\rho}$ has $N$ more beads on runner $i$ than runner
$i-1$ (for $i=1,...,p-1$), we must actually remove an entire $p$-hook
if we are to end our sequence in $\tilde{\cal I}$.
Correspondingly, we must add an entire $p$-hook when we add boxes.
The new additional $p$-hook must appear lower 
in the Young diagram than the old 
removed $p$-hook.

On an abacus, the corresponding motion looks as follows:
move a single bead one place higher on its runner,
move a second bead (necessarily above the first bead) one
place lower on \emph{its} runner.

This corresponds precisely to removing one box from the Young diagram
of $< \lambda >$, and
replacing it lower in the Young diagram to obtain $< \mu >$.
Thus, $< \lambda > \lhd < \mu >$.

Working backwards, we find conversely that 
$< \lambda > \lhd < \mu >$ implies $\mu \lhd \lambda$.
$\Box$

\bigskip

We may deduce from the above proposition 
a number of combinatorial results concerning ideals and
coideals of the poset ${\cal I}$, ordered by the
dominance ordering.

\begin{definition}
For natural numbers $a_1,...,a_{p-1}$ 
such that $\sum_{i=0}^{p-1} a_i = w$, let 
$${\cal I}_{(a_0,a_1,..,a_{p-1})}= 
\left\{  
\begin{array}{cc} 
\lambda \in {\cal I} \quad | & \quad \lambda \unlhd \mu   
\textrm{ for some } \mu \in {\cal I} \textrm{ with }p 
\textrm{-quotient}\\  
& [\mu^0,...,\mu^{p-1}] \textrm{, such that } 
|\mu^i| = a_i  
\end{array} 
\right\}$$ \index{${\cal I}_{(a_0,a_1,..,a_{p-1})}$}
$${\cal J}_{(a_0,a_1,..,a_{p-1})}=
\left\{  
\begin{array}{cc} 
\lambda \in {\cal I} \quad | & \quad \mu \unlhd \lambda   
\textrm{ for some } \mu \in {\cal I} \textrm{ with }p 
\textrm{-quotient}\\ 
& [\mu^0,...,\mu^{p-1}] \textrm{, such that } 
|\mu^i| = a_i  
\end{array} 
\right\}$$ \index{${\cal J}_{(a_0,a_1,..,a_{p-1})}$}
For $i=0,...,p-1$, let 
$${\cal I}_i = {\cal I}_{(0,...,0,w,0,...,0)}$$ \index{${\cal I}_i$}
$${\cal J}_i = {\cal I}_{(0,...,0,w,0,...,0)},$$ \index{${\cal J}_i$}
where the $w$ appears as the $i^{th}$ entry in $(0,...,0,w,0,...,0)$.
Let
$${\cal I}_{res} = {\cal I}_{p-2}, \hspace{0.5cm}
{\cal I}_{unres} = {\cal I} - {\cal I}_{res} \hspace{0.5cm}
{\cal I}_{reg} = {\cal J}_1, \hspace{0.5cm}
{\cal I}_{sing} = {\cal I} - {\cal I}_{reg}.$$
\index{${\cal I}_{res}, {\cal I}_{unres}, {\cal I}_{reg}, {\cal I}_{sing}$}
\end{definition}

Following these definitions,
proposition~\ref{ohdear} has a number of corollaries, which are
easily checked:

\begin{corollary} \label{mark1}
These subsets of ${\cal I}$ are ideals:
$${\cal I}_{(a_0,a_1,..,a_{p-1})},
{\cal I}_i,  {\cal I}_{res},  {\cal I}_{sing}.$$
These subsets of ${\cal I}$ are coideals:
$${\cal J}_{(a_0,a_1,..,a_{p-1})},
{\cal J}_i, {\cal I}_{unres}, {\cal I}_{reg}. \Box$$
\end{corollary}

\begin{corollary}
For $i=0,...,p-1$,
$${\cal I}_i =
\left\{  
\begin{array}{cc}  
\lambda \in {\cal I} &  \textrm{ with }p \textrm{-quotient }  
[\lambda^0,...,\lambda^{p-1}],  
\textrm{ such that } \\ 
& |\lambda^{i+1}|=|\lambda^{i+2}|=...=|\lambda^{p-1}|=0   
\end{array}  
\right\}$$  
$${\cal J}_i = 
\left\{  
\begin{array}{cc} 
\lambda \in {\cal I} &  \textrm{ with }p \textrm{-quotient } 
[\lambda^0,...,\lambda^{p-1}], 
\textrm{ such that } \\ 
& |\lambda^0|=|\lambda^1|=...=|\lambda^{i-1}|=0  
\end{array} 
\right\}. \Box$$
\end{corollary}

\begin{corollary} \label{mark2}
$${\cal I}_{res} = \{ \lambda \in {\cal I} \quad | \quad  
\lambda \textrm{ is }
p \textrm{-restricted } \}.$$ 
$${\cal I}_{unres} = \{ \lambda \in {\cal I} \quad | \quad 
\lambda \textrm{ is } p \textrm{-nonrestricted }
\}.$$
$${\cal I}_{reg} = \{ \lambda \in {\cal I} \quad | \quad
\lambda \textrm{ is }p \textrm{-regular } \}.$$
$${\cal I}_{sing} = \{ \lambda \in {\cal I} \quad | \quad 
\lambda \textrm{ is } p \textrm{-singular  } \}. \Box$$
\end{corollary}

\begin{corollary} \label{spse}
The intersection
${\cal I}_{(a_0,a_1,..,a_{p-1})} \cap
{\cal J}_{(a_0,a_1,..,a_{p-1})}$ is equal to the set,
$$\{ \lambda \in {\cal I} \quad | \quad \lambda \textrm{ has }p \textrm{-quotient }
[\lambda^0,...,\lambda^{p-1}], \textrm{ where } |\lambda^i| = a_i \}. \Box$$
\end{corollary}

The intersection of ideals in the above lemma is analogous to
the classical intersection 
$\{ \mu \unlhd \lambda \} \cap \{ \mu \unrhd \lambda \} = \{ \lambda \}$.
Whilst this classical intersection may be used to index the characters
of symmetric groups by partitions, we use the intersection of ideals
above to study Rock blocks runner by runner.

\begin{definition}
For natural numbers $a_0,...,a_{p-1}$ such that 
$\sum_{i=0}^{p-1} a_i=w$, let   
$${\cal K}_{(a_0,...,a_{p-1})} = 
{\cal I}_{(a_0,...,a_{p-1})} \cap 
{\cal J}_{(a_0,...,a_{p-1})} =$$
\index{${\cal K}_{(a_0,...,a_{p-1})}$}
$$\{ \lambda \in {\cal I} \quad | \quad \lambda \textrm{ has }p \textrm{-quotient }
[\lambda^0,...,\lambda^{p-1}], \textrm{ where }
|\lambda^i| = a_i \}.$$
For $i=0,...,p-1$, let  
$${\cal K}_i = {\cal I}_i \cap {\cal J}_i =$$
\index{${\cal K}_i$}
$$\left\{   
\begin{array}{cc} 
\lambda \in {\cal I} \quad | & \lambda \textrm{ has }p \textrm{-quotient } 
[\lambda^0,...,\lambda^{p-1}], 
\textrm{ such that }  \\
& |\lambda^0|=...=|\lambda^{i-1}|= 
|\lambda^{i+1}|=...=|\lambda^{p-1}|=0, |\lambda^i|=w  
\end{array} 
\right\}.$$
\end{definition}

\begin{center}
{\large Quasi-hereditary subquotients of 
${\cal O}{\bf B}_{\rho,w} ^\Sigma$.} 
\end{center}

We introduce pairs of quasi-hereditary quotients of 
${\cal O}{\bf B}_{\rho,w} ^\Sigma$, 
which are isomorphic under the signature automorphism.

\bigskip
Let $\omega = (1^v)$, a partition of $v$.
Let ${\cal S}(v,v)$ be the Schur algebra associated to polynomial
representations of $GL_v$ of degree $v$,
defined over ${\cal O}$.
Recall that $\xi_\omega {\cal S}(v,v) \xi_\omega$ is isomorphic to 
${\cal O}\Sigma_v$ (theorem~\ref{tensorspace2}).

\bigskip
 
To an ideal $\Gamma$ of ${\cal I}$, let us associate the ideal
$X_{\Gamma}$ of ${\cal S}(v,v)$, the quotient by which, is 
the generalised Schur algebra, ${\cal S}(\Gamma)$. 
Let $I_\Gamma = \xi_\omega b_{\rho,w}^{\cal S}X_{\Gamma} \xi_\omega$ 
be the corresponding ideal of ${\cal O}{\bf B}_{\rho,w}^\Sigma$.

\begin{lemma}
Suppose that $\Gamma$ is an ideal of ${\cal I}_{res}$.
Then ${\cal O}{\bf B}_{\rho,w}^\Sigma / I_\Gamma$ is a quasi-hereditary algebra, with
indexing poset 
$\Gamma$.
The decomposition matrix of the algebra 
${\cal O}{\bf B}_{\rho,w}^\Sigma/I_\Gamma$ is the
square  
submatrix of the decomposition matrix of ${\cal O}{\bf B}_{\rho,w}^\Sigma$, whose
rows are indexed by elements of $\Gamma$.
\end{lemma}
Proof:

The generalised Schur algebra
${\cal O}{\bf B}_{\rho,w}^{\cal S}/ X_{\Gamma}$ is a quasi-hereditary algebra,
whose indexing poset is $\Gamma$.

Over the field $k$, the idempotent
$\xi_\omega$ sends to zero precisely those simple modules indexed
by unrestricted partitions.
Thus, $({\cal O}{\bf B}_{\rho,w}^{\cal S}/X_{\Gamma}) \xi_{\omega}$ 
is a progenerator
for ${\cal O}{\bf B}_{\rho,w}^{\cal S}/X_{\Gamma}$.
It follows that, 
$$\xi_\omega ({\cal O}{\bf B}_{\rho,w}^{\cal S}/X_{\Gamma}) \xi_\omega 
\cong {\cal O}{\bf B}_{\rho,w}^\Sigma/ \xi_\omega X_{\Gamma} \xi_\omega
= {\cal O}{\bf B}_{\rho,w}^\Sigma / I_\Gamma,$$ 
is Morita
equivalent to ${\cal O}{\bf B}_{\rho,w}^{\cal S}/X_{\Gamma}$. 
$\Box$

\bigskip
Let $a_1,...,a_{p-1}$ be natural numbers such that  
$\sum_{i=1}^{p-1} a_i = w$.

\bigskip

We write  $I_{(a_1,...,a_{p-1})}$ \index{$I_{(a_1,...,a_{p-1})}$} for 
$I_\Gamma$, where $\Gamma$ is the ideal
${\cal I}_{(a_1,...,a_{p-1},0)}$ in ${\cal I}$.

Thus, 
${\cal O}{\bf B}_{\rho,w}^\Sigma / I_{(a_1,...,a_{p-1})}$ is
a quasi-hereditary algebra whose poset is 
the ideal ${\cal I}_{(a_1,...,a_{p-1},0)}$, and 
$\{ D_\lambda | \lambda \in {\cal I}_{(a_1,...,a_{p-1},0)} \}$
is a complete set of non-isomorphic simple 
$k{\bf B}_{\rho,w}^\Sigma/I_{(a_1,...,a_{p-1})}$-modules.
 
\bigskip

So long as $i=1,...,p-1$, let us write $I_i$ \index{$I_i$}
for $I_{(0,...,0,w,0,...,0)}$, where 
$w$ appears as the $i-1^{th}$ entry in $(0,...,0,w,0,...,0)$.

Thus, 
${\cal O}{\bf B}_{\rho,w}^\Sigma/ I_{i}$ is
a quasi-hereditary algebra whose poset is 
${\cal I}_{i-1}$.

\bigskip

Let us write $I_{unres}$ \index{$I_{unres}$}for $I_{p-2}$.
Thus, ${\cal O}{\bf B}_{\rho,w}^\Sigma / I_{unres}$ is a quasi-hereditary algebra
whose poset is ${\cal I}_{res}$.

\bigskip

Let $\Omega \subset {\cal I}$ be a coideal. 
Let $x_\Omega$ \index{$x_\Omega$}
be an idempotent in ${\cal O}\Sigma_v$, such
that ${\cal O}\Sigma_v x_\Omega$ is a maximal summand of 
${\cal O}\Sigma_v$, whose indecomposable summands have tops in
$\{ D_\lambda | \lambda \in 
{\cal I}_{res} \cap \Omega \}$. 

\begin{lemma}
Suppose that $\Omega$ is a coideal of ${\cal I}$, and that
$\Gamma$ is an ideal of ${\cal I}_{res}$.
Let $A_{\Gamma \cap \Omega} = 
x_\Omega ({\cal O}{\bf B}_{\rho,w} ^\Sigma / I_\Gamma) x_\Omega$.
Then $A_{\Gamma \cap \Omega}$ is a quasi-hereditary algebra, with
indexing poset
$\Gamma \cap \Omega$.
\end{lemma}
Proof:

$A_{\Gamma \cap \Omega}$ is Morita equivalent to the generalised Schur
algebra with indexing poset 
$\Gamma \cap \Omega$, via the bimodule
$\xi_\Omega ({\cal S}(v,v) / J_{\Gamma \cap \Omega}) \xi_\omega.x_\Omega$.
$\Box$

\bigskip
Suppose that $a_1,...,a_{p-1}$ are natural numbers such that  
$\sum_{i=1}^{p-1} a_i = w$.

\bigskip

We write $A_{(a_1,...,a_{p-1})}$ \index{$A_{(a_1,...,a_{p-1})}$} for 
$A_{\Gamma \cap \Omega}$, where 
$\Omega = {\cal J}_{(a_1,...,a_{p-1},0)}$, and
$\Gamma = {\cal I}_{(a_1,...,a_{p-1},0)}$.
Indeed,
$A_{(a_1,...,a_{p-1})}$ is
a quasi-hereditary algebra whose poset is 
${\cal K}_{(a_1,...,a_{p-1},0)}$, by corollary~\ref{spse}.

\bigskip

So long as $i=1,...,p-1$, we write 
$A_i$ \index{$A_i$} for $A_{(0,...,0,w,0,...,0)}$, where 
$w$ appears as the $i-1^{th}$ entry in $(0,...,0,w,0,...,0)$.
Thus, 
$A_i$ is
a quasi-hereditary algebra whose poset is 
${\cal K}_{i-1}$.

\bigskip

For natural numbers $a_1,...,a_{p-1}$ such that  
$\sum_{i=1}^{p-1} a_i = w$,
let $J_{(a_1,...,a_{p-1})}$ \index{$J_{(a_1,...,a_{p-1})}$} 
be the ideal $I_{(a_{p-1},...,a_1)} ^\#$.
Thus, 
${\cal O}{\bf B}_{\rho,w}^\Sigma / J_{(a_1,...,a_{p-1})}$ is
a quasi-hereditary algebra whose poset is 
${\cal J}_{(0,a_1,...,a_{p-1})} ^{op}$.

We write 
$\{ D^\lambda | \lambda \in {\cal J}_{(0,a_1,...,a_{p-1})} \}$
for the set of simple 
$k{\bf B}_{\rho,w}^\Sigma/J_{(a_1,...,a_{p-1})}$-modules.
 
\bigskip

For $i=1,...,p-1$, let $J_i = I_{p-i} ^\#$. \index{$J_i$}
Thus, 
${\cal O}{\bf B}_{\rho,w}^\Sigma / J_{i}$ is
a quasi-hereditary algebra whose poset is 
${\cal J}_i ^{op}$.

Let $J_{sing} = I_{unres} ^\#$. \index{$J_{sing}$}
Thus, 
${\cal O}{\bf B}_{\rho,w}^\Sigma / J_{sing}$ is
a quasi-hereditary algebra whose poset is 
${\cal I}_{reg} ^{op}$.

\bigskip

Let $y_{(a_1,...,a_{p-1})} = x_{(a_{p-1},...,a_1)} ^\#$. 
\index{$y_{(a_1,...,a_{p-1})}$}
Let 
$$B_{(a_1,...,a_{p-1})} =
y_{(a_1,...,a_{p-1})} \left(
{\cal O}{\bf B}_{\rho,w}^\Sigma / J_{(a_1,...,a_{p-1})} \right) 
y_{(a_1,...,a_{p-1})}.$$ \index{$B_{(a_1,...,a_{p-1})}$}
{\noindent Thus, 
$B_{(a_1,...,a_{p-1})}$ is
a quasi-hereditary algebra whose poset is 
${\cal K}_{(0,a_1,...,a_{p-1})} ^{op}$.}

\bigskip

So long as $i=1,...,p-1$, let us write $B_i$ \index{$B_i$} for 
$B_{(0,...,0,w,0,...,0)}$, where 
$w$ appears as the $i^{th}$ entry in $(0,...,0,w,0,...,0)$.
Thus, 
$B_i$ is
a quasi-hereditary algebra whose poset is 
${\cal K}_{i} ^{op}$.

\begin{remark}
Simple $kA_{(a_1,...,a_{p-1})}$-modules
are in natural one-one
correspondence with the set
$\{ D_\lambda $ $|$
$\lambda \in {\cal K}_{(a_1,...,a_{p-1},0)} \}$.
By theorem~\ref{Rowena}, this set is equal to the set
$\{ D^\lambda$ $|$
$\lambda \in {\cal K}_{(0,a_1,...,a_{p-1})} \}$.

Simple $kB_{(a_1,...,a_{p-1})}$-modules are also
in natural one-one
correspondence with the set
$\{ D^\lambda$ $|$
$\lambda \in {\cal K}_{(0,a_1,...,a_{p-1})} \}$.
\end{remark}

By balancing the algebras ${\cal O}{\bf B}_{\rho,w}^\Sigma/I_{unres}$ and 
${\cal O}{\bf B}_{\rho,w} ^\Sigma/J_{sing}$ on the Mullineux map,
we reveal Ringel dualities between different runners of
${\cal O}{\bf B}_{\rho,w}^\Sigma$, in the following section.

\begin{center}
{\large Walking along the abacus.}
\end{center}

Let $p \geq 3$.
For $\sum_{i=1} ^{p-2} a_i = w$,
consider the ${\cal O}$-lattice 
$$N_{(a_1,...,a_{p-2})} = x_{(0,a_1,...,a_{p-2})} 
{\cal O}\Sigma_v y_{(a_1,...,a_{p-2},0)}.$$ \index{$N_{(a_1,...,a_{p-2})}$}
In this section we prove that 
$N_{(a_1,...,a_{p-2})}$ provides a Ringel duality
between the quasi-hereditary subquotients 
$A_{(0,a_1,...,a_{p-2})}$ and $B_{(a_1,...,a_{p-2},0)}$
of ${\cal O}{\bf B}_{\rho,w}^\Sigma$.

These Ringel dualities should be viewed as internal symmetries of the
Rock block. For simple 
$kA_{(0,a_1,...,a_{p-2})}$-modules are in natural correspondence with
simple $k{\bf B}_{\rho,w}^\Sigma$-modules $D^\lambda$ indexed by elements of 
${\cal K}_{(0,0,a_1,...,a_{p-2})}$.
At the same time, 
simple $kB_{(a_1,...,a_{p-2},0)}$-modules are in natural
correspondence with simple $k{\bf B}_{\rho,w} ^\Sigma$-modules $D^\mu$ indexed by
elements of
${\cal K}_{(0,a_1,...,a_{p-2},0)}$. 

These symmetries therefore enable us to translate module-theoretic information 
along the abacus.

\bigskip

Here's a technical lemma: 
 
\begin{lemma} \label{bum}
Suppose that $p \geq 3$.
Let $\sum_{i=1}^{p-2} a_i = w$. 
Then 
$$J_{sing}. x_{(0,a_1,...,a_{p-2})} = 0,$$
$$J_{(a_1,...,a_{p-2},0)}. 
x_{(0,a_1,...,a_{p-2})} = 
x_{(0,a_1,...,a_{p-2})}. 
J_{(a_1,...,a_{p-2},0)} = 0,$$
$$I_{unres}. y_{(a_1,...,a_{p-2},0)} = 0,$$
$$I_{(0,a_1,...,a_{p-2})}. y_{(a_1,...,a_{p-2},0)} =
y_{(a_1,...,a_{p-2},0)}.I_{(0,a_1,...,a_{p-2})} = 0.$$
\end{lemma}
Proof:

The character of ${\cal O}\Sigma_v. x_{(0,a_1,...,a_{p-2})}$ 
has irreducible
components $\chi^\lambda$, where $\lambda$ lies in the
coideal ${\cal J}_{(0,a_1,...,a_{p-2},0)}$ of ${\cal I}$.
This character has no
components which are $p$-singular, and no components which lie in
${\cal I}-{\cal J}_{(0,a_1,...,a_{p-2},0)}$. 

Over $K$, the ideal
$J_{sing}$ is equal to the Wedderburn component of
$K{\bf B}_{\rho,w}^\Sigma$ with $p$-singular components, whilst 
$J_{(a_1,...,a_{p-2},0)}$ is equal
to the
Wedderburn component with components in 
${\cal I}-{\cal J}_{(0,a_1,...,a_{p-2},0)}$. 
This, along with the fact that irreducible characters for symmetric groups
are self-dual, proves the first two parts of the 
lemma.
The third and fourth parts follow analogously.
$\Box$

\bigskip

We now show that 
$N_{(a_1,...,a_{p-2})}$ is an 
$A_{(0,a_1,...,a_{p-2})}$-$B_{(a_1,...,a_{p-2},0)}$-bimodule: 

\begin{lemma}
The kernel of the natural map from 
$x_{(0,a_1,...,a_{p-2})} {\cal O}\Sigma_v x_{(0,a_1,...,a_{p-2})}$ 
to $End(N_{(a_1,...,a_{p-2})})$ contains the ideal 
$x_{(0,a_1,...,a_{p-2})} I_{(0,a_1,...,a_{p-2})} 
x_{(0,a_1,...,a_{p-2})}$.

The kernel of the natural map from 
$y_{(a_1,...,a_{p-2},0)} {\cal O}\Sigma_v y_{(a_1,...,a_{p-2},0)}$ 
to the ring
$End(N_{(a_1,...,a_{p-2})})$, contains the ideal 
$y_{(a_1,...,a_{p-2},0)} J_{(a_1,...,a_{p-2},0)} 
y_{(a_1,...,a_{p-2},0)}$.

Thus, $N_{(a_1,...,a_{p-2})}$ is an 
$A_{(0,a_1,...,a_{p-2})}$-$B_{(a_1,...,a_{p-2},0)}$-bimodule.
\end{lemma}
Proof:
$$x_{(0,a_1,...,a_{p-2})} I_{(0,a_1,...,a_{p-2})} 
x_{(0,a_1,...,a_{p-2})} {\cal O}\Sigma_v 
y_{(a_1,...,a_{p-2},0)}$$ 
$$\subseteq 
x_{(0,a_1,...,a_{p-2})} I_{(0,a_1,...,a_{p-2})} 
y_{(a_1,...,a_{p-2},0)}=0,$$
by lemma~\ref{bum}.
Likewise, 
$$x_{(0,a_1,...,a_{p-2})} {\cal O}\Sigma_v 
y_{(a_1,...,a_{p-2},0)} J_{(a_1,...,a_{p-2},0)} 
y_{(a_1,...,a_{p-2},0)}$$ 
$$\subseteq 
x_{(0,a_1,...,a_{p-2})} J_{(a_1,...,a_{p-2},0)} 
y_{(a_1,...,a_{p-2},0)} =
0. \Box$$

\begin{theorem} \label{Ringel}
(``Walking along the abacus'')
The bimodule $N_{(a_1,...,a_{p-2})}$ defines a Ringel
duality between $kA_{(0,a_1,...,a_{p-2})}$ and 
$k{\bf B}_{(a_1,...,a_{p-2},0)} ^{op}$.
\end{theorem}
Proof:

We first show that
$N_{(a_1,...,a_{p-2})}$ is a $K$-$k$-tilting module both 
as a left $A_{(0,a_1,...,a_{p-2})}$-module, and as a right
$B_{(a_1,...,a_{p-2},0)}$-module.

\bigskip

Let $R \in \{ K, k \}$.
Recall that under Schur-Weyl duality, Specht modules correspond to 
costandard modules.
Therefore,
the costandard modules for $R{\bf B}_{\rho,w}^\Sigma/I_{unres}$ are those
Specht modules indexed by restricted partitions.
The costandard modules for 
$R{\bf B}_{\rho,w}^\Sigma/I_{(0,a_1,...,a_{p-2})}$ are those Specht modules
indexed by elements of ${\cal I}_{(0,a_1,...,a_{p-2},0)}$.
The costandard modules for
$A_{(0,a_1,...,a_{p-2})}$ are those modules 
$x_{(0,a_1,...,a_{p-2})}S$, where $S$ is a Specht module indexed by
an element of ${\cal K}_{(0,a_1,...,a_{p-2},0)}$.

Since $R\Sigma_v y_{(a_1,...,a_{p-2})}$ is a projective module, it has
a filtration by Specht modules.
Suppose that $S$ is a Specht module in this filtration.
Then $S$ is a costandard
module for $R{\bf B}_{\rho,w}^\Sigma/I_{(0,a_1,...,a_{p-2})}$,
by lemma~\ref{bum}.
Thus, $x_{(0,a_1,...,a_{p-2})} S$ is a costandard module for
$A_{(0,a_1,...,a_{p-2})}$.
So $N_{(a_1,...,a_{p-2})}$ has a filtration by costandard modules.

Since $R\Sigma_v y_{(a_1,...,a_{p-2})}$ is self-dual, it may
also be filtered by dual Specht modules.
The same argument as above now shows that
$N_{(a_1,...,a_{p-2})}$ has a filtration by standard modules.

Thus, $N_{(a_1,...,a_{p-2})}$ is a left $K$-$k$-tilting module for
$A_{(0,a_1,...,a_{p-2})}$.
In the same way, $N_{(a_1,...,a_{p-2})}$ is a right $K$-$k$-tilting module
for $B_{(a_1,...,a_{p-2},0)}$.

\bigskip

We would like to apply theorem~\ref{Brother}.

Recall that the $K\Sigma_v$-$K\Sigma_v$- bimodule $K\Sigma_v$ has character,
$\bigoplus_\lambda \chi^\lambda \otimes \chi^\lambda$.
In this way, the $x_{(0,a_1,...,a_{p-2})} K\Sigma_v x_{(0,a_1,...,a_{p-2})}$-
$y_{(a_1,...,a_{p-2},0)} K\Sigma_v y_{(a_1,...,a_{p-2},0)}$-module
$KN_{(a_1,...,a_{p-2})}$ has character,
$$\bigoplus_{\lambda \in {\cal K}_{(0,a_1,...,a_{p-2},0)}} 
x_{(0,a_1,...,a_{p-2})}.\chi^\lambda \otimes
\chi^\lambda. y_{(a_1,...,a_{p-2},0)}.$$
Note that,
$$\{ x_{(0,a_1,...,a_{p-2})}.\chi^\lambda \quad | \quad \lambda \in 
{\cal K}_{(0,a_1,...,a_{p-2},0)} \},$$ 
is a complete set of irreducible left
$KA_{(0,a_1,...,a_{p-2})}$-modules.
And that, 
$$\{ \chi^\lambda.y_{(a_1,...,a_{p-2},0)} \quad | \quad \lambda \in
{\cal K}_{(0,a_1,...,a_{p-2},0)} ^{op} \},$$ 
is a complete set of irreducible
right $KB_{(a_1,...,a_{p-2},0)}$-modules.

Thus, $KN_{(a_1,...,a_{p-2})}$ induces a Morita equivalence between
$KA_{(0,a_1,...,a_{p-2})}$ and $KB_{(a_1,...,a_{p-2},0)}$, which reverses 
order on the indexing posets.
It is now a consequence of theorem~\ref{Brother}, that
$kA_{(0,a_1,...,a_{p-2})}$, and 
$kB_{(a_1,...,a_{p-2},0)} ^{op}$ are in Ringel duality.
$\Box$

\begin{remark} \label{corr}
Under the Morita equivalence provided by
$KN_{(a_1,...,a_{p-2})}$, a simple
$KA_{(0,a_1,...,a_{p-2})}$-module $S^\lambda _A$ 
corresponds to a simple $KB_{(a_1,...,a_{p-2},0)}$-module
$S^\lambda _B$. 
\end{remark}

\begin{note}
Both $A_{(a_1,...,a_{p-1})}$ and  
$B_{(a_1,...,a_{p-1})}$ have simple modules  
in natural correspondence 
with $\{ D^\lambda | \lambda \in
{\cal K}_{(0,a_1,...,a_{p-1})} \}$.
Let ${\cal L}(a_1,...,a_{p-1})$ \index{${\cal L}(a_1,...,a_{p-1})$}
be the Serre subcategory
of $k\Sigma_v-mod$ generated by
$\{ D^\lambda | \lambda \in
{\cal K}_{(0,a_1,...,a_{p-1})} \}$. 
\end{note}

\newpage
\begin{center}
{\bf \large 
Chapter VII

James adjustment algebras for Rock blocks of symmetric groups.}
\end{center}

Let $k{\bf B}_{\rho,w}^\Sigma$ 
be a Rock block of a symmetric group, whose weight is $w$.  
We show that 
there is a nilpotent
ideal ${\cal N}$ of $k{\bf B}_{\rho,w}^\Sigma$, such that
$k{\bf B}_{\rho,w}^\Sigma /{\cal N}$ is Morita equivalent to a direct sum,
$$\bigoplus_{\substack{a_1,...,a_{p-1} \in \mathbb{Z}_{\geq 0} \\ \sum a_i =w}} 
\left( \bigotimes _{i=1} ^{p-1} {\cal S}(a_i,a_i) \right) ,$$
of tensor products of Schur algebras (theorem~\ref{nilpotent}).
The decomposition matrix of this quotient 
is equal to the James adjustment matrix of $k{\bf B}_{\rho,w}^\Sigma$  \cite{James10}.

In chapter five, we proved the existence of a quotient of 
$k{\bf B}_{\rho,w} ^\Sigma$,
equivalent to ${\cal S}(w,w)$.
In this chapter, we show that the Ringel dualities of chapter VI,
may be applied simultaneously with the ideas of chapter V,
to set up an induction, 
proving the existence of a quotient $k{\bf B}_{\rho,w}^\Sigma/{\cal N}$, 
described above.

Although we choose not explicitly to describe its proof here, 
there exists a generalisation of the main result of this chapter 
to arbitrary Hecke algebras of type $A$.
Indeed, there exists a nilpotent
ideal ${\cal N}$ of $k{\bf B}_{\rho,w}^{{\cal H}_q}$, such that
$k{\bf B}_{\rho,w}^{{\cal H}_q} /{\cal N}$ is Morita equivalent to a direct sum,
$$\bigoplus_{\substack{a_1,...,a_{p-1} \in \mathbb{Z}_{\geq 0} \\ \sum a_i =w}} 
\left( \bigotimes _{i=1} ^{p-1} {\cal S}(a_i,a_i) \right) ,$$
of tensor products of \emph{unquantized} Schur algebras.
The decomposition matrix of this quotient 
is equal to the James adjustment matrix of $k{\bf B}_{\rho,w}^{{\cal H}_q}$.

\begin{center}
{\large The James adjustment algebra of a Hecke algebra.} 
\end{center}

Let $(K, {\cal O},k)$ be an $l$-modular system.
Let $\wp$ be the maximal ideal of ${\cal O}$, so that ${\cal O}/\wp \cong k$.
Let $q \in {\cal O}$ be a primitive $p^{th}$ root of unity, whose
image in $k$, is non-zero.

This section is devoted to the Hecke algebra ${\cal H}_q(\Sigma_n)$,
We define a quotient of this algebra, whose representation theory
controls the James adjustment matrix of ${\cal H}_q(\Sigma_n)$
(see \cite{Geck}).

\bigskip

Let us label the simple $K{\cal H}_q(\Sigma_n)$-modules, 
$$\{ D_q ^\lambda \quad | \quad \lambda \textrm{ a }p \textrm{-regular partition of }
n \},$$ \index{$D_q^\lambda$}
in contrast with the simple $k{\cal H}_q(\Sigma_n)$-modules, which we label,
$$\{ D^\lambda \quad | \quad \lambda \textrm{ a }p \textrm{-regular partition of }
n \}.$$   
The algebra $K{\cal H}_q(\Sigma_n)$ is far from 
semisimple in general, and
thus has a non-trivial radical.
Let ${\cal O}{\cal N}_q$ \index{${\cal O}{\cal N}_q$}
be the intersection of this radical with the
subalgebra ${\cal O}{\cal H}_q(\Sigma_n)$ of $K{\cal H}_q(\Sigma_n)$.
Thus, ${\cal O N}_q$ is equal to the annihilator 
in ${\cal OH}_q(\Sigma_n)$ of all simple
$K{\cal H}_q(\Sigma_n)$-modules.
The ideal
${\cal ON}_q$ is an ${\cal O}$-pure sublattice of 
${\cal OH}_q(\Sigma_n)$.

\bigskip

The algebra 
${\cal G}_q(\Sigma_n) =
{\cal OH}_q(\Sigma_n)/ {\cal ON}_q$ \index{${\cal G}_q(\Sigma_n)$}
is an 
${\cal O}$-free algebra, whose square decomposition matrix, 
$([D_q ^\lambda : D^\mu])$ is equal to the
James adjustment matrix of ${\cal H}_q(\Sigma_n)$.
We call it the \emph{James adjustment algebra} of ${\cal H}_q(\Sigma_n)$.
Indeed, we have (cf. \cite{Geck}, 2.3),
$$Dec_{\wp}({\cal OH}_q(\Sigma_n)) = 
Dec_{<t-q>}(K[t]_{(t-q)} {\cal H}_t(\Sigma_n)) \times
Dec_{\wp}({\cal G}_q(\Sigma_n)).$$
Here, we write 
$Dec_{\cal J}(A)$ for the decomposition
matrix of an algebra $A$, defined over a ring $R$, 
relative to a maximal ideal ${\cal J}$.

\bigskip

Cutting ${\cal G}_q(\Sigma_n)$ at a block $b_{\tau,w}^{{\cal H}_q}$ of ${\cal H}_q(\Sigma_n)$, 
we obtain the James adjustment algebra of the block, which we denote ${\cal G}b_{\tau,w} ^{{\cal H}_q}$.
If $q=1$, we label this block ${\cal G}b_{\tau,w} ^\Sigma$.

\begin{center}
{\large Preliminaries on adjustment algebras for Rock blocks.} 
\end{center}

We now concentrate on Rock blocks of symmetric groups.
Thus, we assume that the image of $q$ in $k$, is equal to $1$, 
and $k$ has characteristic $l=p$.
And we adopt the notation of chapters 4-6.

\bigskip

Let $x_{(a_1,...,a_{p-1}), q}$ \index{$x_{(a_1,...,a_{p-1}), q}$}
be the $q$-analogue of
$x_{(a_1,...,a_{p-1})}$, an idempotent of $K{\bf B}_{\rho,w} ^{{\cal H}_q}$.
Let $A_{(a_1,...,a_{p-1}), q}$ \index{$A_{(a_1,...,a_{p-1}), q}$}
be the $q$-analogue
of $A_{(a_1,...,a_{p-1})}$, and let 
$B_{(a_1,...,a_{p-1}), q}$ \index{$B_{(a_1,...,a_{p-1}), q}$} 
be the $q$-analogue of $B_{(a_1,...,a_{p-1})}$,.
These are subquotients of $K{\bf B}_{\rho,w} ^{{\cal H}_q}$.

For natural numbers $a_i$ such that $\sum_{i=1} ^{p-1} a_i = w$, 
let ${\cal L}_q(a_1,...,a_{p-1})$ \index{${\cal L}_q(a_1,...,a_{p-1})$}
be the Serre subcategory of
$K{\bf B}_{\rho,w}^{{\cal H}_q}-mod$ generated by simple modules
$$\{ D_q^\lambda \quad | \quad \lambda \in {\cal K}_{(0,a_1,...,a_{p-1})} \}.$$ 

Let ${\cal N}_{a_1,...,a_{p-1},q}$ \index{${\cal N}_{a_1,...,a_{p-1},q}$}
be the ideal of 
$K{\bf B}_{\rho,w} ^{{\cal H}_q}$, 
the quotient by which, is equal to
the Wedderburn component of the quotient $K{\cal
G}_q(\Sigma_v)$ with components in ${\cal L}_q(a_1,...,a_{p-1})$.

\begin{proposition} \label{anima}
There are isomorphisms of algebras, 
$$K{\bf B}_{\rho,w} ^{{\cal H}_q}/K{\cal N}_{a_1,...,a_{p-1}, q}
\cong KA_{(a_1,...,a_{p-1}), q} \cong 
KB_{(a_1,...,a_{p-1}), q}.$$
\end{proposition}

Proof:
$$K{\bf B}_{\rho,w} ^{{\cal H}_q}/K{\cal N}_{a_1,...,a_{p-1}, q}$$
$$= x_{(a_1,...,a_{p-1}), q}
\left( K{\bf B}_{\rho,w} ^{{\cal H}_q}/K{\cal N}_{a_1,...,a_{p-1},q} \right)
x_{(a_1,...,a_{p-1}), q},$$ 
which is equal to the
Wedderburn component of,
$$x_{(a_1,...,a_{p-1}), q}
\left( K{\bf B}_{\rho,w} ^{{\cal H}_q}/Radical \right)
x_{(a_1,...,a_{p-1}), q},$$
whose simple components correspond 
to simple objects of 
${\cal L}_q(a_1,...,a_{p-1})$.}

But the simple $KA_{(a_1,...,a_{p-1}), q}$-modules are in one-one
correspondence  
with simple objects of 
${\cal L}_q(a_1,...,a_{p-1})$.
Thus,
$KA_{(a_1,...,a_{p-1}), q}$ surjects onto the quotient
$$K{\bf B}_{\rho,w} ^{{\cal H}_q}/K{\cal N}_{a_1,...,a_{p-1},q}.$$
However,
$KA_{(a_1,...,a_{p-1}), q}$ is semisimple by 
proposition~\ref{backdoor}(1), so this
surjection is an isomorphism.

Applying $\#$ to this isomorphism, 
we discover that in addition, 
$$K{\bf B}_{\rho,w} ^{{\cal H}_q}/K{\cal N}_{a_{p-1},...,a_1, q}
\cong KB_{(a_{p-1},...,a_1), q}.
\Box$$

\begin{center}
{\large Quotients of $k{\bf B}_{\rho,w} ^\Sigma$.}
\end{center}

The algebras $A_{(a_1,..,a_{p-1})}$ and 
$B_{(a_1,...,a_{p-1})}$ may be realised as 
quotients of $k{\bf B}_{\rho,w} ^\Sigma$, and not merely as subquotients.
This will be shown in general in proposition~\ref{Polly}, 
and is crucial to the proof of our main result, theorem~\ref{nilpotent}.
As a preliminary, in this section we prove that
$A_{(0,a_2,...,a_{p-2},0)}$ and
$B_{(0,a_2,...,a_{p-2},0)}$ may be realised as quotients of 
$k{\bf B}_{\rho,w} ^\Sigma$.

We also prove here a baby version 
of theorem~\ref{nilpotent},  
so you may see how these ideas fit together
with those of chapter 6 to provide information on
$k{\bf B}_{\rho,w} ^\Sigma-mod$.

\begin{lemma} \label{oooh}
Suppose that $p \geq 3$.
Let $\sum_{i=1} ^{p-2} a_i = w$. 

Then $$k{\bf B}_{\rho,w} ^\Sigma/I_{(a_1,...,a_{p-2},0)}-mod \cong$$ 
$$\{ M \in k\Sigma_v-mod | M \textrm{ has composition factors }
D^\lambda, \lambda \in {\cal I}_{(0,a_1,...,a_{p-2},0)} 
\cap {\cal I}_{reg} \}.$$
Let $i=2,...,p-1$. Then 
$$k{\bf B}_{\rho,w} ^\Sigma/J_{(0,a_1,...,a_{p-2})}-mod \cong$$
$$\{ M \in k\Sigma_v-mod \quad | \quad M \textrm{ has composition factors }
D^\lambda, \lambda \in {\cal J}_{(0,0,a_1,...,a_{p-2})} \}.$$
\end{lemma}
Proof:

Let $p \geq 3$, and $\sum_{i=1} ^{p-2} a_i = w$. Then,
$$\{ M \in k\Sigma_v-mod | M \textrm{ has simple factors }
D^\lambda, 
\lambda \in {\cal I}_{(0,a_1,...,a_{p-2},0)} 
\cap {\cal I}_{reg}\}$$
$$\subseteq \{
M \in k\Sigma_v-mod \quad | \quad M \textrm{ is generated by }
k\Sigma_v y_{(a_1,...,a_{p-2},0)} \} \cong$$
$$\{ M \in k{\bf B}_{\rho,w} ^\Sigma/I_{unres}-mod \quad | \quad
M \textrm{ is generated by } 
k\Sigma_v y_{(a_1,...,a_{p-2},0)} \},$$
{\noindent where the latter isomorphism is by lemma~\ref{bum}.
Thus,}
$$\{ M \in k\Sigma_v-mod | M \textrm{ has simple factors }D^\lambda, 
\lambda \in {\cal I}_{(0,a_1,...,a_{p-2},0)} 
\cap {\cal I}_{reg} \}$$
$$= 
\left\{  
\begin{array}{cc} 
M \in k{\bf B}_{\rho,w} ^\Sigma /I_{unres}-mod \quad | &  
M \textrm{ has simple factors }D^\lambda, \\
& \lambda \in {\cal I}_{(0,a_1,...,a_{p-2},0)} 
\cap {\cal I}_{reg}
\end{array}
\right\}$$
$$= 
\left\{
\begin{array}{cc}
M \in k{\bf B}_{\rho,w} ^\Sigma /I_{unres}-mod \quad | & 
M \textrm{ has simple factors }  
D_\lambda, \\ 
& \lambda \in {\cal I}_{(a_1,...,a_{p-2},0,0)}
\end{array}
\right\}$$
$$\cong 
\left\{
\begin{array}{cc}
M \in k{\bf B}_{\rho,w} ^\Sigma/I_{(a_1,...,a_{p-2},0)}-mod \quad | & 
M \textrm{ has simple factors }
D_\lambda \\ 
& \lambda \in {\cal I}_{(a_1,...,a_{p-2},0,0)}
\end{array}
\right\},$$
{\noindent
where the latter isomorphism holds thanks to the quasi-heredity of 
the quotient
$k{\bf B}_{\rho,w} ^\Sigma/I_{unres}$.
We deduce that, for $i=1,...,p-2$,
$$k{\bf B}_{\rho,w} ^\Sigma/I_{(a_1,...,a_{p-2},0)}-mod \cong$$
$$\{ M \in k\Sigma_v-mod | M \textrm{ has composition factors }
D^\lambda, \lambda \in {\cal I}_{(0,a_1,...,a_{p-2},0)} 
\cap {\cal I}_{reg} \}.$$
The second part of the lemma follows on an application of $\#$.
$\Box$

\bigskip

Note the following obvious fact:

\begin{lemma} \label{triv}
Suppose that $I,J$ are ideals in an algebra $A$.
Then 
$$A/(I+J)-mod \cong (A/I-mod) \cap (A/J-mod). \Box$$
\end{lemma}

Note that $A_i$ and $B_i$ both have simple modules 
in natural correspondence
with the set $\{ D^\lambda$ $|$ $\lambda \in {\cal K}_i \}$.
Let 
${\cal L}(i) = {\cal L}(0,...,0,w,0,...,0)$ be 
a ``single runner subcategory'' of $k\Sigma_v-mod$,
associated to runner $i$, for $i=1,...,p-1$.
Thus, ${\cal C}(i)$ is the Serre subcategory of $k \Sigma_v-mod$ generated by
$$\{D^\lambda | \lambda \in {\cal K}_i \}.$$
We may now produce ``single runner quotients'' of 
$k{\bf B}_{\rho,w} ^\Sigma$, for 
$p \geq 3$.

\begin{proposition} \label{subcat}
Let $p \geq 3$.

For $i=1,...,p-1$, 
there exists an ideal $\alpha_i$ of $k{\bf B}_{\rho,w} ^\Sigma$, such that 
$$k{\bf B}_{\rho,w} ^\Sigma / \alpha_i-mod \cong {\cal L}(i).$$
There are natural algebra isomorphisms, 
$$kA_i \cong kB_i \cong k{\bf B}_{\rho,w} ^\Sigma / \alpha_i.$$
\end{proposition}
Proof:

For ideals $\alpha_i$, we take $I_i + J_i$, so long as 
$2 \leq i \leq p-2$. 
We take $\alpha_1 = I_1$, and $\alpha_{p-1} = J_{p-1}$.
From lemmas~\ref{oooh} and \ref{triv} above it is clear that
$k{\bf B}_{\rho,w} ^\Sigma / \alpha_i-mod \cong {\cal L}(i)$.

Let $i=1,...,p-1$.
We have $k{\bf B}_{\rho,w} ^\Sigma/\alpha_i 
= y_i \left( k{\bf B}_{\rho,w} ^\Sigma/\alpha_i \right) y_i$, since any 
projective summand of $k{\bf B}_{\rho,w} ^\Sigma$, maximal 
subject to the restriction that
its top lies in ${\cal C}(i)$, must be
isomorphic to a summand of $k\Sigma_v y_i$. 

\bigskip

Note that $\alpha_i = I_{i, q} + J_{i, q} ($mod $p)$, where
$I_{i,q}$ (resp. $J_{i,q}$, $\alpha_{i,q}$) is the $q$-analogue of 
$I_i$ (resp. $J_i$, $\alpha_i$), an ideal of ${\cal O}{\bf B}_{\rho,w} ^{{\cal H}_q}$.

By proposition~\ref{anima}, and lemma~\ref{triv}, we know that
${\cal O}{\bf B}_{\rho,w}^{{\cal H}_q}/\alpha_{i, q} = KB_{i,q}$.
We therefore witness the inclusion,
$y_{i,q} I_{i, q} y_{i,q}
\subseteq y_{i,q} J_{i,q} y_{i,q}$,
over the field $K$.
Taking intersections with the natural ${\cal O}$-form for 
$K{\cal H}_q(\Sigma_v)$, we reveal the inclusion
$y_{i, q} I_{i, q} y_{i, q}
\subseteq y_{i, q} J_{i, q} y_{i, q}$, over ${\cal O}$.
 
Reducing modulo $p$, we find that
$y_i J_i y_i$ contains $y_i I_i y_i$.
In conclusion,
$$k{\bf B}_{\rho,w} ^\Sigma/\alpha_i = y_i \left( k{\bf B}_{\rho,w} ^\Sigma/\alpha_i \right) y_i
\cong y_i \left( k{\bf B}_{\rho,w} ^\Sigma/J_i \right) y_i = kB_i,$$ 
{\noindent for $i=1,...,p-1$.
Likewise, $k{\bf B}_{\rho,w} ^\Sigma/\alpha_i \cong kA_i$, for $i=1,...,p-1$.
$\Box$}

\bigskip

Generalising the above proposition and its proof, we have, 
 
\begin{proposition} \label{Varagn}
Let $p >3$, and let $\sum_{i=2}^{p-2} a_i =w$.

There exists an ideal ${\cal N}_{0,a_2,...,a_{p-2},0}$ 
of $k{\bf B}_{\rho,w} ^\Sigma$ such that, 
$$k{\bf B}_{\rho,w} ^\Sigma / {\cal N}_{0,a_2,...,a_{p-2},0}-
mod \cong {\cal L}(0,a_2,...,a_{p-2},0).$$
There are algebra isomorphisms, 
$$kA_{(0,a_2,...,a_{p-2},0)} 
\cong kB_{(0,a_2,...,a_{p-2},0)}
\cong k{\bf B}_{\rho,w} ^\Sigma / {\cal N}_{0,a_2,...,a_{p-2},0}. \Box$$
\end{proposition}

To complete this section, we bear an infant
theorem~\ref{nilpotent}.

\begin{proposition} \label{Schur}

(a) Let $p=2$. There is a nilpotent ideal $\alpha_1$ of
    $k{\bf B}_{\rho,w} ^\Sigma$, such
    that $$k{\bf B}_{\rho,w} ^\Sigma/\alpha_1-mod \cong {\cal S}(w,w)-mod.$$

(b) Let $p \geq 3$, and let $i=1,...,p-1$. There are equivalences of
    abelian categories, 
    $${\cal L}(i) \cong k{\bf B}_{\rho,w} ^\Sigma/\alpha_{i}-mod \cong
    {\cal S}(w,w)-mod.$$  
\end{proposition}
Proof:

(a) This is theorem~\ref{Comeon}, in case $p=2$.

\bigskip

(b)
Let $i=p-1$.
The first equivalence is then a particular
case of proposition~\ref{subcat}.
To see the second equivalence recall from theorem~\ref{Comeon},
that there is a $k{\bf B}_{\rho,w} ^\Sigma$-module $kM$, 
with composition factors in
${\cal L}(p-1)$, such that
$k{\bf B}_{\rho,w} ^\Sigma/Ann(kM)$ is Morita equivalent to
the Schur algebra ${\cal S}(w,w)$.

We have $\alpha_{p-1} \subseteq Ann(kM)$, by the first part of the
proposition. 
A dimension count yields an isomorphism between 
$k{\bf B}_{\rho,w} ^\Sigma/Ann(kM)$ and the quotient 
$k{\bf B}_{\rho,w} ^\Sigma/\alpha_{p-1}$.

\bigskip

Let $i=2,...,p-1$.
Theorem~\ref{Ringel} provides a Ringel duality between $kA_i$ and and
$kA_{i-1}$.
In the light of proposition~\ref{subcat}, and the knowledge
that ${\cal S}(w,w)$ is
Ringel self-dual, we discover that 
${\cal L}(i) \cong {\cal S}(w,w)-mod$, for $i=1,...,p-1$. 
$\Box$

\begin{center}
{\large A global-local theorem, and the
James adjustment matrix of a Rock block.} 
\end{center}

In the final sections of this chapter, we prove the following:

\begin{theorem} \label{nilpotent}
There is a Morita equivalence between  
$k{\cal G} b_{\rho,w}^\Sigma$ and a direct sum,
$$\bigoplus_{\substack{a_1,...,a_{p-1} \in \mathbb{Z}_{\geq 0} \\ \sum a_i =w}} 
\left( \bigotimes _{i=1} ^{p-1} {\cal S}(a_i,a_i) \right) ,$$
of tensor products of Schur algebras.

Under this Morita equivalence,
$${\cal L}(a_1,...,a_{p-1})
\cong \left( \bigotimes_{i=1}^{p-1} {\cal S}(a_i,a_i) \right) -mod,$$
and the correspondence of
simple modules is:
$$D^{[\emptyset,...\lambda_{p-3},\lambda_{p-2},\lambda_{p-1}]} 
\leftrightarrow
... \otimes L(\lambda_{p-3}') \otimes L(\lambda_{p-2}) 
\otimes L(\lambda_{p-1} ').$$ 

The James adjustment matrix 
of $k{\bf B}_{\rho,w} ^\Sigma$ is equal to the decomposition
matrix of, 
$$\bigoplus_{\substack{a_1,...,a_{p-1} \in \mathbb{Z}_{\geq 0} \\ \sum a_i =w}} 
\left( \bigotimes _{i=1} ^{p-1} {\cal S}(a_i,a_i) \right) .$$
That is to write,
$$[ D^{[\emptyset,...,\lambda_{p-3}',
\lambda_{p-2},\lambda_{p-1}']} 
_q :
D^{[\emptyset,...,\mu_{p-3}',\mu_{p-2},\mu_{p-1}']} ] =$$
$$= \left\{
\begin{array}{ll}
\prod_{i=1} ^{p-1} [\Delta(\lambda_i) : L(\mu_i)], &
\textrm{ if } |\lambda_i|=|\mu_i|, i=1,...,p-1 \\
0, & {otherwise}
\end{array}
\right.$$
\end{theorem}

Here, we set ${\cal S}(0,0) \cong k$.

\begin{remark}
When $w<p$, the quotient $k{\cal G} b_{\rho,w}^\Sigma$ 
is merely the quotient of
$k{\bf B}_{\rho,w} ^\Sigma$ by its radical.
We are therefore far from the strength of theorem~\ref{structsymm}. 
\end{remark}

Our proof of theorem~\ref{nilpotent}
is inductive, on one hand
applying the Ringel dualities of chapter 6, 
and on the other generalising the theory of chapter 5.

\bigskip

There falls an elegant
description of the decomposition
matrix of $k{\bf B}_{\rho,w} ^\Sigma$, in terms of Littlewood-Richardson
coefficients, and decomposition matrices of 
Schur algebras which are bounded in degree by $w$.

\begin{corollary}
The decomposition matrix of $k{\bf B}_{\rho,w} ^\Sigma$ 
is equal to the matrix product,
$$Dec_{<t-q>} \left( K[t]_{(t-q)}b_{\emptyset,1}^{{\cal H}_t} \wr \Sigma_w 
\right) \times$$
$$Dec  \left(
\bigoplus_{\substack{a_1,...,a_{p-1} \in \mathbb{Z}_{\geq 0} \\ \sum a_i =w}} 
\left( \bigotimes _{i=1} ^{p-1} {\cal S}(a_i,a_i) \right)
\right).$$
Here, $Dec_{<t-q>} \left(
K[t]_{(t-q)} b_{\emptyset,1}^{{\cal H}_t} \wr \Sigma_w) \right)$ 
is the decomposition matrix of
a wreath product of the principal block of $K[t]_{(t-q)} {\cal H}_t(\Sigma_p)$, 
with $\Sigma_w$.  
Formulae for the entries in this matrix are
given in terms of
Littlewood-Richardson coefficients (theorem~\ref{stolen}).
$\Box$ 
\end{corollary}

Note that
$kA_{(a_1,...,a_{p-1})}$ is the reduction modulo $p$ of
${\cal O}A_{(a_1,...,a_{p-1}), q}$.
The proof of theorem~\ref{nilpotent} 
rests upon the following proposition:

\begin{proposition} \label{Polly}
Let $a_i$ be natural numbers such that
$\sum_{i=1}^{p-1} a_i = w$.
Then there exists an ideal 
${\cal N}_{a_1,...,a_{p-1}}$ \index{${\cal N}_{a_1,...,a_{p-1}}$}
in $k{\bf B}_{\rho,w} ^\Sigma$,
such that,

(a) ${\cal N}_{a_1,...,a_{p-1}}$ is equal to the reduction modulo $p$
    of ${\cal ON}_{a_1,...,a_{p-1},q}$.

\bigskip

(b) There are isomorphisms, 
$$k{\bf B}_{\rho,w} ^\Sigma/{\cal N}_{a_1,...,a_{p-1}}
\cong A_{(a_1,...,a_{p-1})} \cong B_{(a_1,...,a_{p-1})}.$$

(c) $k{\bf B}_{\rho,w} ^\Sigma/{\cal N}_{a_1,...,a_{p-1}}$ is Morita equivalent to 
$$\bigotimes_{i=1} ^{p-1} {\cal S}(a_i,a_i).$$

(d) Under the Morita equivalence of (3), the simple module 
$$...\otimes L(\lambda_{p-3}') \otimes L(\lambda_{p-2}) 
\otimes L(\lambda_{p-1} ')$$
for
$\bigotimes_{i=1} ^{p-1} {\cal S}(a_i,a_i)$,
corresponds to the simple $k{\bf B}_{\rho,w} ^\Sigma$-module
$D^{[\emptyset,...,\lambda_{p-3},\lambda_{p-2},\lambda_{p-1}]}$ in 
${\cal L}(a_1,...,a_{p-1})$.

\bigskip

(e) The decomposition matrix of the ${\cal O}$-algebra
    $A_{(a_1,...,a_{p-1})}$ is equal to the 
    decomposition matrix of the quasi-hereditary algebra,
    $\bigotimes_{i=1}^{p-1} {\cal S}(a_i,a_i)$.
\end{proposition}

It is the concern of the last section of this chapter to prove 
proposition~\ref{Polly}.
First, let us give
\emph{a proof of theorem~\ref{nilpotent} 
from proposition~\ref{Polly}:}

\bigskip

Let 
$${\cal N} = \bigcap_{\substack{a_1,...,a_{p-1} \in \mathbb{Z}_{\geq 0} \\ \sum a_i =w}} {\cal N}_{a_1,...,a_{p-1}},$$ \index{${\cal N}$}
$$\Omega = \{ (a_i)_{i=1,...,p-1} | a_i \in \mathbb{Z}_{\geq 0},
\sum_{i=1} ^{p-1} a_i = w \}.$$
Then, for $\alpha \in \Omega$, we know that
$k{\bf B}_{\rho,w} ^\Sigma/(( \bigcap_{\omega \in \Omega-\alpha} 
{\cal N}_\omega) + {\cal N}_\alpha)-mod$ is empty, by theorem~\ref{Polly}(4) and lemma~\ref{triv}. 
We deduce that $k{\bf B}_{\rho,w} ^\Sigma =
(( \bigcap_{\omega \in \Omega-\alpha} 
{\cal N}_\omega) + {\cal N}_\alpha)$.
By linear algebra, it follows that,
$$k{\bf B}_{\rho,w} ^\Sigma/{\cal N} \cong \bigoplus_{\substack{a_1,...,a_{p-1} \in \mathbb{Z}_{\geq 0} \\ \sum a_i =w}} 
k{\bf B}_{\rho,w} ^\Sigma/{\cal N}_{a_1,...,a_{p-1}}.$$
{\noindent
Note that proposition~\ref{Polly}(1) implies that
${\cal N}$ is equal to the $p$-modular reduction of
$$\bigcap_{\substack{a_1,...,a_{p-1} \in \mathbb{Z}_{\geq 0} \\ \sum a_i =w}} {\cal ON}_{a_1,...,a_{p-1}, q}
= b_{\rho,w}^{{\cal H}_q}.{\cal ON}_q.$$
Thus, 
$k{\bf B}_{\rho,w} ^\Sigma/{\cal N}$ is isomorphic to
$k{\cal G} b_{\rho,w}^\Sigma.$}

These isomorphisms, along with proposition~\ref{Polly}(3),
complete the proof of the Morita equivalence of
theorem~\ref{nilpotent}, and the correspondence between simple modules
under this Morita equivalence.

To see the correspondence between decomposition numbers, 
first note that the
decomposition matrix of ${\cal G}b_{\rho,w}^\Sigma$ 
is equal to the decomposition matrix of
$$\bigoplus_{\substack{a_1,...,a_{p-1} \in \mathbb{Z}_{\geq 0} \\ \sum a_i =w}} 
{\cal O}A_{(a_1,...,a_{p-1}),q},$$ 
by proposition~\ref{anima}.
Secondly, note that the decomposition matrix of 
${\cal O}A_{(a_1,...,a_{p-1}),q}$ is equal 
to the decomposition matrix of 
${\cal O}A_{(a_1,...,a_{p-1})}$, since both of these
algebras are semisimple over $K$.
Thirdly, note that the decomposition matrix of 
${\cal O}A_{(a_1,...,a_{p-1})}$
is equal to the decomposition matrix of
$\bigotimes_{i=1} ^{p-1} {\cal S}(a_i,a_i)$,
by proposition~\ref{Polly}(5).
$\Box$

\begin{center}
{\large Induction.}
\end{center}

The intent of this section is to convince the Reader
of the truth of proposition~\ref{Polly}.
Let $w$ be a natural number. We 
assume proposition~\ref{Polly} is proven for Rock blocks of weight
strictly less than $w$, and deduce the same result for a Rock block
$k{\bf B}_{\rho,w} ^\Sigma$ of weight $w$. 

\bigskip

Let $a_{p-1}$ be a natural number, $0< a_{p-1} \leq w$.
Let $u=n-a_{p-1}p$.

Let $L_{a_{p-1}} = \Sigma_p^1 \times ... \times \Sigma_p^{a_{p-1}} \times 
\Sigma_u^0 \leq \Sigma_v$, \index{$L_{a_{p-1}}$}
where $\Sigma_p^i = Sym \{(i-1)p+1,...,ip \}$, 
and $\Sigma_u^0 = Sym \{a_{p-1}p+1,...,wp+r \}$ \index{$\Sigma_u^0$}.

Let $e_{a_{p-1}}$ \index{$e_{a_{p-1}}$}
be an idempotent of $k\Sigma_v$, defined to be  
the product of block idempotents
of $L_i$ with cores 
$\emptyset ,..., \emptyset , \rho$, for $i=0,..,a_{p-1}$.

Let $\zeta_{a_{p-1}} = \sum_{x \in \Sigma_p^1 \times ... \times \Sigma_p^{a_{p-1}}}x$.
\index{$\zeta_{a_{p-1}}$}
Let $kM_{a_{p-1}} = {\cal O}\Sigma_v e_{a_{p-1}} \zeta_{a_{p-1}}$.
\index{$kM_{a_{p-1}}$}

\bigskip

Suppose that 
$a_1,...,a_{p-1}$ are natural numbers, whose sum is $w$.
By our inductive assumption,
there is an ideal ${\cal N}_{a_1,...,a_{p-2},0}$ of $k{\bf B}_{\rho,w-a_{p-1}}^\Sigma$,
such that $k{\bf B}_{\rho,w-a_{p-1}}^\Sigma/{\cal N}_{a_1,...,a_{p-2},0}$ 
is Morita equivalent
to $\bigotimes_{i=1} ^{p-2} {\cal S}(a_i,a_i)$, and
whose simple modules are
the simple objects of
${\cal L}(a_1,...,a_{p-2},0)$.

\bigskip

Let $kM_{(a_1,....,a_{p-1})}$ \index{$kM_{(a_1,....,a_{p-1})}$}
be equal to the quotient
$(kM_{a_{p-1}} / kM_{a_{p-1}} {\cal N}_{a_1,...,a_{p-2},0})$.
Since $kM_{a_{p-1}}$ is projective as a right $k\Sigma_u^0$-module, we have 
$$kM_{(a_1,...,a_{p-1})} \cong
kM_{a_{p-1}} \bigotimes_{\Sigma_u^0} \left( k{\bf B}_{\rho,w-a_{p-1}}^\Sigma / 
{\cal N}_{a_1,...,a_{p-2},0} \right).$$
Therefore $kM_{(a_1,...,a_{p-1})}$ is a $k\Sigma_v$-$k\Sigma_u^0$-bimodule.

\bigskip
The following proposition, and its proof, 
generalise theorem~\ref{Comeon}.

\begin{proposition} \label{Comeon2} 

(a) Consider the $k{\bf B}_{\rho,w} ^\Sigma$-module $kM_{(a_1,...,a_{p-1})}$.
Its endomorphism ring $k{\cal E}_{(a_1,...,a_{p-1})}$ 
\index{$k{\cal E}_{(a_1,...,a_{p-1})}$}
is isomorphic to
$k{\bf B}_{\rho,w-a_{p-1}} ^\Sigma/{\cal N}_{a_1,...,a_{p-2},0} \otimes k\Sigma_{a_{p-1}}$.

\bigskip

(b) The quotient of $k{\bf B}_{\rho,w} ^\Sigma$, by the annihilator 
$kI_{(a_1,...,a_{p-1})}$ \index{$kI_{(a_1,...,a_{p-1})}$} 
of $kM_{(a_1,...,a_{p-1})}$,  
is Morita equivalent to the tensor product, 
$$\bigotimes_{i=1} ^{p-1} {\cal S}(a_i,a_i),$$
of Schur algebras.

\bigskip

(c) The 
$\left( \bigotimes_{i=1} ^{p-1} {\cal S}(a_i,a_i) \right)$-
$\left( \bigotimes_{i=1} ^{p-2} {\cal S}(a_i,a_i) \right) \otimes k\Sigma_{a_{p-1}}$-
bimodule which corresponds 
via Morita equivalence to the $k{\bf B}_{\rho,w} ^\Sigma/kI_{(a_1,...,a_{p-1})}$-
$k{\cal E}_{(a_1,...,a_{p-1})}$ bimodule 
$kM_{(a_1,...,a_{p-1})}$ is isomorphic to,
$$\left( \bigotimes_{i=1} ^{p-2} {\cal S}(a_i,a_i) \right)
\otimes E^{\otimes a_{p-1} \#}.$$ 

(d) Under the Morita equivalence between 
$k{\bf B}_{\rho,w} ^\Sigma/ kI_{(a_1,...,a_{p-1})}$ and the tensor product
$\bigotimes_{i=1} ^{p-1} {\cal S}(a_i,a_i)$,
the correspondence between
simple modules is:
$$D^{[\emptyset,...\lambda_{p-3},\lambda_{p-2},\lambda_{p-1}]} 
\leftrightarrow
... \otimes L(\lambda_{p-3}') \otimes L(\lambda_{p-2}) 
\otimes L(\lambda_{p-1} '). \Box$$ 
\end{proposition}

Let 
$KM_{(a_1,...,a_{p-1}),q}$ \index{$KM_{(a_1,...,a_{p-1}),q}$}
be the $q$-analogue of
$kM_{(a_1,...,a_{p-1})}$.

The proposition below 
is a $q$-analogue of proposition \ref{Comeon}, valid in characteristic
zero.

\begin{proposition} \label{backdoor}
Let $K$ be a splitting field for ${\cal H}_q(\Sigma_v)$.
Let $a_1,...,a_{p-1}$ be natural numbers whose sum is $w$.

(a) $KA_{(a_1,...,a_{p-1}), q}$ is a semisimple algebra.

\bigskip

(b) $KM_{(a_1,...,a_{p-2},a_{p-1}),q}$ is a semisimple 
$K{\bf B}_{\rho,w} ^{{\cal H}_q}$-module.
Its endomorphism ring $K{\cal E}_{(a_1,...,a_{p-1}),q}$ 
\index{$K{\cal E}_{(a_1,...,a_{p-1}),q}$}
is isomorphic to
$\left( K{\bf B}_{\rho,w-a_{p-1}} ^
{{\cal H}_q}/K{\cal N}_{a_1,...,a_{p-2},0,q} \right)
\otimes K\Sigma_{a_{p-1}}$.

\bigskip

(c) The quotient of $K{\bf B}_{\rho,w} ^{{\cal H}_q}$, by the annihilator 
$KI_{(a_1,...,a_{p-1}), q}$ \index{$KI_{(a_1,...,a_{p-1}), q}$}
of $M_{(a_1,...,a_{p-1}), q}$  
is Morita equivalent to the tensor product, 
$\bigotimes_{i=1} ^{p-1} {\cal S}(a_i,a_i)$
of semisimple Schur algebras.

\bigskip

(d) The 
$\left( \bigotimes_{i=1}^{p-1}
{\cal S}(a_i,a_i) \right)$-
$\left( \bigotimes_{i=1} ^{p-2} {\cal S}(a_i,a_i) \right)
\otimes K\Sigma_{a_{p-1}}$-bimodule, which corresponds 
(under Morita equivalence) to 
$KM_{(a_1,...,a_{p-1}), q}$ is isomorphic to
$$\left( \bigotimes_{i=1} ^{p-2} {\cal S}(a_i,a_i) \right)
\otimes E^{\otimes a_{p-1} \#}.$$

(e) The annihilator $KI_{(a_1,...,a_{p-1}),q}$ is precisely 
${\cal N}_{a_1,...,a_{p-1},q}$. 
\end{proposition}
Proof:

Dirichlet's theorem guarantees the existence of infinitely
many prime numbers $l'$, 
such that $l'=1$ (modulo $p$) (i.e. such that $\mathbb{F}_{l'}$ 
contains primitive $p^{th}$ roots of unity).
Let us choose such a prime, such that $w<l'$. \index{$l'$}

Let $q$ be a primitive $p^{th}$ root of unity.
A second application of Dirichlet's theorem provides a prime number
$q'$, such that $q' =q$ (modulo $l'$). 

\bigskip

Let $(K_{l'},{\cal O}_{l'},k_{l'})$ be an $l'$-modular system, 
such that $K_{l'}$ is a splitting 
field for $GL_v(q')$, and such that $q \in {\cal O}_{l'}$.

By theorem~\ref{mainnew}, there is an equivalence between
$k_{l'}b_{\rho,w}^{{\cal H}_q}$, and 
$k_{l'}b_{\emptyset,1}^{{\cal H}_q} \wr \Sigma_w$.
The decomposition matrix of this algebra is equal to the decomposition
matrix of $K_{l'}b_{\rho,w}^{{\cal H}_q}$, by theorem~\ref{stolen}.
Therefore, the James adjustment algebra 
$K{\cal G}b_{\rho,w} ^{{\cal H}_q}$ is semisimple, over
ANY splitting field $K$.

The proposition is now visible, by induction on $w$. 
$\Box$

\bigskip

Here goes the induction.
\emph{Proof of proposition~\ref{Polly}:}

\bigskip

How to define the ideals ${\cal N}_{a_1,...,a_{p-1}}$ ?

\bigskip

Case 1: if $a_1=a_{p-1}=0$, the ideal of
proposition~\ref{Varagn} 
suffices.

Case 2:
If $a_{p-1} \neq 0$, proposition~\ref{Comeon} 
above provides the ideal : set 
${\cal N}_{a_1,...,a_{p-1}} = kI_{(a_1,...,a_{p-1})}$.

Case 3:
If $a_{p-1} = 0$, and $a_1 \neq 0$,
set ${\cal N}_{a_1,....,a_{p-1}} =
{\cal N}_{a_{p-1},...,a_1}^{\#}$.

\bigskip

(a)
Case 1:
Note that, by a $q$-analogue of 
proposition~\ref{Varagn}, we know that 
$I_{(0,a_2,...,a_{p-2},0),q} +
J_{(0,a_2,...,a_{p-2},0),q}$ is contained in
${\cal N}_{0,a_2,...,a_{p-2},0,q}$, over the field $K$, and hence also
over ${\cal O}$.
Thus,
$${\cal N}_{0,a_2,...,a_{p-2},0} \subseteq
{\cal N}_{0,a_2,...,a_{p-2},0,q} \quad (\textrm{modulo } p).$$
However, proposition~\ref{anima} and proposition~\ref{Varagn} 
imply that both quotients,
$K{\bf B}_{\rho,w} ^{{\cal H}_q}/{\cal N}_{0,a_2,...,a_{p-2},0, q}$
and $k{\bf B}_{\rho,w} ^\Sigma/ {\cal N}_{0,a_2,...,a_{p-2},0}$ have the same
dimension, equal to the dimension of
$A_{(a_1,...,a_{p-1})}$.
 
\bigskip

Case 2:
Let $KM_{(a_1,...,a_{p-1}),q}$ be the $q$-analogue of 
$kM_{(a_1,...,a_{p-1})}$.
By proposition~\ref{backdoor}, 
$KM_{(a_1,...,a_{p-1}),q}$ is a semisimple
$K{\cal H}_q(\Sigma_v)$-module, and the quotient
${\cal O}{\bf B}_{\rho,w} ^{{\cal H}_q}/{\cal ON}_{a_1,...,a_{p-1},q}$
surjects onto $k{\bf B}_{\rho,w} ^\Sigma/{\cal N}_{a_1,...,a_{p-1}}$.

Furthermore, on
writing $k{\cal E}_{(a_1,...,a_{p-1})}$ (respectively $K{\cal E}_{(a_1,...,a_{p-1}),q}$) 
for the endomorphism ring
of $kM_{(a_1,...,a_{p-1})}$ (respectively $KM_{(a_1,...,a_{p-1}),q}$),
we have
$$k{\bf B}_{\rho,w} ^{\Sigma}/{\cal N}_{a_1,...,a_{p-1}}
= End_{k{\cal E}_{(a_1,...,a_{p-1})}} (kM_{(a_1,...,a_{p-1})}),$$
$$K{\bf B}_{\rho,w} ^{{\cal H}_q}/K{\cal N}_{a_1,...,a_{p-1},q} 
= End_{K{\cal E}_{(a_1,...,a_{p-1}),q}} (KM_{(a_1,...,a_{p-1}),q}).$$
By propositions~\ref{Comeon2} and \ref{backdoor},
these two endomorphism rings have the same dimensions.
Thus,
${\cal N}_{a_1,...,a_{p-1}} =
{\cal ON}_{a_1,...,a_{p-1},q} ($modulo $p)$.

\bigskip

Case 3:
Note that 
$${\cal N}_{a_1,...,a_{p-1}} = {\cal N}_{a_{p-1},...,a_1} ^{\#} =$$
$${\cal N}_{a_{p-1},...,a_1, q} ^{\#} ( \textrm{modulo } p)
={\cal N}_{a_1,...,a_{p-1},q} (\textrm{modulo } p).$$

(b) follows from proposition~\ref{anima} and (a), 
by $p$-modular reduction.

\bigskip

(c)
What is the Morita type of the quotient 
$k{\bf B}_{\rho,w} ^\Sigma/{\cal N}_{a_1,...,a_{p-1}}$ ?

In case 2, $a_{p-1} \neq 0$, we know that
$k{\bf B}_{\rho,w} ^{\Sigma}/{\cal N}_{a_1,....,a_{p-1}}$ is Morita equivalent to
$\bigotimes_{i=1} ^{p-1} {\cal S}(a_i,a_i)$, 
by proposition~\ref{Comeon2}.
Recall that ${\cal S}(n,n)$ is Ringel self-dual for any $n$ (theorem~\ref{donut}).
The Ringel dualities of theorem~\ref{Ringel}, and the isomorphisms of
part (b), show that 
$k{\bf B}_{\rho,w} ^{\Sigma}/{\cal N}_{a_1,....,a_{p-1}}$ is Morita equivalent to
$\bigotimes_{i=1} ^{p-1} {\cal S}(a_i,a_i)$ in general.
\bigskip

(d)
Tracing back through proposition \ref{Comeon2},
along remark \ref{corr}, past theorem \ref{Ringel} to proposition 
\ref{Rowena},
the correspondence between simple modules is visible.

\bigskip

(e)  
The decomposition matrix 
of $A_{(a_1,...,a_{p-1})}$ as an 
${\cal O}$-algebra is equal to the decomposition matrix 
$([\Delta(\lambda) : L(\mu)])$ of
$kA_{(a_1,...,a_{p-1})}$ as a quasi-hereditary algebra.
(e) is now clear from (b) and (c).
$\Box$

\newpage
\begin{center}
{\bf \large 
Chapter VIII

Doubles, Schur super-bialgebras, and Rock blocks of Hecke algebras.}
\end{center}

J. Alperin's weight conjecture \cite{Alperin}, 
now along with a number of examples
(eg. theorem~\ref{Swd}, theorem~\ref{nilpotent}, \cite{DeVisscher}),
suggests that the category of modules for a finite group algebra,
over a field of prime characteristic, should 
resemble a highest weight category.
However, group algebras are symmetric algebras, and therefore far from quasi-hereditary.
This chapter presents a conjectural resolution
to this problem, for symmetric groups.

Indeed, we associate 
symmetric associative algebras to certain
bialgebras, via a double construction (theorem~\ref{symmetric}).
To any super-algebra, we then assign a ``Schur super-bialgebra''.
From the algebra of $n \times n$ matrices, concentrated in parity zero, 
we thus recover the classical Schur bialgebra, ${\cal S}(n)$.
Applying the aforementioned 
double construction to certain Schur super-bialgebras, 
which correspond to quivers of type $A$,  
we reveal symmetric algebras which
should be Morita equivalent to Rock blocks for Hecke algebras
(conjecture~\ref{Rock}).

\begin{center}
{\large A double construction.}
\end{center}
 
Let $k$ be a field.
Let $B$ be an
bialgebra over $k$, endowed with a $k$-endomorphism $\sigma$,
which is an algebra
anti-homomorphism, and a coalgebra anti-homomorphism.
Suppose that $B$ is graded, with finite dimensional graded pieces.
Let $B^*$ be the graded dual of $B$ \index{$B$, $B^*$}.
Then $B^*$ is a bialgebra, whose product is dual to the coproduct on $B$, and whose coproduct is dual to the product on $B$.

Let us write comultiplication as
$\Delta(x) = \sum x_{(1)} \otimes x_{(2)}$. \index{$\Delta(x)$}

\begin{theorem} \label{symmetric}
The tensor product $D(B) = B \otimes B^*$ \index{$D(B)$}
is a $k$-algebra, with
associative product given by,
$$(a \otimes \alpha).(b \otimes \beta) =
\sum a_{(2)} b_{(1)} \otimes \beta_{(2)}  \alpha_{(1)}
< a_{(1)} ^\sigma, \beta_{(1)} > 
< \alpha_{(2)}, b_{(2)} ^\sigma >.$$
 
Furthermore, $D(B)$ possesses a
symmetric associative bilinear form,
$$< a \otimes \alpha , b \otimes \beta>
= < a^\sigma, \beta><  \alpha, b^\sigma>.$$
Therefore, if $\sigma$ is invertible, 
then $D(B)$ is a symmetric algebra.
 
So long as $B$ is cocommutative, there are algebra homomorphisms,
$$\Delta_l: D(B) \rightarrow D(B) \otimes B,$$ \index{$\Delta_l$}
$$\Delta_l: a \otimes \alpha \mapsto 
\sum a_{(1)} \otimes \alpha \otimes a_{(2)},$$
$$\Delta_r: D(B) \rightarrow B \otimes D(B)$$ \index{$\Delta_r$}
$$\Delta_r: a \otimes \alpha \mapsto 
\sum a_{(1)} \otimes a_{(2)} \otimes \alpha.$$
\end{theorem}
Beneath is a picture of the product $a \otimes \alpha$ and $b \otimes \beta$ in $D(B)$.
We discovered this product, upon studying the group algebra of the principal block of $\Sigma_5$, in characteristic two.

\bigskip

\hspace{2.8cm} $a$ 
\hspace{1.3cm} $\alpha$
\hspace{1.3cm} $b$
\hspace{1.2cm} $\beta$

\begin{center}

\epsfig{figure=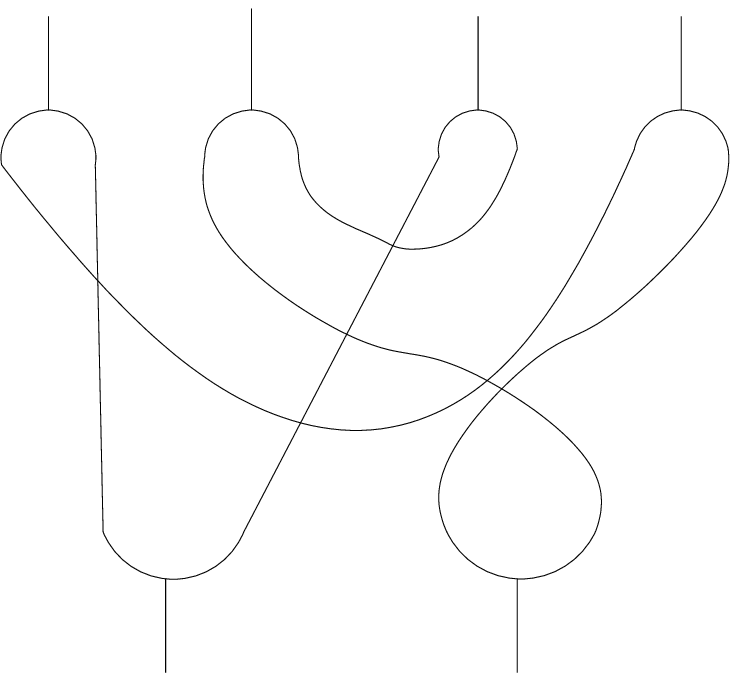, scale=0.8, angle=0} 

\end{center}

Proof:

We first check associativity:

\bigskip

$((a \otimes \alpha).(b \otimes \beta)).(c \otimes \gamma)$

$=(\sum a_{(2)} b_{(1)} \otimes  \beta_{(2)} \alpha_{(1)}
<a_{(1)} ^\sigma, \beta_{(1)}>
< \alpha_{(2)}, b_{(2)} ^\sigma> ).(c \otimes \gamma)$

$= \sum (a_{(2)} b_{(1)})_{(2)} c_{(1)} \otimes 
(\beta_{(2)} \alpha_{(1)})_{(1)} \gamma_{(2)}
<a_{(1)} ^\sigma, \beta_{(1)}>
<\alpha_{(2)}, b_{(2)} ^\sigma>$

$\hspace{2cm}
<(a_{(2)}b_{(1)})_{(1)} ^\sigma, \gamma_{(1)}>
<(\beta_{(2)} \alpha_{(1)})_{(2)}, c_{(2)} ^\sigma>$

$= \sum a_{(3)} b_{(2)} c_{(1)} \otimes
\gamma_{(2)} \beta_{(2)} \alpha_{(1)} 
<a_{(1)} ^\sigma, \beta_{(1)}>
<\alpha_{(3)}, b_{(3)} ^\sigma>$

$\hspace{2cm}
<b_{(1)} ^\sigma a_{(2)} ^\sigma, \gamma_{(1)}>
<\beta_{(3)} \alpha_{(2)}, c_{(2)} ^\sigma >$

$= \sum a_{(3)} b_{(2)} c_{(1)} \otimes
 \gamma_{(2)} \beta_{(2)} \alpha_{(1)}
<a_{(1)} ^\sigma, \beta_{(1)}>
< \alpha_{(3)}, b_{(3)} ^\sigma>$

$\hspace{2cm}
<b_{(1)} ^\sigma, \gamma_{(1)(1)}>
<a_{(2)} ^\sigma, \gamma_{(1)(2)}>
< \alpha_{(2)}, c_{(2)(1)} ^\sigma>
< \beta_{(3)}, c_{(2)(2)} ^\sigma>$

$= \sum a_{(3)} b_{(2)} c_{(1)} \otimes
\gamma_{(3)} \beta_{(2)} \alpha_{(1)} 
<a_{(1)} ^\sigma, \beta_{(1)}>
<\alpha_{(3)}, b_{(3)} ^\sigma>$

$\hspace{2cm}
<b_{(1)} ^\sigma, \gamma_{(1)}>
<a_{(2)} ^\sigma, \gamma_{(2)}>
<\alpha_{(2)}, c_{(2)} ^\sigma>
<\beta_{(3)}, c_{(3)} ^\sigma>$.

\bigskip

This final symmetric expression may similarly
be shown to equal 
$(a \otimes \alpha).((b \otimes \beta).(c \otimes \gamma))$.
Associativity is proven !
\bigskip

Now suppose that $\sigma$ is invertible.
There can be little doubt of the symmetry, nor the non-degeneracy of
the bilinear form we have defined on $D(B)$.
What about associativity ?

\bigskip

$< (a \otimes \alpha).(b \otimes \beta), (c \otimes \gamma) >$

$= \sum 
< a_{(2)} b_{(1)} \otimes \beta_{(2)}  \alpha_{(1)} , c \otimes \gamma >
<a_{(1)} ^\sigma, \beta_{(1)}>
< \alpha_{(2)}, b_{(2)} ^\sigma>$

$= \sum 
< b_{(1)}^\sigma a_{(2)} ^\sigma , \gamma >
< \beta_{(2)} \alpha_{(1)},  c^\sigma >
<a_{(1)} ^\sigma, \beta_{(1)}>
<\alpha_{(2)}, b_{(2)} ^\sigma>$

$=\sum 
< b_{(1)}^\sigma , \gamma_{(1)} >
< a_{(2)} ^\sigma , \gamma_{(2)} >
<  \alpha_{(1)}, c_{(1)} ^\sigma >$

$\hspace{2cm}
< \beta_{(2)}, c_{(2)} ^\sigma >
<a_{(1)} ^\sigma, \beta_{(1)}>
<\alpha_{(2)}, b_{(2)} ^\sigma>$.

\bigskip

This final symmetric expression may similarly
be shown to equal 
$< (a \otimes \alpha) , (b \otimes \beta).(c \otimes \gamma) >$.
Associativity of $<,>$ is proven !

\bigskip
The last check we make is that $\Delta_r$ is an algebra homomorphism,
so long as $B$ is cocommutative
(a similar calculation can be written down for $\Delta_l$):

\bigskip

$\Delta_r( (a \otimes \alpha).(b \otimes \beta))$

$= \sum (a_{(2)} b_{(1)})_{(1)}
\otimes (a_{(2)} b_{(1)})_{(2)} \otimes  \beta_{(2)} \alpha_{(1)}
< a_{(1)} ^\sigma , \beta_{(1)} >
< \alpha_{(2)}, b_{(2)} ^\sigma >$

$= \sum a_{(2)} b_{(1)}  
\otimes a_{(3)} b_{(2)} \otimes  \beta_{(2)} \alpha_{(1)}
< a_{(1)} ^\sigma , \beta_{(1)} >
<  \alpha_{(2)}, b_{(3)} ^\sigma >$

$= \sum a_{(1)} b_{(1)}  
\otimes a_{(3)} b_{(2)} \otimes \beta_{(2)} \alpha_{(1)} 
< a_{(2)} ^\sigma , \beta_{(1)} >
< \alpha_{(2)}, b_{(3)} ^\sigma >$

$= \sum a_{(1)} b_{(1)}  
\otimes a_{(2)(2)} b_{(2)(1)} \otimes  \beta_{(2)} \alpha_{(1)}
< a_{(2)(1)} ^\sigma , \beta_{(1)} >
< \alpha_{(2)}, b_{(2)(2)} ^\sigma >$

$= (a_{(1)} \otimes a_{(2)} \otimes \alpha).
(b_{(1)} \otimes 1 \otimes b_{(2)} \otimes \beta)$

$= \Delta_r(a \otimes \alpha). \Delta_r(b \otimes \beta)$.
$\Box$

\begin{remark}
Since $B$ possesses an algebra anti-automorphism $\sigma$, the dual 
of a left/right $B$-module may be given the structure of a left/right 
$B$-module as well.
The left/right regular action of $B$ on itself, 
may thus be dualised to define a left/right action of $B$ on $B^*$.
We obtain a simpler expression,
$$(a \otimes \alpha).(b \otimes \beta) =
\sum a_{(2)} b_{(1)} \otimes (a_{(1)} \circ \beta ) 
(\alpha \circ b_{(2)}),$$
for the associative product on $D(B)$.
\end{remark}

\begin{remark}
When $B$ is cocommutative, $\Delta_l$ and $\Delta_r$ both give $D(B)$ the structure of a $B$-comodule.
The coproducts $\Delta_l$ and $\Delta_r$ satisfy the following property:

Let $M$ be a $B$-module, and let $N$ be a $D(B)$-module.
The $D(B)$-module $M \otimes N$ (formed via $\Delta_r$) is
isomorphic to the $D(B)$-module $N \otimes M$ (formed via
$\Delta_l$).
\end{remark}

In the examples of this article, 
we find ourselves in the situation of the following lemma.
Its proof is a routine check.

\begin{lemma} \label{setup}
Suppose that $B = \oplus_{r \in \mathbb{Z}_+} B(r)$ is a bialgebra, 
which is direct sum of finite   
dimensional pieces, $B(r)$. 
Suppose that $B$ possesses a degree preserving
$k$-endomorphism $\sigma$,
which is an algebra antiautomorphism, and a coalgebra automorphism.
We write $B^*$ for the graded dual,
$\oplus_{r \in \mathbb{Z}_+} B(r)^*$, of $B$.
Suppose further, that

i. $B(r)$ is a subalgebra of $B$, for $r \in \mathbb{Z}_+$.

ii. $B(0) cong k$, 
and  the projection and embedding maps between $B(0)$ and $B$ give $B$ the structure of an augmented coalgebra.

iii. $B = \oplus_{r \in \mathbb{Z}_+} B(r)$ is a graded coalgebra. 
Thus,
$$\Delta: B(r) \rightarrow \bigoplus_{d=0} ^r B(r-d) \otimes B(d).$$ 
Then the degree $r$ part of $D = D(B)$,
$$D(r) = \bigoplus_{d=0} ^r B(r-d) \otimes B^*(d),$$
\index{$D(r)$}
is a finite-dimensional, graded, symmetric algebra summand of 
$D$, where $B(r-d) \otimes B^*(d)$ 
is given degree $d$.

The ideal,
$${\cal N}(r) =
\bigoplus_{d=1} ^r B(r-d) \otimes B^*(d),$$
\index{${\cal N}(r)$}
of $D(r)$ is nilpotent.

The quotient $D(r) / {\cal N}(r)$ is isomorphic to the
degree zero part $D^0(r) = B(r)$, of $D(r)$.

Irreducible $D(r)$-modules are in natural correspondence
with irreducible $B(r)$-modules.

In this way, $D$ is a graded associative algebra, whose degree zero part
is isomorphic to $B$, as an algebra. 

Upon writing ${\cal N}$ \index{${\cal N}$} for the ideal 
$\oplus_{r \in \mathbb{Z}_+} {\cal N}(r)$ of $D$, a splitting of 
the natural algebra monomorphism $B \rightarrow D$ 
becomes visible:
$$D \rightarrow D / {\cal N} \cong B.$$

The degree $d$ part
$B(r-d) \otimes B^*(d)$ of $D(r)$
inherits a $D^0(r)$-$D^0(r)$-bimodule structure from $D(r)$,
for $d=0,...,r$.
This is nothing but
the natural $B(r)$-$B(r)$-bimodule structure on 
$B(r-d) \otimes B^*(d)$.
$\Box$
\end{lemma}

We require a
super- generalisation of theorem~\ref{symmetric}.

\bigskip

Suppose that $B$ is a super-bialgebra.
Thus, $B$ is $\mathbb{Z}/2$-graded algebra and coalgebra, 
so that the product $m$, and the coproduct
$\Delta$, preserve the grading:
$$m: B^i \otimes B^j \rightarrow B^{i+j},$$
$$\Delta: B^k \rightarrow \bigoplus_{i+j=k} B^i \otimes B^j,
\textrm{ for } i,j,k \in \mathbb{Z}/2.$$
In addition,
$$\Delta(a.b) = \sum (-1)^{|a_{(2)}||b_{(1)}|} 
a_{(1)}. b_{(1)} \otimes a_{(2)}.b_{(2)}.$$
Suppose that $B$ is endowed with 
a parity preserving 
endomorphism $\sigma$, which is a coalgebra anti-automorphism, 
and an algebra anti-automorphism.
\begin{theorem} \label{mod}
The tensor product $D(B) = B \otimes B^*$ is a super-algebra, 
with associative product given by,
$$(a \otimes \alpha).(b \otimes \beta) =$$
$$\sum
(-1)^{s(a, \alpha, b, \beta)}
a_{(2)} b_{(1)} \otimes  \beta_{(2)} \alpha_{(1)}
< a_{(1)} ^\sigma, \beta_{(1)} > 
<  \alpha_{(2)}, b_{(2)} ^\sigma>,$$
where
$$s(a, \alpha, b, \beta) =
|a_{(1)}|(|a_{(2)}| + |b_{(1)}|) + |b_{(1)}||\alpha| + 
|\alpha_{(1)}| |\beta|,$$
and $\mathbb{Z}/2$-grading, given by
$$|a \otimes \alpha| = |a| + |\alpha|.$$
In fact, $D(B)$ is endowed with a symmetric associative bilinear form,
$$< a \otimes \alpha , b \otimes \beta> =
< a^\sigma, \beta>< \alpha, b^\sigma>.$$
Therefore, if $\sigma$ is invertible, $D(B)$ is a symmetric super-algebra.
\end{theorem}
Proof:

Write out the proof of theorem~\ref{symmetric} diagrammatically,
rather than algebraically (thus, a variable should be represented by a 
string, a product by the fraying of a string into two parts, a coproduct
by the joining of two strings together, etc.).
To generalise this proof to the super- situation, 
we need only introduce the sign $(-1)^{|a||b|}$ whenever two strings 
(corresponding to variables $a$,$b$, in degrees $|a|, |b|$) cross.

The sign allocated to our product diagram is $-1$, raised to the power
$$|a_{(1)}|(|a_{(2)}| + |b_{(1)}|) + |b_{(1)}|(|\alpha_{(1)}| + | \alpha_{(2)}|) + 
|\alpha_{(1)}(|\beta_{(1)}| + |\beta_{(2)}|).$$
A slightly simpler expression is $(-1)^{s(a,\alpha,b,\beta)}$.
$\Box$

\begin{remark} The $\mathbb{Z}_+$-grading 
on $B$ of lemma~\ref{setup} is necessarily
a different grading from the $\mathbb{Z}/2$-grading on 
$B$ of theorem~\ref{mod}.
Indeed, lemma~\ref{setup} generalises to apply to a
bialgebra which is $\mathbb{Z}_+ \times \mathbb{Z}/2$-graded, so that 
the $\mathbb{Z}$-grading 
is compatible with lemma~\ref{setup}, whilst the 
$\mathbb{Z}/2$-grading is compatible
with theorem~\ref{mod}.
\end{remark}

\begin{center}
{\large Examples.}
\end{center}

\begin{example}
Let $B = {\cal S}(1)$ be the Schur bialgebra 
associated to $GL_1(k)$, with trivial 
antiautomorphism. 
As a coalgebra, 
${\cal S}(1)$ is the graded
dual of the polynomial ring $k[X]$ in one variable.
Each homogeneous component ${\cal S}(1,r)$ is isomorphic (as an algebra)
to a copy of the
field $k$.

Then we are in the context of lemma~\ref{setup}, 
and $D(r)$ is isomorphic to 
the uniserial algebra $k[Y]/(Y^r)$.
$\Box$
\end{example}

\begin{example} \label{BT}
Let $B = B(0) \oplus B(1) = k \oplus T_n(k)$, 
be the direct sum of a copy of the field $k$ 
(in degree zero), and the algebra of $n \times n$ upper triangular matrices
(in degree one).
On writing $\epsilon$ for the unit in $B(0)$ 
(\emph{not a unit for $B$}), 
we see that $B$ is a 
cocommutative bialgebra via the coproduct,
$$\Delta:
x \mapsto x \otimes \epsilon + \epsilon \otimes x, \hspace{1cm} x \in B(1)$$
$$\Delta: \epsilon \mapsto \epsilon \otimes \epsilon.$$
The bialgebra
$B$ possesses an algebra anti-automorphism $\sigma$, 
acting trivially on $B(0)$, but non-trivially on 
$B(1)$ taking $E_{i,j}$ to $E_{n-j+1,n-i+1}$.
This map $\sigma$ is a coalgebra anti-automorphism.

We are in the setup of lemma~\ref{setup}, and
$D(1)$ is isomorphic to the path algebra of the circular quiver with $n$ vertices, and clockwise orientation, modulo the ideal of 
paths of length $\geq n+1$.
This is a uniserial algebra, otherwise known as the Brauer tree 
algebra of a star, with multiplicity one (\cite{Alp}, chapter 5).
$\Box$
\end{example}

\begin{example} 
Let $Q$ be a quiver, without loops or
multiple edges, equipped with an orientation reversing automorphism.
Let $B = k \oplus kQ/I_2$ be the direct sum of a 
copy of the field $k$ 
(in degree zero), and the path algebra $kQ$, modulo the ideal of 
paths of length $\geq 2$ (in degree one).
Just as in example~\ref{BT}, 
$B$ may be given the structure of a bialgebra, equipped 
with an algebra anti-automorphism which is a coalgebra anti-automorphism.

This time,
$D(1)$ is isomorphic to the zigzag algebra (see \cite{Khov}),
whose graph is the underlying graph of $Q$.
If the underlying graph is an ordinary Dynkin diagram of type $A$, 
the zigzag algebra is otherwise known as a linear Brauer tree 
algebra, with multiplicity one.
$\Box$
\end{example}

Brauer tree algebras appear naturally in the block theory of finite groups
with cyclic defect groups.
There appear to be more mysterious 
instances of doubles appearing in
finite group theory: 

\begin{example} \label{Karin}
Let $k$ be a field of characteristic two.
Let $B = {\cal S}(2)$ be the Schur 
bialgebra associated to $GL_2(k)$, with transpose 
antiautomorphism. 
Then $D(2)$ is
Morita equivalent to the Rock block $k{\bf B}_{\rho,2}^\Sigma$,
otherwise known as the principal block of 
$k\Sigma_5$.
This follows from Erdmann's description of the basic algebras
for tame blocks of group algebras \cite{Erd}.
$\Box$ 
\end{example}

\begin{example} \label{Holm}
Let $k$ be a field of characteristic two.
Let $B = {\cal S}(2)$ be the Schur bialgebra 
associated to $GL_2(k)$, with transpose 
antiautomorphism. Let 
$$C = 
(0 \rightarrow D(2) \xi_{(2)} \otimes \xi_{(2)} D(2) \xi_{(1^2)}
\rightarrow D(2) \xi_{(1^2)} \rightarrow 0).$$
Here, the differential is given by the product map in $D(2)$.
Then $C$ is a tilting complex for $D(2)$,  
and its endomorphism ring ${\cal E}$
in the homotopy category is Morita equivalent to $k\Sigma_4$.
This follows from Holm's description \cite{Holm}
of derived equivalences between
tame blocks of group algebras.
$\Box$
\end{example}

\begin{remark} The equivalences of 
examples~\ref{Karin} and~\ref{Holm} both lift to Hecke 
algebras at $-1$, over fields of arbitrary characteristic.
\end{remark}

\begin{center}
{\large Notation.}
\end{center}

Let $V$ be a vector space.

We write $\bigwedge(V)$ for the exterior algebra on $V$, the coinvariants of the signature action of $\times_{r \geq 0} \Sigma_r$ on the tensor algebra 
$T(V) = \bigoplus_{r \geq 0} V^{\otimes r}$.
If $v_1,...,v_n$ is a basis for $V$, then $\{ v_{i_1} \wedge ... \wedge v_{i_r} | i_1<... < i_r, r \geq 0 \}$
is a basis for $\bigwedge(V)$, where $v_{i_1} \wedge ... \wedge v_{i_r}$ is the image in $\bigwedge(V)$ of $v_{i_1} \otimes ... \otimes v_{i_r}$.

We write $\bigvee(V)$ for the invariants of the signature action of $\times_{r \geq 0} \Sigma_r$ on $T(V)$.
If $v_1,...,v_n$ is a basis for $V$, then $\{ v_{i_1} \vee ... \vee v_{i_r} | i_1<... < i_r, r \geq 0 \}$
is a basis for $\bigwedge(V)$, where $v_{i_1} \vee ... \vee v_{i_r}$ is the anti-symmetrisation of $v_{i_1} \otimes ... \otimes v_{i_r}$.

We write ${\cal A}(V)$ for the symmetric algebra on $V$, the coinvariants of the permutation action of $\times_{r \geq 0} \Sigma_r$ on the tensor algebra 
$T(V) = \bigoplus_{r \geq 0} V^{\otimes r}$.
If $v_1,...,v_n$ is a basis for $V$, then $\{ v_{i_1} ... v_{i_r} | i_1 \leq ... \leq i_r, r \geq 0 \}$
is a basis for ${\cal A}(V)$, where $v_{i_1} ... v_{i_r}$ is the image in ${\cal A}(V)$ of $v_{i_1} \otimes ... \otimes v_{i_r}$.

We write ${\cal S}(V)$ for the invariants of the permutation action of $\times_{r \geq 0} \Sigma_r$ on $T(V)$.
If $v_1,...,v_n$ is a basis for $V$, then $\{ v_{i_1} *...* v_{i_r} | i_1 \leq... \leq i_r, r \geq 0 \}$
is a basis for ${\cal S}(V)$, where $v_{i_1} ... v_{i_r}$ is the symmetrisation of $v_{i_1} \otimes ... \otimes v_{i_r}$.

Now suppose $M$ is the algebra of $n \times n$ matrices.
Then we write $\bigwedge(n) = \bigwedge(M^*)$, $\bigvee(n) = \bigvee(M)$, and ${\cal A}(n) = {\cal A}(M^*)$, ${\cal S}(n) = {\cal S}(M)$.
We write $\bigwedge(n,r)$, $\bigvee(n,r)$, ${\cal A}(n,r)$, ${\cal S}(n,r)$ for the $r^{th}$ homogeneous components of these various spaces.
 
\begin{center}
{\large Schur super-bialgebras.}
\end{center}

Let $A,B$ be super-algebras (i.e. $\mathbb{Z}/2$-graded 
associative algebras). 
Their tensor product, $A \otimes B$,
becomes a super-algebra, with parity,
$$| a \otimes b | = |a| + |b|,$$
and super-product,
$$(a \otimes b).(a' \otimes b') = (-1)^{|b||a'|} (aa' \otimes bb').$$

For a super-algebra $A$, we define the super-algebra 
$A^{\otimes r}$ inductively,
to be $A^{\otimes r} = A^{\otimes r-1} \otimes A$, with super-product
as above.

The symmetric group $\Sigma_r$ acts naturally
as parity-preserving automorphisms
on $A^{\otimes r}$.
A simple reflection $(i, i+1) \in \Sigma_r$ acts as:
$$( a_1 \otimes ... \otimes a_r)^{(i, i+1)} =
(-1)^{|a_i||a_{i+1}| }
a_1 \otimes ... \otimes a_{i-1}
\otimes a_{i+1}\otimes a_{i}\otimes a_{i+2} \otimes ... \otimes
a_{r}.$$ 

Let $A$ be a super-algebra.
Let $T(A)$ be the direct sum, \index{$T(A)$}
$$T(A) = \bigoplus_{d \geq 0} A^{\otimes d},$$
of super-algebras.
This algebra becomes a super-bialgebra, with parity,
$$| a_1 \otimes... \otimes a_r | = |a_1| + ... + |a_r|,$$
and coassociative comultiplication,
$$\Delta(a_1 \otimes ... \otimes a_r) =
\sum_{i=0} ^r (a_1 \otimes ... \otimes a_i) \otimes (a_{i+1} 
\otimes ... \otimes a_r).$$

\begin{definition}
Let $A$ be a super-algebra.
The \emph{Schur super-bialgebra} associated to $A$ is the graded
sub-super-bialgebra,
$${\cal S}(A) = \bigoplus_{r \geq 0} {\cal S}(A)(r) =
\bigoplus_{r \geq 0} (A^{ \otimes r})^{\Sigma_r},$$
\index{${\cal S}(A)$}
of $\Sigma_r$-fixpoints on $T(A)$.
\end{definition}

\begin{remark}
In case $A = M = M_n(k)$ is the algebra of $n \times n$ matrices, concentrated 
in parity zero, we recover the classical Schur bialgebra 
${\cal S}(n)$ in this way.
\end{remark}

\begin{remark}
If $A$ is a super-algebra, equipped with a parity-preserving
anti-automorphism $\sigma$,
then $A^{\otimes r}$ may also be equipped with a parity-preserving 
anti-automorphism, 
$$\sigma: a_1 \otimes ... \otimes a_r \mapsto 
\sigma(a_r) \otimes ... \otimes \sigma(a_1).$$
Indeed, this map acts as a coalgebra anti-automorphism on $T(A)$, and
restricts to a coalgebra anti-automorphism on ${\cal S}(A)$.  

Therefore, if $A$ is a super-algebra, equipped with a parity-preserving 
anti-automorphism $\sigma$, then ${\cal S}(A)$ is a super-bialgebra,
equipped with an endomorphism $\sigma$, which is an algebra and 
coalgebra anti-automorphism.
Under such circumstances, we may apply theorem~\ref{mod}, and 
form the double, $D({\cal S}(A))$.
\end{remark}

\begin{center}
{\large Schiver super-bialgebras.}
\end{center}

The example which concerns us most in this booklet, is the special case when
$A$ is Morita equivalent to the
path algebra of a quiver $Q$, modulo the ideal of paths of length $>1$.
We are particularly interested in this when $Q$ is a Dynkin quiver
of type $A$.

\bigskip

Let $Q$ be a quiver \index{$Q$}
(that is, a locally finite, oriented graph).
We note its vertex set $V$, \index{$V$} and its set of edges $E$. \index{$E$}
For any edge $e \in E$, we denote its source $s(e) \in V$, and its tail 
$t(e) \in V$.

Let $P_Q$ be the path algebra of $Q$, modulo the ideal 
of paths of length $>1$.
To any natural number $n$, let us assign the algebra 
$P_Q(n)$, \index{$P_Q, P_Q(n)$}
which is Morita equivalent to $P_Q$, and whose simple modules
all have dimension $n$.  
Thus, $P_Q(n) = End_{P_Q}(P_Q ^{\oplus n})$, and as vector spaces we have
$$P_Q(n) \cong
M^{\oplus V} \oplus  
M^{\oplus E}, $$
where $M$ is the algebra of $n \times n$ matrices over $k$.
The space $P_Q(n)$
is naturally an algebraic affine super-variety, where
paths in $P_Q(n)$ of length $0$ and $1$ are given
parities $0$ and $1$ respectively.

\begin{definition} 
The \emph{Schur quiver super-bialgebra}, or \emph{Schiver super-bialgebra} associated
to $(Q,n)$, is
the Schur super-bialgebra,
${\cal S}_Q(n) = {\cal S}(P_Q(n))$,
\index{${\cal S}_Q(n)$}
associated to $P_Q(n)$. Its graded dual, the ring of functions,
$${\cal A}_Q(n) \cong \left( 
{\cal A}(M^*) \right) ^{\otimes V}
\otimes \left(
\bigwedge( M^* ) 
\right)^{\otimes E},$$
\index{${\cal A}_Q(n)$}
on $P_Q(n)$,
is isomorphic as an algebra,
to a tensor product of symmetric and exterior algebras. 
\end{definition}

Our sole motivation for this 
definition is the apparent emergence of such structure in type A 
representation theory.
But the alternating structure also bears encouraging  
homological consequences.
For example, the super-symmetric aspect of this definition
looks rather becoming, when one
considers Koszul duality for these super-algebras - 
see remark~\ref{Koszul}.

\begin{remark}
When $Q$ has merely one vertex, and no arrows, we recover
the classical Schur bialgebras ${\cal S}(n)$ according to this construction.
In general the tensor product,
$${\cal S}_{V}(n) = ({\cal S}(n)) ^{\otimes V},$$
of Schur algebras, is naturally a sub-bialgebra of 
$${\cal S}_Q(n) \cong
\left( 
{\cal S}(n)  \right)^{\otimes V}
\otimes \left(
\bigvee( n) 
\right)^{\otimes E}.$$
The inclusion map splits as an algebra homomorphism, via
$${\cal S}_Q(n) \rightarrow {\cal S}_Q(n) / {\cal J} \cong 
({\cal S}(n)) ^{\otimes V},$$
where ${\cal J}$ is the direct sum of subspaces,
$$\left( \bigotimes_{v \in V} 
{\cal S}(n, a_v) 
\right) \otimes \left(
\bigotimes_{e \in E} \bigvee (n, b_e) \right)$$
of ${\cal S}_Q(n)$, such that $b_e>0$ for some $e \in E$.
Thus, ${\cal J}$ is a nilpotent ideal of the algebra ${\cal S}_Q(n)$. 
\end{remark}

If $Q$ is a quiver, then the disjoint union
$Q \coprod Q^{op}$ of quivers possesses an obvious orientation-reversing
automorphism, exchanging $Q$ and $Q^{op}$.
Thus, ${\cal S}_{Q \coprod Q^{op}}(n)$ may be
equipped with a $k$-endomorphism $\sigma$, which is an algebra and 
coalgebra anti-automorphism.
Under such circumstances, we may apply theorem~\ref{mod}, and 
form the double, $D({\cal S}_{Q \coprod Q^{op}}(n))$.

\begin{definition}
The \emph{Schiver double} associated to $(Q,n)$, 
is the natural algebra summand,
$${\cal D}_Q(n) = {\cal S}_Q(n) \otimes {\cal A}_{Q^{op}}(n),$$
\index{${\cal D}_Q(n)$, ${\cal D}_Q(n,r)$} 
of the double $D({\cal S}_{Q \coprod Q^{op}}(n))$ corresponding to 
$Q \coprod Q^{op}$.

Thus, ${\cal D}_Q(n)$ is a direct sum,
$\bigoplus_{r \geq 0} {\cal D}_Q(n,r)$ of algebras, where
$${\cal D}_Q(n,r) = \bigoplus_{r_1 + r_2 = r}
{\cal S}_Q(n,r_1) \otimes {\cal A}_{Q^{op}}(n,r_2).$$
\end{definition}

\begin{remark} \label{grading}
Each algebra summand ${\cal D}_Q(n,r)$ is 
$\mathbb{Z}_+ \times \mathbb{Z}_+$-graded, where 
the component,
$$\left( \bigotimes_{v \in V} 
{\cal S}(n, a_v) 
\right) \otimes \left(
\bigotimes_{e \in E} \bigvee (n, b_e) \right)$$
$$\otimes
\left( \bigotimes_{v \in V} 
{\cal A}(n, c_v) 
\right) \otimes \left(
\bigotimes_{e \in E} \bigwedge (n, d_e) \right),$$
is given degree $(\sum_E b_e + \sum_V c_v, \sum_V c_v + \sum_E d_e)$.

\bigskip

Each algebra summand ${\cal D}_Q(n,r)$ is 
$\mathbb{Z}_+$-graded, where 
the component above is given degree 
$(\sum_E b_e + 2\sum_V c_v + \sum_E d_e)$.
We write ${\cal D}_Q^i(n,r)$ for the degree $i$ part 
with respect to this grading. Thus,
$${\cal D}_Q(n,r) = \bigoplus_{i=0} ^{2r} {\cal D}_Q^i(n,r),$$
\index{${\cal D}_Q^i(n,r)$}
as a direct sum of graded pieces.
In degree zero, we have 
$${\cal D}_Q^0(n,r) =  {\cal S}_{V(Q)}(n,r).$$ 
\end{remark}

\begin{center}
{\large Schiver doubles: independence of quiver orientation.}
\end{center}

This section is devoted to a proof of the following result...  

\begin{theorem} \label{invariant}
The Schiver double ${\cal D}_Q(n)$ is
independent of the orientation of $Q$, and
as such, is an invariant of the underlying graph of $Q$. 
\end{theorem}

For a locally finite
graph $\Gamma$, we thus write ${\cal D}_\Gamma(n)$ for
the Schiver double ${\cal D}_Q(n)$, where $Q$ is any orientation of $\Gamma$.

We give a proof  of theorem~\ref{invariant}, in case $Q=A_1$ is the quiver with
two vertices, and one arrow connecting those two vertices. 
To say that the corresponding double is 
independent of orientation, is to say that there is an algebra isomorphism 
between the double corresponding to the quiver,
$\circ$-------------$\rhd \circ$,
and the double corresponding to the quiver,
$\circ \lhd$-------------$\circ$.
We therefore reveal an algebra automorphism of
${\cal D}_{A_1}(n)$.

\bigskip

Theorem~\ref{invariant} follows for a general quiver from the case $Q=A_1$.
To see this, 
first observe that distinct arrows do not interact with one another 
when multiplied in ${\cal D}_Q(n)$. Therefore, if $Q'$ is obtained from $Q$, by the
reversing of an arrow, we have ${\cal D}_{Q'}(n) \cong {\cal D}_{Q}(n)$. 
Secondly note that we may obtain one orientation of $Q$ from another by 
reversing a collection of arrows, and so ${\cal D}_Q(n)$ is indeed independent
of the orientation of $Q$. 

\bigskip

We present a triad of preliminary lemmas.

\bigskip

Let $\bigvee(n) = \bigvee(M)$, and let $\bigwedge(n) = \bigwedge(M^*)$.
\index{$\bigvee(n)$, $\bigwedge(n)$}
Theorem~\ref{donut} implies the following lemma.

\begin{lemma} \label{above}
There is an ${\cal S}(n)$-${\cal S}(n)$-bimodule isomorphism, which
exchanges $\bigvee(n)$ and $\bigwedge(n)$, for $n \geq 0$.

When $n=1$, this isomorphism is defined by the structure of a symmetric algebra on $M$.
$\Box$
\end{lemma}

We write $*$ for either of the inverse homomorphisms which describe the 
isomorphism of lemma~\ref{above}.
We have $(x \wedge y)^* = (x^* \vee y^*)$, and
$(x \vee y)^* = (x^* \wedge y^*)$.

The tranpose anti-automorphism $\sigma$ of ${\cal S}_{A_1}(n)$ maps
$s \otimes \lambda \otimes s'$ to $s'^T \otimes \lambda^T \otimes s^T$. 

\begin{lemma}
The left action of ${\cal S}_{A_1}(n)$ on ${\cal A}_{A_1}(n)$, 
is given by,
$$(s \otimes \lambda \otimes t) \circ (a \otimes \mu \otimes b) =$$
$$\sum (-1)^{|\mu_{(2)}| |\lambda|}
(t_{(2)} \circ a) \otimes (t_{(1)} \circ \mu_{(2)}) \otimes
(\mu_{(1)(2)}.(s \circ b)) < \mu_{(1)(1)}, \lambda^T >.$$

The right action of ${\cal S}_{A_1}(n)$ on ${\cal A}_{A_1}(n)$, 
is given by,
$$(a \otimes \mu \otimes b) \circ (s \otimes \lambda \otimes t) =$$
$$\sum (-1)^{|\mu_{(1)}| |\mu_{(2)}|}
((a \circ t). \mu_{(2)(1)}) \otimes (\mu_{(1)} \circ s_{(2)}) \otimes
(b \circ s_{(1)}) < \mu_{(2)(2)}, \lambda^T >.$$
\end{lemma}
Proof:

We record a calculation for the left action:

\bigskip

$<(s \otimes \lambda \otimes t) \circ (a \otimes \mu \otimes b),
(s' \otimes \lambda' \otimes t')>$

$= <a \otimes \mu \otimes b, (s \otimes \lambda \otimes t)^\sigma 
\circ (s' \otimes \lambda' \otimes t')>$

$= <a \otimes \mu \otimes b, (t^T \otimes \lambda^T \otimes s^T) 
\circ (s' \otimes \lambda' \otimes t')>$

$= \sum <a \otimes \mu \otimes b, 
(t^T)_{(1)} s' \otimes 
(\lambda^T \circ t'_{(1)} ) \vee ((t^T)_{(2)} \circ \lambda')
\otimes s^T t'_{(2)} >$

$= \sum 
(-1)^{|\mu_{(2)}| |\lambda|} <a, t_{(2)}^T s'>
<\mu_{(1)}, (\lambda^T \circ t'_{(1)})>$

\hspace{2cm}
$<\mu_{(2)}, (t_{(1)}^T \circ \lambda')>
<b, s^T t'_{(2)} >$

$= \sum (-1)^{|\mu_{(2)}| |\lambda|}
<t_{(2)} \circ a, s'> < \mu_{(1)(1)}, \lambda^T >
< \mu_{(1)(2)}, t'_{(1)} >$

\hspace{2cm}
$< t_{(1)} \circ \mu_{(2)}, \lambda' > <s \circ b, t'_{(2)} >$

$= \sum (-1)^{|\mu_{(2)}| |\lambda|}
< (t_{(2)} \circ a \otimes t_{(1)} \circ \mu_{(2)}
\otimes \mu_{(1)(2)}.(s \circ b)), (s' \otimes \lambda' \otimes t')>$
 
\hspace{2cm}
$< \mu_{(1)(1)}, \lambda^T >$
$\Box$

\bigskip

Dual to the bilinear form  
$$\epsilon: \bigwedge(n) \otimes \bigvee(n) \rightarrow k, $$ 
there is a natural map 
$$\epsilon^*: k \rightarrow \bigwedge(n) \otimes \bigvee(n).$$
Squeezing the identity map on ${\cal S}(n)$ inside $\epsilon$, we obtain a map
$$\phi:
{\cal S}(n) \rightarrow \bigwedge(n) \otimes {\cal S}(n) \otimes \bigvee(n).$$
The right action of ${\cal S}(n)$ on 
$\bigwedge(n)$ may be formulated as a map,
$$m_1:  \bigwedge(n) \otimes {\cal S}(n) \rightarrow \bigwedge(n).$$
The left action of ${\cal S}(n)$ on 
$\bigvee(n)$ 
may be formulated as a map,
$$m_2: {\cal S}(n) \otimes \bigvee(n) \rightarrow \bigvee(n).$$

\begin{lemma} \label{square}
The diagram
$$\xymatrix{ 
{\cal S}(n) \ar[d]^\phi \ar[r]^\phi  
& \bigwedge(n) \otimes {\cal S}(n) \otimes \bigvee(n) \ar[d]^{T \otimes m_2} \\
\bigwedge(n) \otimes {\cal S}(n) \otimes \bigvee(n) \ar[r]^{m_1 \otimes T} 
& \bigwedge(n) \otimes \bigvee(n)}$$
commutes.
\end{lemma}
Proof:

Let $\{ \xi_{ij} \, | \, i, j = 1,...,n \}$ be a basis for $M$. Let
$\{ X_{ij} \}$ be the basis for $M^*$ which is identified with $\{ \xi_{ij} \}$
via the symmetric structure on $M$. 
The Schur algebra may be given basis, 
whose elements have the form $\xi_{a_1 b_1} ....\xi_{a_r b_r}$.
Let $\Sigma(a,b)$ be the stabilizer in $\Sigma_n$ of the sequence $(a_i,b_i)_{i=1}^n$.
We have,
$$(m_1 \otimes T) \circ \phi (\xi_{a_1 b_1} ....\xi_{a_r b_r})$$ 
$$= (m_1 \otimes 1) \left( 
\sum_{i_k,j_k} X_{i_1 j_1} \wedge...\wedge X_{i_n j_n} \otimes
\xi_{a_1 b_1} ....\xi_{a_r b_r} \otimes
\xi_{j_1 i_1} \vee .... \vee \xi_{j_r i_r} \right) $$
$$= (m_1 \otimes 1) \left( 
\begin{array}{c}
\sum_{i_k,j_k, \sigma \in \Sigma(a,b)} 
X_{i_1 j_1} \wedge...\wedge X_{i_n j_n} \otimes 
\xi_{a_{\sigma 1} b_{\sigma 1}} \otimes.... \otimes \xi_{a_{\sigma r} b_{\sigma r}}  
\otimes \\
\xi_{j_1 i_1} \vee .... \vee \xi_{j_r i_r} 
\end{array}
\right)$$
$$= \sum_{i_k,\sigma \in \Sigma(a,b)}  X_{i_1 b_{\sigma 1}} \wedge...\wedge X_{i_n b_{\sigma n}}
\otimes 
\xi_{a_{\sigma 1} i_1} \vee .... \vee \xi_{a_{\sigma r} i_r}.$$
A similar computation shows that $(T \otimes m_2) \circ \phi$ is equal to the same sum.
$\Box$

\bigskip

The duals of $m_1, m_2$ are maps, 
$$m_1 ^* : \bigvee(n) \rightarrow  \bigvee(n) \otimes {\cal A}(n).$$
$$m_2 ^* : \bigwedge(n) \rightarrow {\cal A}(n) \otimes \bigwedge(n).$$
We write 
$m_1 ^*(\gamma) = \sum \gamma_{(1)} \otimes \gamma_{(2)}$ for 
$\gamma \in \bigvee(n)$, and 
$m_2 ^*(\delta) = \sum \delta_{(1)} \otimes \delta_{(2)}$ for $\delta \in \bigwedge(n)$.

\begin{lemma} \label{help}
For $\alpha \in \bigvee(n), \beta \in \bigwedge(n)$, we have,
$$\sum \alpha_{(2)} < \alpha_{(1)}, \beta^T >
= \sum < \alpha^{*T}, \beta_{(2)} ^*, >  \beta_{(1)} ^*$$
\end{lemma}
Proof: 

By lemma~\ref{square}, the following diagram commutes:
$$\xymatrix{
& {\cal S}(n) \ar[dl]^\phi \ar[dd]_\phi \ar[dr]^{(* \otimes 1 \otimes *) \phi} \\
\bigwedge(n) \otimes {\cal S}(n) \otimes \bigvee(n) \ar[dd]_{m_1 \otimes T} & &
\bigvee(n) \otimes {\cal S}(n) \otimes \bigwedge(n) \ar[dd]^{T \otimes m_2} 
\ar[dl]_{* \otimes 1 \otimes *} \\
& \bigwedge(n) \otimes {\cal S}(n) \otimes \bigvee(n) \ar[dl]^{T \otimes m_2} \\
\bigwedge(n) \otimes \bigvee(n)  \ar[rr]_{* \otimes *}
& & \bigvee(n) \otimes \bigwedge(n) 
}$$
Therefore, the diagram dual to this one commutes.
The two passages from $\bigvee(n) \otimes \bigwedge(n)$ to ${\cal A}(n)$
around the boundary of this dual diagram describe the two sides to the 
formula of the lemma.
$\Box$

\begin{theorem} \label{autm}
There is an involutory
algebra automorphism $\theta$ of ${\cal D}_{A_1}(n)$,
given by,
$$\theta(s \otimes \lambda \otimes t \otimes 
a \otimes \xi \otimes b) = (-1)^{|\lambda||\xi|}
(t \otimes \xi^* \otimes s \otimes b \otimes \lambda^* \otimes a).$$
\end{theorem}
Proof:

\bigskip

$[s \otimes \lambda \otimes t \otimes a \otimes \xi \otimes b].
[u \otimes \mu \otimes v \otimes c \otimes \eta \otimes d]$

$= \sum \pm (s \otimes \lambda \otimes t)_{(2)}.(u \otimes \mu \otimes v)_{(1)}
\otimes [(s \otimes \lambda \otimes t)_{(1)} \circ
(c \otimes \eta \otimes d)] \otimes$

\hspace{2cm}
$[(a \otimes \xi \otimes b) \circ (u \otimes \mu \otimes v)_{(2)}]$

$= \sum \pm (s_{(2)} \otimes \lambda_{(2)} \otimes t_{(2)}).
(u_{(1)} \otimes \mu_{(1)} \otimes v_{(1)})
\otimes$

\hspace{2cm}
$[(s_{(1)} \otimes \lambda_{(1)} \otimes t_{(1)}) \circ
(c \otimes \eta \otimes d)] \otimes$

\hspace{2cm}
$[(a \otimes \xi \otimes b) \circ 
(u_{(2)} \otimes \mu_{(2)} \otimes v_{(2)})]$

$= \sum \pm s_{(2)(1)}u_{(1)} \otimes 
(\lambda_{(2)} \circ v_{(1)(1)}) \vee (s_{(2)(2)} \circ \mu_{(1)})
\otimes t_{(2)} v_{(1)(2)} \otimes$

\hspace{2cm}
$[t_{(1)(2)} \circ c \otimes
t_{(1)(1)} \circ \eta_{(2)} \otimes 
\eta_{(1)(2)}.(s_{(1)} \circ d)].$

\hspace{2cm}
$[(a \circ v_{(2)}). \xi_{(2)(1)}
\otimes \xi_{(1)} \circ u_{(2)(2)}  
\otimes c \circ u_{(2)(1)}]$

\hspace{2cm}
$< \eta_{(1)(1)}, \lambda_{(1)}^T >
< \xi_{(2)(2)}, \mu_{(2)}^T >$

$= \sum \pm s_{(2)}u_{(1)} \otimes 
(\lambda_{(2)} \circ v_{(1)}) \vee (s_{(3)} \circ \mu_{(1)})
\otimes t_{(3)} v_{(2)} \otimes$

\hspace{2cm}
$[t_{(2)} \circ c \otimes
t_{(1)} \circ \eta_{(3)} \otimes 
\eta_{(2)}.(s_{(1)} \circ d)].$

\hspace{2cm}
$[(a \circ v_{(3)}). \xi_{(2)}
\otimes \xi_{(1)} \circ u_{(3)}  
\otimes c \circ u_{(2)}]$

\hspace{2cm}
$< \eta_{(1)}, \lambda_{(1)}^T >
< \xi_{(3)}, \mu_{(2)}^T >$

$= \sum \pm s_{(2)}u_{(1)} \otimes 
(\lambda_{(2)} \circ v_{(1)}) \vee (s_{(3)} \circ \mu_{(1)})
\otimes t_{(3)} v_{(2)} \otimes$

\hspace{2cm}
$(t_{(2)} \circ c). (a \circ v_{(3)}). \xi_{(2)} \otimes
(t_{(1)} \circ \eta_{(3)}) \wedge (\xi_{(1)} \circ u_{(3)})  \otimes$ 

\hspace{2cm}
$\eta_{(2)}.(s_{(1)} \circ d).(c \circ u_{(2)})
< \eta_{(1)}, \lambda_{(1)}^T >
< \xi_{(3)}, \mu_{(2)}^T >.$

\bigskip

\noindent{Thus,}

\bigskip

$\theta([s \otimes \lambda \otimes t \otimes a \otimes \xi \otimes b]).
\theta([u \otimes \mu \otimes v \otimes c \otimes \eta \otimes d])$

$= [t \otimes \xi^* \otimes s \otimes b \otimes \lambda^* \otimes a].
[v \otimes \eta^* \otimes u \otimes d \otimes \mu^* \otimes c]$

$= \sum \pm t_{(2)}v_{(1)} \otimes 
(\xi^* _{(2)} \circ u_{(1)}) \vee (t_{(3)} \circ \eta^* _{(1)})
\otimes s_{(3)} u_{(2)} \otimes$

\hspace{2cm}
$(s_{(2)} \circ d). (b \circ u_{(3)}). \lambda^* _{(2)} \otimes
(s_{(1)} \circ \mu^* _{(3)}) \wedge (\lambda^* _{(1)} \circ v_{(3)})  \otimes$ 

\hspace{2cm}
$\mu^* _{(2)}.(t_{(1)} \circ c).(d \circ v_{(2)})
< \mu^* _{(1)}, \xi_{(1)}^{*T} >
< \lambda^* _{(3)}, \eta_{(2)}^{*T} >.$

\bigskip

{\noindent On the other hand,}

\bigskip

$\theta([s \otimes \lambda \otimes t \otimes a \otimes \xi \otimes b].
[u \otimes \mu \otimes v \otimes c \otimes \eta \otimes d])$

$= \sum \pm  t_{(3)} v_{(2)} \otimes 
[(t_{(1)} \circ \eta_{(3)}) \wedge (\xi_{(1)} \circ u_{(3)})]^*
\otimes s_{(2)}u_{(1)} \otimes$

\hspace{2cm}
$\eta_{(2)}.(s_{(1)} \circ d).(c \circ u_{(2)}) \otimes
[(\lambda_{(2)} \circ v_{(1)}) \vee (s_{(3)} \circ \mu_{(1)})]^* \otimes$ 

\hspace{2cm}
$(t_{(2)} \circ c). (a \circ v_{(3)}). \xi_{(2)}
< \eta_{(1)}, \lambda_{(1)}^T >
< \xi_{(3)}, \mu_{(2)}^T >.$

$= \sum \pm  t_{(3)} v_{(2)} \otimes 
(\xi_{(1)} ^* \circ u_{(3)}) \vee (t_{(1)} \circ \eta_{(3)} ^*)
\otimes s_{(2)}u_{(1)} \otimes$

\hspace{2cm}
$\eta_{(2)}.(s_{(1)} \circ d).(c \circ u_{(2)}) \otimes
(s_{(3)} \circ \mu_{(1)} ^*) \wedge (\lambda_{(2)} ^* \circ v_{(1)}) \otimes$ 

\hspace{2cm}
$(t_{(2)} \circ c). (a \circ v_{(3)}). \xi_{(2)}
< \eta_{(1)}, \lambda_{(1)}^T >
< \xi_{(3)}, \mu_{(2)}^T >.$

\bigskip

{\noindent Given the cocommutativity of the classical Schur algebra, it is 
clear from lemma~\ref{help}, that up to a sign, the terms in 
our expression for $$\theta([s \otimes \lambda \otimes t \otimes a \otimes \xi \otimes b].
[u \otimes \mu \otimes v \otimes c \otimes \eta \otimes d])$$ agree with the terms in
our expression for
$$\theta([s \otimes \lambda \otimes t \otimes a \otimes \xi \otimes b]).
\theta([u \otimes \mu \otimes v \otimes c \otimes \eta \otimes d]).$$
The difference in sign of each term is precisely
$$(-1)^{|\lambda_{(1)}||\lambda_{(2)}| + |\mu_{(1)}||\mu_{(2)}|
|\xi_{(1)}||\xi_{(2)}| + |\eta_{(1)}||\eta_{(2)}|}.$$

Diagrammatically comparing a term of
$$\theta(\theta([s \otimes \lambda \otimes t \otimes a \otimes \xi \otimes b]).
\theta([u \otimes \mu \otimes v \otimes c \otimes \eta \otimes d]))$$ with 
the relevant term of $$[s \otimes \lambda \otimes t \otimes a \otimes \xi \otimes b].
[u \otimes \mu \otimes v \otimes c \otimes \eta \otimes d],$$ one sees 
that their difference is also 
$$(-1)^{|\lambda_{(1)}||\lambda_{(2)}| + |\mu_{(1)}||\mu_{(2)}|
|\xi_{(1)}||\xi_{(2)}| + |\eta_{(1)}||\eta_{(2)}|}.$$

Therefore, our expressions for
$$\theta([s \otimes \lambda \otimes t \otimes a \otimes \xi \otimes b].
[u \otimes \mu \otimes v \otimes c \otimes \eta \otimes d]),$$ and
$$\theta([s \otimes \lambda \otimes t \otimes a \otimes \xi \otimes b]).
\theta([u \otimes \mu \otimes v \otimes c \otimes \eta \otimes d])$$ agree.
This completes the proof of theorem~\ref{autm}.
$\Box$

\begin{center}
{\large Schiver doubles and wreath products.}
\end{center}
 
Schiver bialgebras and their doubles may be understood to be  
generalisations of certain 
wreath products, as is illustrated by the following 
example: 

\bigskip

Let ${\cal S}(n)$ be the Schur bialgebra associated to $GL_n(k)$, and let 
the Schur algebra ${\cal S}(n,r)$ be the subalgebra of degree $r$.
Let the double of ${\cal S}(n)$ be denoted ${\cal D}(n)$, 
and let its degree $r$ part be 
written ${\cal D}(n,r)$.

Let $n \geq r$, and
let $\omega = (1^r)$ be the partition of $r$ with $r$ parts.
According to Green's presentation of Schur-Weyl duality
(theorem \ref{tensorspace2}), 
$$\xi_{\omega} {\cal S}(n,r) \xi_{\omega} \cong k\Sigma_r.$$
We have the following generalisation (for a further generalisation, see \cite{Turner2}, theorem 3):

\begin{proposition} \label{Weyl}

(a)
The endomorphism ring 
$\xi_{\omega} {\cal D}(n,r) \xi_{\omega}$ is isomorphic to
the wreath product
$k[x]/(x^2) \wr \Sigma_r$.

\bigskip

(b) If $char(k)=2$, then
$\xi_{\omega} {\cal D}(n,r) \xi_{\omega}$ is isomorphic to
the wreath product
$k\Sigma_2 \wr \Sigma_r$.

\bigskip

(c)
If $char(k)=0$, or $w< char(k)$, 
then ${\cal D}(n,r)$ is Morita equivalent to
$k[x]/(x^2) \wr \Sigma_r$. 
\end{proposition}
Proof:

Let $\Omega = \{ 1,...,r \}$.
Let us identify a subset $\Pi \subset \Omega$ of size $d$ with
the $d$-tuple $(\pi_1,...,\pi_d) \in I(n,d)$, where $\Pi$ is 
ordered just as $\Omega$ is ordered.

Thus, $\xi_{\Omega, \Omega} = \xi_\omega$.
The set 
$$\{ \xi_{\Pi \sigma, \Pi} \otimes 
X_{(\Omega-\Pi) \sigma, (\Omega-\Pi)} | 
\Pi \subset \Omega, \sigma \in \Sigma_r\}$$
is a basis for 
$\xi_{\omega} D(n,r) \xi_{\omega}$.  
The size of this basis set is $2^r.r!$.

The subspace spanned by
$\{ \xi_{\Omega \sigma, \Omega} | \sigma \in \Sigma_r \}$
is a subalgebra, naturally isomorphic to $k\Sigma_r$.

For $\Pi \subset \Omega$, let $\xi_\Pi = \xi_{\Pi,\Pi}$.
We show that the subspace spanned by 
$\{ \xi_{\Pi} \otimes 
X_{\Omega-\Pi} | 
\Pi \subset \Omega \}$ is also a subalgebra, isomorphic to
$(k[x]/(x^2))^{\otimes r}$ via
$$\Phi:
\xi_{\Pi} \otimes 
X_{\Omega-\Pi} \mapsto
z_1 \otimes...\otimes z_r,$$
where $z_i = 1$ if $i \in \Pi$, and $z_i = 0$ otherwise.

The isomorphism of the
lemma is then quite plain from the formulas,
$$\xi_{\Pi} \otimes X_{\Omega-\Pi}. 
\xi_{\Omega \sigma, \Omega} =
\xi_{\Pi \sigma,\Pi} \otimes X_{\Omega-\Pi \sigma, \Omega-\Pi},$$
$$\xi_{\Omega \sigma^{-1}, \Omega}.
\xi_{\Pi} \otimes X_{\Omega-\Pi}. 
\xi_{\Omega \sigma, \Omega} =
\xi_{\Pi} \otimes
X_{(\Omega-\Pi) \sigma}.$$ 
Let $\Gamma, \Pi$ be subsets of $\Omega$.
Let $s = \xi_\Pi$, and let $a=X_{\Omega-\Pi}$.
Let $t = \xi_\Gamma$, and let $b = X_{\Omega-\Gamma}$.
We compute the product of $s \otimes a$ and $t \otimes b$.
$$\Delta(\xi_\Pi) = 
\sum_{\alpha \subset \Pi} \xi_\alpha \otimes \xi_{\Pi-\alpha}$$
$$\Delta(\xi_\Gamma) = 
\sum_{\beta \subset \Pi} \xi_\beta \otimes \xi_{\Pi-\beta}.$$
Thus, non-zero terms of $(s \otimes a)(t \otimes b)$ only appear when
$b_{(1)} = X_\alpha$,
$s_{(2)} = \xi_{\Pi - \alpha}$,
$a_{(2)} = X_{\Gamma - \beta}$,
and $t_{(1)} = \xi_\beta$.

This implies, in turn, that only $b = b_{(1)} = b_{(2)} = X_\alpha$ gives
possible non-zero terms in $(s \otimes a).(t \otimes b)$.
Thus, $\alpha = \Omega-\Gamma$.

In addition, only 
$a = a_{(1)} = a_{(2)} = X_{\Gamma - \beta}$
gives possible non-zero terms in
$(s \otimes a).(t \otimes b)$.
Thus, $\Gamma-\beta = \Omega - \Pi$.

For the product $s_{(2)} t_{(1)}$ to be non-zero, we now require
$\beta = \Pi - \alpha$.
Thus, $\beta = \Gamma \cap \Pi$ and $\alpha = \Omega - \Gamma$.

From all this, we conclude that 
$(s \otimes a).(t \otimes b)$ is zero unless 
$\Pi \cup \Gamma = \Omega$, and that in this case
$(s \otimes a).(t \otimes b)$ is equal to
$\xi_{\Pi \cap \Gamma} \otimes X_{\Omega - \Pi \cap \Gamma}$. 
It follows that 
$\Phi$ is indeed an algebra isomorphism, and (a) is proven !

\bigskip
(b) is immediate from (a), since $k\Sigma_2 \cong k[x]/(x^2)$, 
in characteristic two.

(c) also
follows immediately from (a), because $k\Sigma_r$ is semisimple so long as    
$r!$ is invertible in $k$. Thus, ${\cal D}(n,r)$ and $k[x]/(x^2) \wr \Sigma_r$ 
have the same number of simple modules. 
$\Box$

\begin{center}
{\large Schiver doubles and blocks of Hecke algebras.}
\end{center}

Here is a summary of our education concerning Rock blocks of symmetric
groups, in characteristic two, in chapters 4 and 5.

\begin{theorem}
Let $k$ be a field of characteristic two.
Let $k{\bf B}_{\rho,w}^\Sigma$ 
be a Rock block of a symmetric group, whose weight is $w$.
Then $k{\bf B}_{\rho,w}^\Sigma$ contains a nilpotent ideal ${\cal N}$, such that
$k{\bf B}_{\rho,w}^\Sigma/ {\cal N}$ is Morita equivalent to ${\cal S}(w,w)$.  
Let $e$ be an idempotent in $k{\bf B}_{\rho,w} ^\Sigma$, 
such that the indecomposable
summands of $k{\bf B}_{\rho,w}^\Sigma e$ are precisely those indecomposable
summands with tops in the set,
$$\{ D^{[\emptyset, \lambda]} | \lambda \textrm{ is }p
\textrm{-regular }\}.$$
Then $e k{\bf B}_{\rho,w}^\Sigma e$ is Morita equivalent to 
$k\Sigma_2 \wr \Sigma_w$.  
$\Box$
\end{theorem}

This theorem, proposition~\ref{Weyl},
and indeed example~\ref{Karin} of this paper, 
give evidence that Rock blocks $k{\bf B}_{\rho,w}^\Sigma$, are
in fact Morita equivalent
to the Schur doubles ${\cal D}(n,w)$, for $n \geq w$, when the field $k$ has 
characteristic two.
But the Ringel dualities of chapter 6 impose a more general conjecture on us:

\bigskip

Let $A_{p-1}$ be the ordinary Dynkin graph with $p-1$ vertices:

\bigskip

\hspace{1cm}$\circ$--------------$\circ$--------------$\circ$.............$\circ$--------------$\circ$--------------$\circ$

\bigskip

Let $q \in k^\times$, and let $p$ be the least natural 
number such that, 
$$1+q+...+q^{p-1} = 0.$$ 
Let $n,w$ be natural numbers, such that $n \geq w$.
Let ${\cal D}_{A_{p-1}}(n)$ be 
the Schiver double associated to the graph
$A_{p-1}$.

\begin{conjecture} \label{Rock}
The degree $w$ part ${\cal D}_{A_{p-1}}(n,w)$ of ${\cal D}_{A_{p-1}}(n)$,
is Morita equivalent to the Rock block $k{\bf B}_{\rho,w}^{{\cal H}_q}$
of a Hecke algebra, whose weight is $w$.

Indeed,
${\cal D}_{A_{p-1}}(n,w)$ is derived equivalent to any block 
$k{\bf B}_{\tau,w}^{{\cal H}_q}$ of 
a Hecke algebra, whose weight is $w$.
\end{conjecture}

Since any two blocks of the same weight are derived equivalent,
the conjectured derived equivalences would follow immediately 
from the conjectured Morita equivalences. 

\begin{remark}
When $q=1$, the Hecke algebra ${\cal H}_q(\Sigma_n)$ is 
isomorphic to the 
symmetric group algebra $k\Sigma_n$, and $p$ is merely
the characteristic of the field $k$. In this case,
conjecture~\ref{Rock} could be viewed 
as relating to a non-abelian generalisation of 
Brou\'e's abelian defect group
conjecture \cite{Broue} for symmetric groups, as proved by Chuang, 
Kessar, and Rouquier \cite{KessCh}, \cite{ChRou}.

So far, however, I can see 
no natural interpretation of the Schiver double ${\cal D}_{A_{p-1}}(n)$
in terms of $p$-local group theoretic information, 
weights for symmetric groups, etc. etc.
\end{remark}

\begin{center}
{\large Doubles, and graded rings associated to Rock blocks of
Hecke algebras.}
\end{center}

I am unable to prove conjecture~\ref{Rock}. 
The best I can do, is describe a filtration, 
$$k{\bf B}_{\rho,w} ^{{\cal H}_q}= {\cal N}[0] \supset {\cal N}[1] \supset 
{\cal N}[2] \supset ...
\supset {\cal N}[2w] \supset {\cal N}[2w+1] =0,$$
\index{${\cal N}[i]$}
on the Rock block of a Hecke algebra of weight $w$, such that 
${\cal N}[i].{\cal N}[j] \subseteq {\cal N}[i+j]$, and
the associated graded ring,
$$gr_{\rho,w}^{{\cal H}_q} = 
\bigoplus_{i=0} ^{2w} {\cal N}[i]/{\cal N}[i+1],$$
\index{$gr_{\rho,w}^{{\cal H}_q}$}
resembles the Schiver double 
${\cal D}_{A_{p-1}}(w,w)$.
Suppose $p \geq 3$.
We can show that the degree zero part of the graded ring $gr_{\rho,w}^{{\cal H}_q}$ is Morita equivalent
to the degree zero part ${\cal S}_{V(A_{p-1})}(w,w)$ of ${\cal D}_{A_{p-1}}(w,w)$, and that 
the homogeneous components of $gr_{\rho,w}^{{\cal H}_q}$ correspond via Morita equivalence, 
in the category of ${\cal S}_{V(A_{p-1})}(w,w)$-${\cal S}_{V(A_{p-1})}(w,w)$-bimodules, to the homogeneous 
components of ${\cal D}_{A_{p-1}}(w,w)$.

In this section, we give a sketch of a proof of this fact.  
The most obvious obstacle to proving conjecture~\ref{Rock}, 
concerning Rock blocks of Hecke algebras,
is our inability to show that $k{\bf B}_{\rho,w} ^{{\cal H}_q}$ is graded, in a certain way.

\bigskip

\emph{Step 1.} 
If $K$ is a certain 
field of characteristic zero, 
and $q \in K$ is a primitive $p^{th}$ root of unity,
then $K{\bf B}_{\rho,w}^{{\cal H}_q}$ is Morita equivalent to
$K{\bf B}_{\emptyset,1}^{{\cal H}_q} \wr \Sigma_w$.

Proof:

By Dirichlet's theorem, there exists a prime number $l$, such that 
$l=1$ (modulo $p$). Similarly, there exists a 
prime number $\bar{q}$, such that $\bar{q} = q$ (modulo $l$).
By theorem~\ref{mainnew}, there is an $l$-modular system $(K, {\cal O},k)$,
such that $k{\bf B}_{\rho,w}^{{\cal H}_{\bar{q}}}$ is Morita equivalent to
$k{\bf B}_{\emptyset,1}^{{\cal H}_{\bar{q}}} \wr \Sigma_w$.
We can lift this equivalence by the following argument,
due to Joe Chuang.

The bimodule inducing Morita equivalence is a
summand of the $k{\bf B}_{\rho,w}^{{\cal H}_{\bar{q}}}$- 
$k{\bf B}_{\emptyset,1}^{{\cal H}_{\bar{q}}} \wr \Sigma_w$-bimodule,
$$kT = k{\bf B}_{\rho,w}^{{\cal H}_{\bar{q}}} 
\bigotimes_{\otimes^w k{\bf B}_{\emptyset,1}^{{\cal H}_{\bar{q}}}}
k{\bf B}_{\emptyset,1}^{{\cal H}_{\bar{q}}} \wr \Sigma_w.$$
\index{$kT$}
We would like to define a summand of the ${\cal O}{\bf B}_{\rho,w}^{{\cal H}_{q}}$- 
${\cal O}{\bf B}_{\emptyset,1}^{{\cal H}_{q}} \wr \Sigma_w$-bimodule,
$$T = {\cal O}{\bf B}_{\rho,w}^{{\cal H}_{q}} 
\bigotimes_{\otimes^w K{\bf B}_{\emptyset,1}^{{\cal H}_{q}}}
{\cal O}{\bf B}_{\emptyset,1}^{{\cal H}_{q}} \wr \Sigma_w,$$
which induces a Morita equivalence. It is enough lift idempotents from
$End(kT)$ to $End(T)$. 
Algebraically, this translates to the problem of
lifting centralizers of parabolic subalgebras from characteristic $l$ to 
characteristic zero. The arguments of A. Francis 
(\cite{Francis}, 3.6, 3.8)
show that this is possible.

\bigskip

\emph{Step 2.} 
Use the approach of Cline, Parshall, and Scott  (\cite{CPS4}, 5.3) 
to generalise theorem~\ref{Comeon} from blocks
of symmetric groups to blocks of Hecke algebras.
Generalise the results of chapters 6 and 7 from blocks of symmetric groups,
to blocks of Hecke algebras. 
Thus produce a nilpotent ideal ${\cal N}$ of $k{\bf B}_{\rho,w} ^{{\cal H}_q}$,
such that $k{\bf B}_{\rho,w} ^{{\cal H}_q}/{\cal N}$ is Morita equivalent to
${\cal S}_{V(A_{p-1})}(w,w)$. 

\bigskip

\emph{Step 3.}
Define the filtration ${\cal N}[i]$ on $k{\bf B}_{\rho,w}^{{\cal H}_q}$, using
good idempotents for the $q$-Schur algebra, as well as the signature
automorphism on the Hecke algebra. 

The method of definition generalises that of the 
bimodule $N$, of chapter 6. Whilst $N$ is defined to be 
$x k{\bf B}_{\rho,w}^{{\cal H}_q} y$, for fixed idempotents $x,y$, 
the ideal ${\cal N}[i]$ contains
sums of terms $x {\cal H}_q y$, $y {\cal H}_q x$, 
$x {\cal H}_q y {\cal H}_q x$, and $y {\cal H}_q x {\cal H}_q y$,
for various idempotents $x,y$.

By comparison with the characteristic zero case (Step 1), 
it can be seen that the filtration satisfies 
${\cal N}[i].{\cal N}[j] \subseteq {\cal N}[i+j]$, 
and ${\cal N} = {\cal N}[1]$.

\bigskip

\emph{Step 4.}
Show that the $k{\bf B}_{\rho,w} ^{{\cal H}_q}/{\cal N}$-
$k{\bf B}_{\rho,w} ^{{\cal H}_q}/{\cal N}$-bimodule,
${\cal N}[i]/{\cal N}[i+1]$ corresponds, via Morita equivalence, to 
the ${\cal D}^0_{A_{p-1}}(w,w)$-${\cal D}^0_{A_{p-1}}(w,w)$-bimodule, 
${\cal D}_{A_{p-1}}^i(w,w)$.
Prove this by induction on $w$.

\bigskip

Steps 1-4 imply that the graded components of the algebras 
$gr_{\rho,w}^{{\cal H}_q}$ and ${\cal D}_{A_{p-1}}(w,w)$ are in natural correspondence.
A more ambitious project would be to follow through steps 5-8, and thus prove that
$gr_{\rho,w}^{{\cal H}_q}$ and ${\cal D}_{A_{p-1}}(w,w)$ are Morita equivalent.
There are possible difficulties in pushing this through.
In an analogous, but more elementary situation, we have succeeded in overcoming the necessary obstacles, 
in work with V. Miemietz \cite{GL2}.

\emph{Step 5.}
Show that for compositions $(a_1,...,a_{p-1}), (c_1,..,c_{p-1})$ of $w$,
for which there exists a composition $(b_1,...,b_{p-1})$ of $w$, such that,
$$(a_1,...,a_{p-1}) \lhd (b_1,...,b_{p-1}) \lhd (c_1,...,c_{p-1}),$$
we have,
$$Ext^1({\cal L}_q(a_1,...,a_{p-1}), {\cal L}_q(c_1,...,c_{p-1}))=$$
$$Ext^1({\cal L}_q(c_1,...,c_{p-1}), {\cal L}_q(a_1,...,a_{p-1}))= 0.$$
Prove this by induction on $w$, using duality and Frobenius reciprocity,
with base case $w=1$.

\bigskip

\emph{Step 6.}
Show (using Step 5) that the multiplication morphism, 
$$\phi: {\cal N}[i]/{\cal N}[i+1] \bigotimes_{k{\bf B}_{\rho,w}^{{\cal H}_q}} 
{\cal N}[j]/{\cal N}[j+1] \rightarrow {\cal N}[i+j]/{\cal N}[i+j+1],$$
of $k{\bf B}_{\rho,w} ^{{\cal H}_q}/{\cal N}$-
$k{\bf B}_{\rho,w} ^{{\cal H}_q}/{\cal N}$-bimodules, is a surjection.
In other words, $gr_{\rho,w} ^{{\cal H}_q}$ is generated in degrees
$0$ and $1$.
 
\bigskip

\emph{Step 7}
The generalized Koszulity of ${\cal D}_{A_\infty}(w,w)$ (remark~\ref{Koszul})
implies that ${\cal D}_{A_\infty}(w,w)$ is quadratic.
When $p \geq 4$, ${\cal D}_{A_{p-1}}(w,w)$ is also quadratic.

The correspondences of Step 4 may be fixed so that in degree $(1,1)$, 
$\phi$ corresponds to the morphism 
$${\cal D}_{A_{p-1}}^1(w,w) \bigotimes_{{\cal D}^0_{A_{p-1}}(w,w)} 
{\cal D}_{A_{p-1}} ^1(w,w)$$ 
$$\rightarrow 
{\cal D}_{A_{p-1}} ^2(w,w).$$
of ${\cal D}^0_{A_{p-1}}(w,w)$-${\cal D}^0_{A_{p-1}}(w,w)$-bimodules.

\bigskip

\emph{Step 8.} 
Conclude from Steps 6 and 7 that ${\cal D}_{A_{p-1}}(w,w)$ surjects
onto an algebra Morita equivalent to $gr_{\rho,w} ^{{\cal H}_q}$.
By a dimension count, this surjection is an isomorphism. $\Box$

\newpage
\begin{center}
{\bf \large 
Chapter IX

Power sums.}
\end{center}

In this chapter, 
we define certain chain complexes for the Schiver doubles, 
whose Grothendieck character
describes the symmetric functions $p_r = x_1^r + x_2^r + ...$ 
(theorem~\ref{pop}).

Of course, it is not difficult to define such complexes in a naive way: take standard modules for the Schur algebra indexed 
by hook partitions,  place them in homological degree $i$, where $i$ is the number of parts of the partition,
and then give them a zero differential.

However, we describe here a more subtle method, which uses the structure of the doubles, rather than merely the Schur algebra.
The reason for expecting such complexes to exist, and defining them, is the categorification program.
Indeed, the existence of such complexes, which invoke the structure of the Schiver doubles, is consistent 
with the apparent affinity between derived categories of Schiver doubles, and those of blocks of Hecke algebras. 

Whilst blocks of Hecke algebras define a category lifting the Fock space realization of the basic representation
of $\widehat{\mathfrak{sl}_p}$, power sums play a defining role in the combinatorial formulation of the principal homogeneous 
realization of the basic representation for $\widehat{\mathfrak{sl}_p}$
by I. Frenkel, N. Jing and W. Wang.
We expect the equivalences between blocks of Hecke algebras and doubles to be one aspect of a 
categorical realisation of the isomorphism between the Fock space realization and the principal homogeneous realization.
The description of induction and restriction functors between symmetric groups via certain functors between doubles should be 
another aspect. Indeed, we expect functors between doubles which correspond to power sums.
Note that such functors should be realised between doubles, and not their Schur algebra quotients,
because it is the doubles which we expect to categorify the principal homogeneous realization, and not their quotients.

\begin{center}
{\large Complexes for Schiver doubles, and power sums.}
\end{center}

Most of the bases for the ring of symmetric functions given in
I. MacDonald's book \cite{Mac}
have natural interpretations as characters of modules
for the Schur algebra.
Elementary symmetric functions correspond
to exterior powers of the natural module for $M_n(k)$, complete symmetric
functions correspond to symmetric powers of the natural module,
and Schur functions correspond to Weyl modules.
However, the power sums $p_r = x_1^r + x_2^r + ...$ 
have no such interpretation.

In this section, we describe complexes ${\cal P}_r$
for Schiver doubles ${\cal D}_\Gamma(n,r)$, whose
homology describes the power sum $p_r$.
Indeed, we may define one such complex ${\cal P}_r (a)$, 
for every pair $(\gamma, a)$, where $\gamma$ is an vertex of $\Gamma$, and 
$a$ an edge emanating from $\gamma$. 

\bigskip

\begin{lemma} \label{quim}
Let $n \geq r$.
Let $Q$ be a quiver. Suppose that $\gamma_1,\gamma_2$ are vertices
in $Q$.
Suppose that $a$ is an arrow
in $Q$, whose source is $\gamma_1$, and tail is $\gamma_2$,
and that $a$ is the only such arrow.
Let $r_1, r_2$ be natural numbers, whose sum is $r$.
There is a graded module $M_{r_1,r_2} = M_{r_1, r_2}(a)$ 
for ${\cal S}_Q(n,r)$, whose  
graded pieces are the ${\cal S}_{V(Q)}(n,r)$-modules,
$$\Delta_{\gamma_1}(r_1) \otimes E_{\gamma_2} ^{\otimes r_2}$$
$$\Delta_{\gamma_1}(r_1,1) \otimes E_{\gamma_2} ^{\otimes r_2-1}$$
$$...$$
$$\Delta_{\gamma_1}(r_1, 1^{r_2}),$$
in degrees $0,1,...,r_2$.
\end{lemma}
Proof:

Let $\xi_{r_1,r_2}$ be the unit of the algebra 
${\cal S}_{\gamma_1}(n,r_1) \otimes {\cal S}_{\gamma_2}(n,r_2)$, 
an idempotent in ${\cal S}_Q(n,r)$.
Then, 
$$N_{r_1, r_2}(a) = {\cal S}_Q(n,r) \xi_{r_1,r_2} 
\bigotimes_{{\cal S}_{\gamma_1}(n,r_1) \otimes {\cal S}_{\gamma_2}(n,r_2)}
\Delta_{\gamma_1}(r_1) \otimes E_{\gamma_2} ^{\otimes r_2},$$
is a graded module for ${\cal S}_Q(n,r)$, whose  
graded pieces are,
$$\Delta_{\gamma_1}(r_1) \otimes E_{\gamma_2} ^{\otimes r_2}$$
$$\Delta_{\gamma_1}(r_1) \otimes E_{\gamma_1} 
\otimes E_{\gamma_2} ^{\otimes r_2-1}$$
$$...$$
$$\Delta_{\gamma_1}(r_1) \otimes E_{\gamma_1} ^{\otimes r_2},$$
in degrees $0,1,...,r_2$.

For a partition $\lambda$ of $t$, let 
$\xi_{>\lambda} = \sum_{\mu>\lambda} \xi_\mu$ be the sum of
Green's idempotents $\xi_\mu$ \cite{Green}, 
corresponding to partitions $\mu$,
greater than $\lambda$, with respect to the dominance ordering. 

Let $i_{r_1,r_2}^j$ 
be the idempotent $\xi_{>(r_1,1^j)} \otimes \xi_{(1^{r_2-j})}$,
an element of the algebra,
$${\cal S}_{\gamma_1}(n,r_1+j) \otimes
{\cal S}_{\gamma_2}(n,r_2-j).$$
Let $i_{r_1,r_2} = \sum_{j=1}^{r_2} i_{r_1,r_2}^j$, 
an idempotent in the algebra, 
$$\bigoplus_{j=1}^{r_2} {\cal S}_{\gamma_1}(n,r_1+j) \otimes
{\cal S}_{\gamma_2}(n,r_2-j).$$
Let us define,
$$M_{r_1, r_2}(a)= N_{r_1, r_2}(a) / 
{\cal S}_Q(n,r)i_{r_1,r_2}N_{r_1, r_2}(a),$$ 
to be the quotient of $N_{r_1, r_2}(a)$, relative to a trace 
from the projective module ${\cal S}_Q(n,r)i_{r_1,r_2}$.

Since the tensor product of a projective 
${\cal S}_{\gamma_1}(n,r_1+j)$-module with
the natural representation $E_{\gamma_1}$ is projective, 
with $\Delta$-composition factors given by the branching rule, 
we may compute the composition factors of $M_{r_1,r_2}(a)$, and find
them to be,
$$\Delta_{\gamma_1}(r_1) \otimes E_{\gamma_2} ^{\otimes r_2}$$
$$\Delta_{\gamma_1}(r_1,1) \otimes E_{\gamma_2} ^{\otimes r_2-1}$$
$$...$$
$$\Delta_{\gamma_1}(r_1, 1^{r_2}),$$
in degrees $0,1,...,r_2$.
$\Box$

The following theorem generalises the above lemma:

\begin{theorem} \label{pop}
Let $n \geq r$.
Let $Q$ be a quiver. Suppose that $\gamma_1,\gamma_2$ are vertices
in $Q$.
Suppose that $a$ is an arrow
in $Q$, whose source is $\gamma_1$, and tail is $\gamma_2$,
and that $a$ is the only such arrow.
Let $r_1, r_2$ be natural numbers, whose sum is $r$.
There is a graded module $C_{r_1,r_2} = C_{r_1, r_2}(a)$ 
for ${\cal D}_Q(n,r)$, whose  
graded pieces are the ${\cal S}_{V(Q)}(n,r)$-modules,
$$\Delta_{\gamma_1}(r_1) \otimes E_{\gamma_2} ^{\otimes r_2}$$
$$\Delta_{\gamma_1}(r_1-1) \otimes E_{\gamma_2} ^{\otimes r_2+1} +
\Delta_{\gamma_1}(r_1, 1) \otimes E_{\gamma_2} ^{\otimes r_2 -1}$$
$$\Delta_{\gamma_1}(r_1-1,1) \otimes E_{\gamma_2} ^{\otimes r_2} +
\Delta_{\gamma_1}(r_1, 1^2) \otimes E_{\gamma_2} ^{\otimes r_2 -2}$$
$$...$$
$$\Delta_{\gamma_1}(r_1-1,1^{r_2 - 1}) \otimes E_{\gamma_2} ^{\otimes 2} +
\Delta_{\gamma_1}(r_1, 1^{r_2 })$$
$$\Delta_{\gamma_1}(r_1-1, 1^{r_2}) \otimes E_{\gamma_2},$$
in degrees $0,1,...,r_2+1$. 

There is a natural homomorphism $d_{r_1,r_2} =d_{r_1,r_2}(a)$ 
of degree one, from
$C_{r_1-1,r_2+1}$ to $C_{r_1, r_2}$. The sequence of maps
$\{ d_{r_1,r_2} | r_1 + r_2 =r \}$ defines a chain complex, 
${\cal P}_{r}(a)$,
given by,
$$... \rightarrow 0 \rightarrow C_{1,r-1} 
\rightarrow C_{2, r-2} \rightarrow ...
\rightarrow C_{r-1,1} \rightarrow C_{r,0} \rightarrow 0 
\rightarrow ...$$
The homology of this complex at term 
$C_{r_1,r_2}$ is isomorphic as a ${\cal D}_Q(n,r)$-module,
to the standard module
$\Delta_{\gamma_1}(r_1,1^{r_2})$ for ${\cal S}_{\gamma_1}(n,r)$, 
concentrated in degree $r_1-1$. 
\end{theorem}
Proof:

Consider the ${\cal D}_Q(n,r)-
{\cal S}_{\gamma_1}(n,r_1) \otimes {\cal S}_{\gamma_2}(n,r_2)$-bimodule,
$$X_{r_1,r_2} =
{\cal S}_Q(n,r)\xi_{r_1,r_2} \oplus 
({\cal S}_Q(n,r-1)\xi_{r_1-1,r_2} 
\otimes \bigwedge(n,1)_{\gamma_2,\gamma_1}),$$
given as a quotient of the ${{\cal D}_Q(n,r)}$-
${{\cal S}_{\gamma_1}(n,r_1) \otimes {\cal S}_{\gamma_2}(n,r_2)}$-bimodule,
$${\cal D}_Q(n,r) \xi_{r_1,r_2},$$
modulo terms of higher degree.
Note that
$$Y_{r_1,r_2} =
({\cal S}_Q(n,r-1)\xi_{r_1-1,r_2} 
\otimes \bigwedge_{\gamma_2,\gamma_1}(n,1)),$$
is a sub-bimodule of $X_{r_1,r_2}$.
Let us define the ${\cal D}_Q(n,r)$-module,
$$U_{r_1,r_2} = X_{r_1,r_2} \bigotimes_
{{\cal S}_{\gamma_1}(n,r_1) \otimes {\cal S}_{\gamma_2}(n,r_2)}
(\Delta_{\gamma_1}(r_1) \otimes E_{\gamma_2}^{\otimes r_2}),$$
which contains the submodule,
$$V_{r_1,r_2} = Y_{r_1,r_2} \bigotimes_
{{\cal S}_{\gamma_1}(n,r_1) \otimes {\cal S}_{\gamma_2}(n,r_2)}
(\Delta_{\gamma_1}(r_1) \otimes E_{\gamma_2}^{\otimes r_2}).$$
Thus, $U_{r_1,r_2}$ is a graded module, whose graded pieces
are the ${\cal S}_{V(Q)}(n,r)$-modules,
$$\Delta_{\gamma_1}(r_1) \otimes E_{\gamma_2} ^{\otimes r_2}$$
$$\Delta_{\gamma_1}(r_1-1) \otimes E_{\gamma_2} ^{\otimes r_2+1} +
\Delta_{\gamma_1}(r_1) \otimes E_{\gamma_1} 
\otimes E_{\gamma_2} ^{\otimes r_2 -1}$$
$$\Delta_{\gamma_1}(r_1-1) \otimes E_{\gamma_1}  
\otimes E_{\gamma_2} ^{\otimes r_2} +
\Delta_{\gamma_1}(r_1) \otimes E_{\gamma_1}^{\otimes 2} 
\otimes E_{\gamma_2} ^{\otimes r_2 -2}$$
$$...$$
$$\Delta_{\gamma_1}(r_1-1) \otimes E_{\gamma_1} ^{r_2-1}  
\otimes E_{\gamma_2} ^{\otimes 2} +
\Delta_{\gamma_1}(r_1, 1^{r_2 })$$
$$\Delta_{\gamma_1}(r_1-1) \otimes E_{\gamma_1}^{r_2}
\otimes E_{\gamma_2},$$
in degrees $0,1,...,r_2+1$. 
Let us define,
$$O_{r_1,r_2} = {\cal D}_Q(n,r) i_{r_1-1,r_2+1} V_{r_1,r_2},$$
a submodule of $U_{r_1,r_2}$. 
As in lemma~\ref{quim}, 
$D_{r_1,r_2} = V_{r_1,r_2}/O_{r_1,r_2}$ 
is a graded module, whose graded pieces
are the ${\cal S}_{V(Q)}(n,r)$-modules,
$$\Delta_{\gamma_1}(r_1-1) \otimes E_{\gamma_2} ^{\otimes r_2+1}$$
$$\Delta_{\gamma_1}(r_1-1,1)  
\otimes E_{\gamma_2} ^{\otimes r_2}$$
$$...$$
$$\Delta_{\gamma_1}(r_1-1,1^{r_2-2}) \otimes   
\otimes E_{\gamma_2} ^{\otimes 2}$$
$$\Delta_{\gamma_1}(r_1-1, 1^{r_2})
\otimes E_{\gamma_2},$$
in degrees $0,1,...,r_2$. 
Let us define,
$$P_{r_1,r_2} = {\cal D}_Q(n,r)i_{r_1,r_2}(U_{r_1,r_2}/O_{r_1,r_2}),$$
a submodule of $M_{r_1,r_2}/O_{r_1,r_2}$.
We now define, 
$$C_{r_1,r_2} = U_{r_1,r_2}/(O_{r_1,r_2}+P_{r_1,r_2}),$$
which contains $D_{r_1,r_2}$ as a submodule. 
The graded pieces of $C_{r_1,r_2}$ are visibly those described in the 
statement of the theorem.

Let $\eta$ be the element, 
$1_{{\cal S}_{\gamma_1}(n,r_1-1)} \otimes 1_{M_n(k)} \otimes 
1_{{\cal S}_{\gamma_2}(n,r_2+1)}$ of the
${\cal S}_{\gamma_1}(n,r_1-1) \otimes {\cal S}_{\gamma_2}(n,r_2+1)$-
${\cal S}_{\gamma_1}(n,r_1) \otimes {\cal S}_{\gamma_2}(n,r_2)$- bimodule,
$${\cal S}_{\gamma_1}(n,r_1-1) \otimes 
\bigwedge_{\gamma_2, \gamma_1}(n,1) \otimes
{\cal S}_{\gamma_2}(n,r_2).$$ 
Then multiplication by $\eta$ defines a map from
$X_{r_1-1,r_2+1}$ to $Y_{r_1,r_2}$.
By restriction, there is a map, 
$$d_{r_1,r_2}: C_{r_1-1,r_2+1} \rightarrow C_{r_1,r_2}.$$
The kernel of $d_{r_1+1,r_2-1}$ is equal to the submodule 
$D_{r_1-1,r_2+1} + \Delta_{\gamma_1}(r_1,1^{r_2})$ of
$C_{r_1,r_2}$, and the image of $d_{r_1,r_2}$
is isomorphic to $D_{r_1,r_2}$.
The chain complex ${\cal P}_r(a)$ defined by the sequence of maps, thus
has homology $\Delta_{\gamma_1}(r_1, 1^{r_2})$, at term $C_{r_1,r_2}$.
$\Box$

\begin{theorem}
The Grothendieck character of the complex ${\cal P}_{r}(a)$ 
of ${\cal D}_Q(n,r)$-modules describes the power
sum $p_r$, at $\gamma_1$.
\end{theorem}
Proof:

The Grothendieck character of a standard module $\Delta(\lambda)$
is given by the Schur function $s_\lambda$.
There is a formula for the power sum $p_r$, given by 
(\cite{Mac}, I.4, example 10),
$$p_r = s_{(r)} - s_{(r-1,1)} + s_{(r-2,1^2)} -... \pm s_{(1^r)}.$$
The Grothendieck character of ${\cal P}_r(a)$ thus describes the 
power sum $p_r$. $\Box$

\bigskip

I. Frenkel, N. Jing and W. Wang have given a description of
the homogeneous vertex operator construction of the basic representation
of an affine Lie algebra of type $ADE$, via wreath 
products of finite group algebras \cite{Fren}.
The characters of symmetric groups which are used in that paper correspond
in symmetric function theory to elementary symmetric functions,
to complete symmetric functions, and to power sums.
 
We have observed here that there are objects in the derived category 
of a Schiver double which correspond to all these functions.

In fact, upon studying Frenkel, Jing and Wang's construction more carefully, 
one realises that the   
Schiver double afforded to a simply-laced Dynkin diagram $\Gamma$, 
underlies a category for a vertex representation for 
the affinization of the Kac-Moody Lie algebra defined by $\Gamma$, 
at least when $\Gamma$ is ordinary/affine, of type $ADE$.

By this, we only mean that we can describe a category whose complexified 
Grothendieck group is the vertex representation, and we can describe functors 
which correspond to the vertex operators, upon passing to the
Grothendieck group.
The category is a direct product of derived categories for Schiver doubles,
and the functors are described by certain complexes of bimodules.
Note that we have \emph{not} explored the 
extent to which relations in the affine Lie algebra lift to relations
between functors. 

Conjecture~\ref{Rock} is comfortingly consistent with this categorical
perspective.
The blocks of Hecke algebras are already well known to be
categories which describe the basic representation for
$\widehat{\mathfrak{sl}_p}$,
as has most elegantly been described by I. Grojnowski \cite{Groj},
following theory of A. Lascoux, B. Leclerc, and J-Y. Thibon \cite{LLT}, 
as well as S. Ariki \cite{Ariki}, and A. Kleshchev \cite{Klesh1}, \cite{Klesh2}.
In Grojnowski's article, relations in the affine Lie algebra do
lift to relations between functors. 

In general, we expect an equivalence of 
$\widehat{\mathfrak{sl}_p}$-categorifications,
$$D^b \left( \bigoplus_{w \geq 0, s \in \widehat{W}/ W} {\cal D}_{A_{p-1}}(w,w) \right) \cong
D^b \left( \bigoplus_{r \geq 0} {\cal H}_q(\Sigma_r) \right),$$
although it is unclear to me what this means precisely; 
so far, only $\mathfrak{sl}_2$-categorifications possess an axiomatic definition and a general theory
\cite{ChRou}. Above, $W$ is the Weyl group of $\mathfrak{sl}_p$, and
$\widehat{W}$ is the corresponding affine Weyl group.

\newpage

\begin{center}
{\bf \large 
Chapter X

Schiver doubles of type $A_\infty$.}
\end{center}

We consider the infinite Dynkin quiver ${A}_{\infty}$: \index{$A_\infty$}

\bigskip

\hspace{1.5cm}..........$\circ$--------------$\circ$--------------$\circ$--------------$\circ$..........

\bigskip

We prove that the module category of ${\cal D}_{A_\infty}(w,w)$ is
a highest weight category, and then  
conjecture that Rock blocks $k{\bf B}_{\rho,w}^{{\cal S}_q}$, 
of $q$-Schur algebras, are Morita equivalent to certain subquotients of
the algebra ${\cal D}_{A_\infty}(w,w)$, defined in chapter 8
(conjecture~\ref{SchurRock}).

We also describe a walk along ${\cal D}_{A_\infty}(n)$, analogous to 
J.A. Green's
walk along the Brauer tree (theorem~\ref{walk}).

\begin{center}
{\large Schiver bialgebras for $A_\infty$.}
\end{center}

Consider the infinite Dynkin quiver ${A}_{\infty}$, \index{$A_\infty$}
which has
vertex set 
$$V_\infty = \{ v_i, i \in \mathbb{Z} \},$$
and edge set $E_\infty = \{ e_i, i \in \mathbb{Z} \}$.
The source of any edge $e_i$ is the vertex
$v_{i+1}$, and its tail is $v_{i}$.

We present the Schiver bialgebra corresponding to $A_\infty$. 
This bialgebra amplifies the category of chain complexes over $k$, 
as the classical 
Schur bialgebra amplifies the category of vector spaces over $k$.
 
\begin{definition}
Let $\underline{\Lambda}(n,r)$ 
(respectively $\underline{\Lambda}'(n,r)$)
\index{$\underline{\Lambda}(n,r)$}
be the set 
$\{ \underline{\lambda} = (\lambda_i)_ {i \in \mathbb{Z}} \}$ of sequences 
of partitions, with $n$ parts or fewer,
whose sizes sum to $r$.

For two elements $\underline{\lambda}, \underline{\mu}$ of 
$\underline{\Lambda}(n,r)$
(respectively $\underline{\Lambda}'(n,r)$),
let $\underline{\lambda} \unlhd \underline{\mu}$ if and only if,
the sequence $(\lambda_i)_ {i \in \mathbb{Z}}$ can be obtained from 
the sequence
$(\mu_i)_ {i \in \mathbb{Z}}$ in finitely many steps, by repeatedly
either,

(1) removing a box from the Young diagram of $\mu_i$, and replacing it 
lower down on the same Young diagram, to create a new partition, or

(2) removing a box from the Young diagram of $\mu_i$, and adding it 
on to the Young diagram of $\mu_{i-1}$
(respectively $\mu_{i+1}$), to create a new partition.
\end{definition}

The posets $\underline{\Lambda}(n,r), \underline{\Lambda}'(n,r)$ 
generalise the poset of partitions of $r$ with $n$ parts or fewer, with the dominance ordering.
 
Recall that the the degree zero part of the
algebra ${\cal S}_{{A}_{\infty}}(n,r)$ is equal to,
$${\cal S}_{V_\infty}(n,r) \cong
\bigoplus_{(r_i)_{i \in \mathbb{Z}}, \sum r_i =r} \left(
\bigotimes_{i \in \mathbb{Z}} {\cal S}(n, r_i) \right).$$

Recall further,
that the module category ${\cal S}(n, r)-mod$ for the classical 
Schur algebra, is a highest weight
category, whose standard modules $\Delta(\lambda)$ are named Weyl modules,
indexed by partitions of $r$ with $n$ parts or fewer.
We write $\nabla(\lambda)$ for the costandard ${\cal S}(n,r)$-module 
corresponding to a partition $\lambda$ of $r$.
 
\begin{definition} Let 
$\underline{\lambda} = (\lambda_i)_ {i \in \mathbb{Z}}$ be a sequence of 
partitions, finitely many of which are non-empty,
whose sizes are given by the sequence $(r_i)_{i \in \mathbb{Z}}$. 
Let $r = \sum_{i \in \mathbb{Z}} r_i$.

The standard module, $\Delta_1(\underline{\lambda})$ 
for ${\cal S}_{{A}_{\infty}}(n,r)$,  
is given by,
$${\cal S}_{A_{\infty}}(n,r)
\bigotimes_{{\cal S}_{V_\infty}(n,r)}
\left(
\bigotimes_{i \in \mathbb{Z}} \Delta(\lambda_i) \right).$$

The costandard module, $\nabla_1(\underline{\lambda})$
for ${\cal S}_{{A}_{\infty}}(n,r)$,  
is the ${\cal S}_{V_\infty}(n,r)$-module,
$$\bigotimes_{i \in \mathbb{Z}} \nabla(\lambda_i).$$
\index{$\Delta_1(\underline{\lambda}), \nabla_1(\underline{\lambda})$}

The costandard module, $\nabla_2(\underline{\lambda})$
for ${\cal S}_{{A}_{\infty}}(n,r)$,  
is given by,
$$Hom_{{\cal S}_{V_\infty}(n,r)}
\left(
{\cal S}_{A_{\infty}}(n,r),
\bigotimes_{i \in \mathbb{Z}} \nabla(\lambda_i) \right).$$

The standard module, $\Delta_2(\underline{\lambda})$
for ${\cal S}_{{A}_{\infty}}(n,r)$,  
is the ${\cal S}_{V_\infty}(n,r)$-module,
$$\bigotimes_{i \in \mathbb{Z}} \Delta_2(\lambda_i).$$
\index{$\Delta_2(\underline{\lambda}), \nabla_2(\underline{\lambda})$}
\end{definition}

\begin{theorem} \label{hw} Let $n \geq r$ be natural numbers.

The module category
${\cal S}_{{A}_{\infty}}(n,r)-mod$, is a highest weight category 
with respect to the poset $\underline{\Lambda}(n,r)$.
Given an element $\underline{\lambda}$ of $\underline{\Lambda}(n,r)$,
the corresponding standard module is
$\Delta_1(\underline{\lambda})$, and the corresponding
costandard module is $\nabla_1(\underline{\lambda})$.

The module category
${\cal S}_{{A}_{\infty}}(n,r)-mod$, is also a highest weight category 
with respect to the poset $\underline{\Lambda}'(n,r)$.
Given an element $\underline{\lambda}$ of $\underline{\Lambda}'(n,r)$,
the corresponding standard module is
$\Delta_2(\underline{\lambda})$, and the corresponding
costandard module is $\nabla_2(\underline{\lambda})$.
\end{theorem}
Proof:

We describe only the quasi-hereditary structure with respect to
$\underline{\Lambda}(n,r)$. The quasi-hereditary structure with respect to
$\underline{\Lambda}'(n,r)$ may be understood similarly.

Throughout this proof, we indiscriminately use the Ringel self-duality
of ${\cal S}_.(n,t)$, for $n \geq t$ (theorem~\ref{donut}). 

We first show that the projective cover of a simple module 
$L(\underline{\lambda})$ is filtered by standard modules 
$\Delta_1(\underline{\mu})$, with 
$\underline{\lambda} \unlhd \underline{\mu}$.
Note that if $P(\underline{\lambda})$ is the projective 
${\cal S}_{V_\infty}(n,r)$ module, with top $L(\underline{\lambda})$, then
${\cal S}_{A_\infty}(n,r) \otimes_{{\cal S}_{V_\infty}(n,r)} 
P(\underline{\lambda})$ is the projective cover of 
$L(\underline{\lambda})$ as a ${\cal S}_{A_\infty}(n,r)$-module.
Since $P(\underline{\lambda})$ is filtered by $\Delta(\underline{\mu})$,
with $\underline{\lambda} \unlhd \underline{\mu}$,
and for $n \geq t$,
the functor $\otimes_{{\cal S}_{V_\infty}(n,t)} \bigvee(n,t)$
is exact on the category of $\Delta$-filtered ${\cal S}_{V_\infty}(n)$-modules,
we deduce that
${\cal S}_{A_\infty}(n) \otimes_{{\cal S}_{V_\infty}(n)} 
P(\underline{\lambda})$ is filtered by $\Delta_1(\underline{\mu})$,
with $\underline{\lambda} \unlhd \underline{\mu}$.

We secondly remark that the simple composition factors of
$\Delta_1(\underline{\lambda})$ are indexed by elements $\underline{\mu}$
of $\underline{\Lambda}(n,r)$, such that 
$\underline{\mu} \unlhd \underline{\lambda}$. This is a consequence
of the branching rule for classical Schur algebras, as 
well as the quasi-heredity of ${\cal S}_{V_\infty}(n)$. 

Thirdly, that the costandard modules 
relevant to this highest weight structure
are the $\nabla_1(\underline{\lambda})$'s is now visible. 
By duality, we need only observe that $\Delta_2(\underline{\lambda})$ is the 
largest quotient of the projective cover of $L(\underline{\lambda})$ 
for whose composition factors $L(\underline{\mu})$ 
(excepting the top $L(\underline{\lambda})$),
the multipartition
$\underline{\mu}$ is strictly smaller than $\underline{\lambda}$,
with respect to the ordering on $\underline{\Lambda}'(n,r)$.
This is apparent from the structure we have already described on
the projective cover of $L(\underline{\lambda})$.
$\Box$

\bigskip

We consider the Schiver double, ${\cal D}_{A_\infty}(n)$. An immediate 
corollary of theorem~\ref{invariant} is,

\begin{theorem} \label{auto}
The action of the infinite dihedral group $D_\infty$ \index{$D_\infty$}
as graph automorphisms
of $A_\infty$ lifts to an action of $D_\infty$ as
algebra automorphisms on
${\cal D}_{A_\infty}(n)$. $\Box$
\end{theorem}

Let 
$$'[1] : \underline{\Lambda}(n,r) \rightarrow \underline{\Lambda}(n,r),$$ 
$$(\underline{\lambda}'[1])_i = \lambda'_{i-1},$$
be the map on $\underline{\Lambda}(n,r)$, which shifts a sequence by $1$, and 
then conjugates each entry in the sequence.

\begin{lemma} \label{last}
Let $n \geq r$ be natural numbers. 

There is an isomorphism,
$\Delta_1(\underline{\lambda}) \cong \nabla_2(\underline{\lambda}'[1])$.
\end{lemma}
Proof:

The natural sequence of homomorphisms,
$$\Delta_1(\underline{\lambda})^{op} \otimes \Delta_1(\underline{\lambda}'[1])
\rightarrow 
\Delta_1(\underline{\lambda})^{op} \otimes_{{\cal S}_{A_\infty}(n,r)} 
\Delta_1(\underline{\lambda}'[1])$$
$$\cong \left(
\bigotimes_{i \in \mathbb{Z}} \Delta(\lambda_i)^{op} \right) 
\otimes_{{\cal S}_{V_\infty}(n,r)}
\left(
\bigotimes_{i \in \mathbb{Z}} \bigvee(n,r_i) \right) 
\otimes_{{\cal S}_{V_\infty}(n,r)}
\left(
\bigotimes_{i \in \mathbb{Z}} \nabla(\lambda_i ') \right)$$
$$\cong \left(
\bigotimes_{i \in \mathbb{Z}} \Delta(\lambda_i)^{op} \right) 
\otimes_{{\cal S}_{V_\infty}(n,r)}
\left(
\bigotimes_{i \in \mathbb{Z}} \Delta(\lambda_i) \right)
\cong k,$$
defines a non-degenerate bilinear form,
$$<,>: \Delta_1(\underline{\lambda})^{op} \times 
\Delta_1(\underline{\lambda}'[1]) \rightarrow k,$$
such that
$<x \circ s, y> = <x, s \circ y>$, for $s \in {\cal S}_{A_\infty}(n,r)$.

Since the dual of the ${\cal S}_{A_\infty}(n,r)$-module 
$\Delta_1(\underline{\lambda})$, is isomorphic to
$\nabla_2(\underline{\lambda})$,
the existence of such a bilinear form seals the proof of the lemma. 
$\Box$

\begin{theorem} \label{exch}
Let $n \geq r$ be natural numbers. 

The Schiver algebra
${\cal S}_{{A}_{\infty}}(n,r)$ is Ringel self-dual.
Indeed, Ringel duality exchanges the two highest 
weight structures we have introduced
on ${\cal S}_{{A}_{\infty}}(n,r)-mod$.

The module category
${\cal D}_{A_{\infty}}(n,r)-mod$ is a highest weight category, 
with respect to the poset $\underline{\Lambda}(n,r)$.
Given an element $\underline{\lambda}$ of $\underline{\Lambda}(n,r)$,
the corresponding standard module is
$\Delta_1(\underline{\lambda})$.

Furthermore, ${\cal D}_{A_{\infty}}(n,r)-mod$ is a highest weight category, 
with respect to the poset $\underline{\Lambda}'(n,r)$.
Given an element $\underline{\lambda}$ of $\underline{\Lambda}'(n,r)$,
the corresponding costandard module is
$\nabla_2(\underline{\lambda})$.

Indeed, ${\cal D}_{{A}_{\infty}}(n,r)$ is Ringel self-dual,
and Ringel duality exchanges these two highest 
weight structures on ${\cal D}_{{A}_{\infty}}(n,r)-mod$.
\end{theorem}
Proof:

As a consequence of lemma~\ref{last}, 
and theorem~\ref{hw}, the regular representation of
${\cal S}_{{A}_{\infty}}(n,r)$ can be filtered by 
$\nabla_2$'s, as well as filtered by $\Delta_2$'s.

Thus, the regular representation is a full tilting module for
${\cal S}_{{A}_{\infty}}(n,r)$ with respect to $\underline{\Lambda}'(n,r)$, 
and indeed ${\cal S}_{{A}_{\infty}}(n,r)$
is Ringel self-dual. Ringel duality exchanges the two highest weight
structures we have defined on ${\cal S}_{{A}_{\infty}}(n,r)$, 
because the functor 
$Hom_{{\cal S}_{{A}_{\infty}}(n,r)}({\cal S}_{{A}_{\infty}}(n,r), -)$ must 
affect costandard modules to become standard modules.
$\Box$

\begin{remark} 
Let ${\cal C}$ be a highest weight category, with poset $\Lambda$, and 
let $\Pi = \Gamma \cap \Omega$ be the 
intersection of an ideal $\Gamma \subset \Lambda$ and a coideal 
$\Omega \in \Lambda$. Then there is a canonically 
defined highest weight category 
${\cal C}(\Pi)$, whose poset is $\Pi$ (see theorem~\ref{lick}).
So long as $\Pi$ is a finite set, 
${\cal C}$ is the module category of a 
quasi-hereditary algebra.
\end{remark}

Let $p$ be a natural number. We define $\Gamma_p(n,r)$ \index{$\Gamma_p(n,r)$}
(respectively $\Gamma'_p(n,r)$) to be the 
ideal (respectively coideal),
of sequences $(\lambda_i) \in \underline{\Lambda}(n,r)$ 
(respectively $\underline{\Lambda}'(n,r)$), all
of whose entries $\lambda_i$ are zero, for $i>p-1$.

Let us define $\Omega_p(n,r)$ \index{$\Omega_p(n,r)$}
(respectively  $\Omega'_p(n,r)$), to be the 
coideal (respectively ideal), of sequences 
$(\lambda_i) \in \underline{\Lambda}(n,r)$
(respectively $\underline{\Lambda}'(n,r)$), all
of whose entries $\lambda_i$ are zero, for $i<0$.

Let $\Pi_p(n,r) = \Gamma_p(n,r) \cap \Omega_p(n,r)$  \index{$\Pi_p(n,r)$}
(respectively $\Pi'_p(n,r) = \Gamma'_p(n,r) \cap \Omega'_p(n,r)$ ),
be the set of sequences 
$(\lambda_i) \in \underline{\Lambda}(n,r)$
(respectively $\underline{\Lambda}'(n,r)$), all
of whose entries $\lambda_i$ are zero, for $i<0$, and $i>p-1$.

\bigskip

For $n \geq r$, let ${\cal Q}_p(n,r)$  \index{${\cal Q}_p(n,r)$}
(respectively ${\cal Q}'_p(n,r)$), be the quasi-hereditary
subquotient of
${\cal D}_{A_\infty}(n,r)$, whose poset
is $\Pi_p(n,r)$ (respectively $\Pi'_p(n,r)$), and whose module category 
is the highest weight category,
$({\cal D}_{A_\infty}(n,r)-mod)(\Pi_p(n,r))$
(respectively $({\cal D}_{A_\infty}(n,r)-mod)(\Pi'_p(n,r))$). 

\begin{remark}
By theorem~\ref{exch}, Ringel duality exchanges the quasi-hereditary algebras,
${\cal Q}_p(n,r)$ and ${\cal Q}'_p(n,r)$.

The quiver $A_\infty$ possesses an orientation-reversing automorphism, which
exchanges vertex $v_0$, and vertex $v_{p-1}$.
By theorem~\ref{auto}, this automorphism lifts to an 
automorphism $\Theta$ 
of ${\cal D}_{A_\infty}(n,r)$. 

The automorphism $\Theta$ provides an isomorphism between
${\cal Q}_p(n,w)$ and ${\cal Q}'_p(n,w)$.
Therefore, ${\cal Q}_p(n,w)$ is Ringel self-dual. 
\end{remark}

We may now formulate a generalisation of
conjecture~\ref{Rock} to $q$-Schur algebras.

Let $q \in k^\times$, and let $p$ be the least natural 
number such that, 
$$1+q+...+q^{p-1} = 0.$$ 
Let $n,w$ be natural numbers, such that $n \geq w$.

\begin{conjecture} \label{SchurRock}
The quasi-hereditary algebra
${\cal Q}_p(n,w)$ 
is Morita equivalent to any Rock block $k{\bf B}_{\rho,w}^{{\cal S}_q}$
of a $q$-Schur algebra, whose weight is $w$.

Indeed,
${\cal Q}_p(n,w)$ is derived equivalent 
to any block $k{\bf B}_{\tau,w}^{{\cal S}_q}$ of 
a $q$-Schur algebra, whose weight is $w$.
\end{conjecture}

How far does this conjecture generalise ?

\begin{question}
Can all blocks of $q$-Schur algebras of weight $w$ be $\mathbb{Z}_+$-graded,
so that the degree zero part is Morita equivalent to the James adjustment algebra ?
\end{question}
 
\begin{remark}
Conjecture~\ref{Rock} and conjecture~\ref{SchurRock} are related
as follows:

For $n \geq w$, the double
${\cal D}_{A_{p-1}}(n,w)$ is obtained from ${\cal Q}_p(n,w)$ by cutting
at the idempotent $j_p \in {\cal Q}_p(n,r)$ corresponding to the subset
$\Pi_{p-1}(n,w) \subset \Pi_p(n,w)$.

We should 
thus define the Specht modules for ${\cal D}_{A_{p-1}}(n,w)$, 
to be those modules $j_p.\Delta$, 
where $\Delta$ is a standard module for ${\cal Q}_p(n,w)$.
\end{remark}

It is now possible to deduce the 
following result from theorem~\ref{nilpotent},
the definition of standard modules for ${\cal D}_{A_\infty}(w,w)$,
and formula~\ref{stolen}, for the
decomposition matrix of $K{\bf B}_{\rho,w}^{{\cal H}_q}$, 
where $q$ is a $p^{th}$ root of unity.

\begin{corollary}
Let $k$ be a field of characteristic $p$.
Then the symmetric group Rock block $k{\bf B}_{\rho,w}^\Sigma$ has
the same decomposition matrix as ${\cal D}_{A_{p-1}}(w,w)$. $\Box$
\end{corollary}

One proof uses the Littlewood-Richardson rule, concerning 
tensor products of modules for the Schur algebra.
Conjecture~\ref{Rock} thus structurally clarifies formula~\ref{stolen}
of Chuang-Tan, and Leclerc-Miyachi.

\begin{center}
{\large Walking along $A_\infty$.}
\end{center}
 
The super-algebra $P_{A_{\infty}}$, is endowed 
with a natural differential $d$ of degree $1$,
given by the infinite sum, 
$\sum_{i \in \mathbb{Z}} e_i$, of all edges.
Indeed, the complex,
$$......... \rightarrow P_{A_{\infty}}  
\rightarrow P_{A_{\infty}} \rightarrow P_{A_{\infty}}
\rightarrow .........,$$
with differential given by right multiplication by $d$,
is a linear exact sequence of
left $P_{A_{\infty}}$-modules.

\bigskip

In the last passage of this letter, we lift this elegant differential
structure on $P_{A_{\infty}}$, to the super-bialgebra
${\cal S}_{A_{\infty}}(n)$, and its double
${\cal D}_{A_{\infty}}(n)$.
We call the resulting chain
complex a ``walk along $A_{\infty}$'', since 
it generalises a homological structure discovered by J.A.
Green on blocks of finite groups of cyclic defect: 
the ``walk around the Brauer tree''.

\bigskip

Let $d$ be the differential on $P_{A_\infty}(n)$ of degree 1,
given by
$$d = (0^{\times V}) \times (1^{\times E}) \in
\left( End_k(k^n) \right) ^{\times V} \times 
\left( End_k(k^{n}) 
\right)^{\times E}.$$
Let $d_r$ be the differential on $P_{A_\infty}(n)^{\otimes r}$
of degree 1, given by
$$d \otimes 1 \otimes... \otimes 1
+ 1 \otimes d \otimes 1 \otimes... \otimes 1
+....+ 1 \otimes ... \otimes 1 \otimes d.$$
Note that $d_r$ is invariant under the action of the symmetric group
$\Sigma_r$, and so
$d_r$ is a differential on the Schiver super-bialgebra 
${\cal S}_{A_\infty}(n,r)$.

\begin{theorem} (``Walk along $A_\infty$'') \label{walk}
The chain complex,
$$......... \rightarrow {\cal S}_{A_\infty}(n,r)  
\rightarrow {\cal S}_{A_\infty}(n,r) \rightarrow 
{\cal S}_{A_\infty}(n,r)
\rightarrow .........,$$
with differential given by right multiplication by $d_r$, is
a linear exact sequence of left ${\cal S}_{A_\infty}(n,r)$-modules.
The chain complex,
$$......... \rightarrow {\cal D}_{A_\infty}(n,r)  
\rightarrow {\cal D}_{A_\infty}(n,r) \rightarrow 
{\cal D}_{A_\infty}(n,r)
\rightarrow .........,$$
with differential given by right multiplication by $d_r$, is
a linear exact sequence of left ${\cal D}_{A_\infty}(n,r)$-modules.
\end{theorem}
Proof:

The classical Koszul complex on $End(k^n)^*$ is the acyclic
chain complex,
$${\cal A}(n) \otimes \bigwedge(n),$$ 
whose differential is given by,
$$d: {\cal A}(n,l) \otimes \bigwedge(n,m)  \rightarrow
{\cal A}(n,l+1) \otimes \bigwedge(n,m-1),$$
$$y_1...y_l \otimes x_1 \wedge ... \wedge x_m \mapsto $$
$$\sum_{i=1} ^m (-1)^{i-1}
y_1...y_l x_i \otimes
x_1 \wedge ... \wedge x_{i-1} \wedge x_{i+1} \wedge ...\wedge x_m.$$ 
Its dual is an acyclic chain complex,
$$ {\cal S}(n) \otimes \bigvee(n).$$

Tensoring together
$\mathbb{Z}$ copies of this dual Koszul complex, and forming
the total complex, in degree $r$
we obtain an acyclic chain complex which corresponds
precisely to the first exact 
sequence of theorem~\ref{walk}.

Tensoring together $\mathbb{Z}$ copies of the Koszul complex, 
along with $\mathbb{Z}$
copies of the dual Koszul complex, in degree $r$
we obtain an exact sequence which corresponds
precisely to the second exact 
sequence of theorem~\ref{walk}.
$\Box$

\begin{remark} \label{Koszul}
The super-algebra, $P_{A_\infty}$, is a Koszul algebra \cite{BGS}.
Its Koszul dual is the path algebra $k A_\infty ^{op}$ on the quiver
$A_\infty$, with opposite orientation.

The Schiver super-algebra ${\cal S}_{A_\infty}(n)$ is not 
a Koszul algebra (unless $k$ is a field of characteristic zero). 
Its degree zero part,
$${\cal S}_{V_\infty}(n) = \bigotimes_{v \in \mathbb{Z}} {\cal S}(n),$$
is not semisimple. However, it does 
possess a linear resolution, and
the algebra,
$$Ext^* _{{\cal S}_{A_\infty}(n)}({\cal S}_{V_\infty}(n),
{\cal S}_{V_\infty}(n)),$$
is isomorphic to ${\cal S}(kA_\infty ^{op})(n)$,
while the path algebra $kA_\infty ^{op}$
is concentrated in parity zero.

The algebra,
$$Ext^* _{{\cal D}_{A_\infty}(n)}({\cal S}_{V_\infty}(n),
{\cal S}_{V_\infty}(n)),$$
is isomorphic to the algebra 
${\cal S}(\Pi_{A_\infty})(n)$, where
$\Pi_{A_\infty}$ is the preprojective algebra on 
the graph $A_\infty$, concentrated in parity zero.

A similar statement is true, relating ${\cal D}_{\tilde{A}_{p-1}}(n)$
and ${\cal S}(\Pi_{\tilde{A}_{p-1}})(n)$.
I prove this, along with various stronger results, 
in my paper ``On seven families of algebras'' \cite{Turner2}.
Sending $p$ to infinity, one obtains theorems for the algebras associated to $A_\infty$.
\end{remark}

\bigskip

\newpage
\printindex
\newpage
\footnotesize{

}

\end{document}